\documentclass[10pt]{article}
\usepackage{amssymb,amsmath,amsfonts,amsthm,enumerate,color}
\usepackage[toc,page,title,titletoc,header]{appendix}
\usepackage{graphicx}
\usepackage{indentfirst}
\usepackage{multicol}
\usepackage{booktabs}
\usepackage{url}
\usepackage{graphicx}
\usepackage{epstopdf}

\numberwithin{equation}{section}

\newtheorem{Theorem}{Theorem}[section]
\newtheorem{Lemma}[Theorem]{Lemma}

\newtheorem{Definition}[Theorem]{Definition}
\newtheorem{Corollary}[Theorem]{Corollary}

\numberwithin{equation}{section}

 \def\p{\partial} 
\def \Vh0{\stackrel{\circ}{V}_h}

\newcommand{\q}{\quad}    
  \def\f{\frac}  

\def\ms{\medskip}  

\def\p{\partial}

\newcommand{\lc}
{\mathrel{\raise2pt\hbox{${\mathop<\limits_{\raise1pt\hbox
{\mbox{$\sim$}}}}$}}}

\newcommand{\gc}
{\mathrel{\raise2pt\hbox{${\mathop>\limits_{\raise1pt\hbox{\mbox{$\sim$}}}}$}}}

\newcommand{\ec}
{\mathrel{\raise2pt\hbox{${\mathop=\limits_{\raise1pt\hbox{\mbox{$\sim$}}}}$}}}

\def\bb{\begin{equation}} \def\ee{\end{equation}}

\def\beqn{\begin{eqnarray}}  \def\eqn{\end{eqnarray}}

\def\beqnx{\begin{eqnarray*}} \def\eqnx{\end{eqnarray*}}

\def\bn{\begin{enumerate}} \def\en{\end{enumerate}}

\def\bd{\begin{description}} \def\ed{\end{description}}

\makeatletter

\newenvironment{figurehere}
  {\def\@captype{figure}}
  {}
\makeatother

\title{Phased and phaseless domain reconstruction \\ 
in inverse scattering problem via scattering coefficients}

\author{
Habib Ammari\thanks{\footnotesize Department of Mathematics,
ETH Z\"urich,
R\"amistrasse 101, CH-8092 Z\"urich, Switzerland. The work of this author was supported  by the
ERC Advanced Grant Project MULTIMOD--267184. (habib.ammari@math.ethz.ch). }
\and Yat Tin Chow\footnote{Department of Mathematics, University of California, Los Angeles, CA 90095-1555, 
USA. (ytchow@math.ucla.edu).}
\and 
Jun Zou\footnote{Department of Mathematics, Chinese University of Hong Kong, Shatin, N.T., Hong Kong. 
The work of this author was substantially supported by Hong Kong RGC grants
(projects 14306814 and 405513). (zou@math.cuhk.edu.hk).
}}

\begin{document}

\date{}
\maketitle

\begin{abstract}
In this work we shall review the (phased) inverse scattering problem and then pursue the phaseless reconstruction from 
far-field data with the help of the concept of scattering coefficients.
We perform sensitivity, resolution and stability analysis of both phased and phaseless problems 
and compare the degree of ill-posedness of the phased and phaseless reconstructions.
The phaseless reconstruction is highly nonlinear and much more severely ill-posed.
Algorithms are provided to solve both the phased and phaseless reconstructions in the linearized case.
Stability is studied by estimating the condition number of the inversion process for both the phased and phaseless cases.
An optimal strategy is suggested to attain the infimum of the condition numbers of the phaseless reconstruction, which may 
provide an important guidance for efficient phaseless measurements in practical applications.
To the best of our knowledge, the stability analysis in terms of condition numbers are new for 
the phased and phaseless inverse scattering problems, and are very important to help us 
understand the degree of ill-posedness of these inverse problems.
Numerical experiments are provided to illustrate the theoretical asymptotic behavior, 
as well as the effectiveness and robustness of the phaseless reconstruction algorithm.
\end{abstract}

\bigskip



\section{Introduction} \label{sec1}

The inverse scattering problems are well known to be severely ill-posed.
It has widespread applications in, e.g., oil/crack detection, target identification,
geophysical prospection, non-destructive testing,
medical imaging, physiological measurement 
\cite{iop1, iop2, iop4, em20, iop6, em16, iop8, iop9, iop10, iop11, kj16, iop13, iop14}.
Due to their applications, inverse scattering problems have attracted much attention, and 
many numerical algorithms have been developed over the recent decades for 
phased reconstruction problems, e.g.,  time-reversal multiple signal classification methods \cite{kj8, kj14}, 
the contrast source inversion methods \cite{kj1, iop1, iop2, iop10, iop11, kj16,kj17}, the continuation method \cite{kj2}, the subspace-based optimization method \cite{kj4, kj5}, the linear sampling or probing methods \cite{kj6, kj11, kj15}, the parallel radial bisection method \cite{kj13}, direct sampling methods \cite{kj7}, multi-level sampling methods 
\cite{kj11, kjself}, etc.

However, in many areas of applied sciences 
it is very difficult and expensive to obtain the phased data of the scattered field, while 
the phaseless data is much easier to acquire.
In addition, 
the phase of the field is more easily polluted by the noise than the
amplitude in many practical situations. 
For instance, the measurement of the phase is extremely difficult when the working frequency is beyond tens of
gigahertz, and one can not expect a good accuracy of the phase measurement \cite{iop6, em16, iop8}.
This motivates the phaseless reconstructions, and attracts huge attention from both 
the physics and mathematics communities.
Nonetheless, the phaseless reconstruction is yet much more severely ill-posed than the phased reconstruction, 
in particular, it appears to be impossible to recover the location of an obstacle only from the modulus of the far-field
pattern owing to the fact that it is invariant under translations \cite{impossible}.
In spite of this drastic difficulty, several approaches have been proposed in literature for the phaseless medium 
reconstruction in optics, acoustic and electromagnetics, e.g., the phaseless
data multiplicative regularized contrast sources inversion method \cite{em20, em19}, 
and several other methods \cite{em15, em17, em18, iop6, em16, iop8, embb}.
Also, the phaseless acoustic (sound-soft) obstacle reconstruction was 
studied in \cite{soundsoft}, where the reconstruction is split into two parts: 
the shape reconstruction from the phaseless data and the location of the obstacle from a few phased measurements.
%Theoretically,  the uniqueness of a phaseless scattering reconstruction in three dimensions 
%was established in \cite{kli2, kli1}, while 
%the phaseless measurements were connected to the Radon transform of the potential under the Born approximation %\cite{born}.
Theoretically,  the uniqueness of a phaseless scattering reconstruction 
was established in \cite{kli2, kli1}, while the phaseless measurements
were connected to the Radon transform of the potential under the Born approximation
\cite{born}, and a new numerical method was proposed in \cite{kli3} for the phaseless
problem using this connection to Radon transform.
There are also other works which address both the theoretical and algorithmic aspects of problems related to phaseless reconstruction of a function or vector, where the phase of a function or vector is recovered 
from the modulus of its evaluation of a special family of functionals 
\cite{Emmanuel, Demanet, Irene}, e.g., the coefficients of a Cauchy wavelet transform.

In this work, we shall study both the phased and phaseless shape reconstructions from 
the far-field data of an acoustic medium scattering problem, so the governing equation of our interest 
is the following Helmholtz system:
\beqn
\Delta\,u+k^2(1+q(x))\,u = 0 \quad \mbox{in} \q \mathbb{R}^2
\label{scattering1}
\eqn
where $u$ is the total field, $q(x) \ge 0 $ is the contrast of the medium and $k$ is the wave number.

Suppose that $D$ is an inclusion contained inside a homogeneous background medium, and 
it is an open bounded connected domain with a $\mathcal{C}^{1,\alpha}$-boundary for $0< \alpha<1$.
We consider the contrast $q$ of the form
\begin{equation}\label{eq:q}
q (x) = \varepsilon^* \chi_{D} (x),
\end{equation}
where $ \chi_{D} $ is the characteristic function of $D$ and $\varepsilon^*>0$ is a constant.
The Helmholtz system (\ref{scattering1}) is often complemented by 
the physical outgoing Sommerfeld radiation condition:
\beqn
   \big| \f{\partial}{\partial |x|} u^s- i k u^s \big|= O(|x|^{-\f{3}{2}})\quad \text{ as }\quad |x| \rightarrow \infty \, ,
    \label{sommerfield}
\eqn
where $u^s:=u-u^i$ is the scattered field and $u^i$ is the incident wave.
Now we can see that the solution $u$ to the system \eqref{scattering1}-\eqref{sommerfield} represents the total field due to the scattering from the inclusion $D$ corresponding to the incident $u^i$.
Then the phased reconstruction is to recover the shape of $D$ from the phased measurements of either the scattered field 
or the far-field, while the phaseless reconstruction is to recover the shape of $D$ 
from only the magnitude of the scattered field or the far-field. 

We shall analyse the sensitivity, resolution and stability of both phased and phaseless reconstructions in the linearized cases under certain measurement strategies, and compare the major differences between these two reconstructions. 
With the help of these analyses, we will propose an efficient measurement method which leads to 
a well-posed inversion process of the phaseless reconstruction.
As demonstrated by our early works \cite{heteroscattering, homoscattering, han}, 
the scattering coefficients provide a powerful and efficient tool for shape classification of a target domain, 
and this concept will also persist in this work to help us establish 
stable reconstruction algorithms and their analysis.

We start by recalling the phased reconstruction in the linearized case so as to provide important insight 
into the highly nonlinear phaseless reconstruction problems. 
Within this framework, we shall provide a resolution analysis on numerical reconstruction with phased data 
in terms of SNR, then propose algorithms for shape reconstructions with the phased measurement.
Another major focus of this work is the stability of the phaseless reconstruction, for which 
we will provide an efficient algorithm, and estimate the condition number of the phaseless inversion process. 
We are able to establish a sharp upper bound for the infimum of the condition numbers of the inversion process 
over all phaseless measurement strategies for a given target resolution, hence propose 
an optimal modulus measurement method.
A similar analysis is carried out for the phased reconstruction to allow a clear comparison between the phased and phaseless reconstructions.
To the best of our knowledge, our stability estimates in terms of condition numbers 
are completely new to inverse medium scattering problems 
and appear to be a very important and effective novel tool 
to help us better understand the degree of ill-posedness and stability 
of both the phased and phaseless reconstructions. 
%\footnote{Prof. Ammari has suggested that "line 4 page 3 should be rewritten."  
%May I know what is a good suggestion to modify that line?}

The remaining part of the work is organized as follows.
In section \ref{sec2}, we review the concept of scattering coefficients and obtain several important results, 
which will be of crucial importance to connect the scattering coefficients to both the phased and phaseless reconstructions, 
and help us develop efficient algorithms and their analysis.
Then we move on to the sensitivity analysis of the phased measurement data in section \ref{sec3}, which will also give 
a link-up between the phaseless data and information about the shape of the domain. 
An important comparison is provided in section \ref{sec3} for the similarities and differences between the phased and phaseless reconstructions.
A phased reconstruction algorithm in the linearized case is then proposed in section \ref{sec4}, 
also with a clear resolution analysis of the algorithm.  This resolution analysis is very helpful 
for us to understand the corresponding resolution constraint in the phaseless reconstruction.
Next, we introduce our phaseless recovery problem in section \ref{sec5}, and provide a phaseless shape reconstruction algorithm in section \ref{sec6}.
A stability analysis is performed for our new phaseless reconstruction algorithm in section \ref{sec_sta}, 
optimal strategies for minimizing the condition number of the inversion process and 
analysis of the differences between the ill-posed natures of the phased and phaseless reconstruction 
are also given. 
Numerical experiments are presented in section \ref{numerical} to confirm the theoretical estimates of the condition number of our inversion process, and illustrate the effectiveness and robustness of our newly-proposed phaseless recovery algorithm.
We emphasize that, although our analyses are performed only for two dimensions, 
similar results and analysis can be extended to higher dimensions as well.

\section{Revisit to the concept of scattering coefficients and its sensitivity analysis} \label{sec2}
In this section, we recall the definition of the scattering coefficient 
\cite{superresolution,heteroscattering,homoscattering} and provide some useful results 
about sensitivity analysis for our subsequent shape reconstruction. 
%Since some of the proofs of these theorems are similar to those in the aforementioned references, only a sketch of the proofs will be provided for these theorems in this section.
%Although some of the following result can be regarded as well-known properties of the scattering amplitude via Lax-Phillips in disguise, we are viewing the problem from another perspective using the concept of the scattering coefficients, which provides us new insight for both the phased and phaseless reconstruction problems.
To do so, we first introduce some useful notation \cite{heteroscattering, homoscattering}.
Let $\Phi_k$  be the fundamental solution to the Helmholtz equation:
\beqn
(\Delta + k^2)\, \Phi_k(x) = \delta_0(x),
\label{solutionfun}
\eqn
where $\delta_0$ is the Dirac mass at $0$, with the outgoing Sommerfeld radiation condition:
$$
    \big| \f{\partial}{\partial |x|}\Phi_k- i  k \Phi_k \big| = O(|x|^{-\f{3}{2}}) \quad \text{ as }\quad |x| \rightarrow \infty \,.
$$
Then $\Phi_k$ can be written in terms of the Hankel function 
$H^{(1)}_0$ of the first kind of order zero:
\beqn
    \Phi_k (x) =
    -\f{i}{4} H^{(1)}_0(k|x|) \,.
    \label{fundamental}
\eqn

Given an incident field $u^i$ satisfying the homogeneous Helmholtz equation:
\beqn
\Delta u^i + k^2 u^i = 0\, ,
\eqn
then the solution $u$ to (\ref{scattering1}) and (\ref{sommerfield})
can be represented by the Lippmann-Schwinger equation as
\beqn
u(x)= u^i(x) - \varepsilon^* k^2 \int_{D}\Phi_k(x-y)u(y)dy \, , \quad x \in \mathbb{R}^2\,,
\label{lipsch}
\eqn
and the scattered field is given by
\beqn
u^s(x)= -  \varepsilon^* k^2 \int_{D}\Phi_k (x-y)u(y)dy \, , \quad x \in \mathbb{R}^2\,.
\label{lipsch2}
\eqn

In what follows, we shall often use the following single-layer potential:
\beqn
 S_{\p D}^k[\phi] (x) = \int_{\p D} \Phi_k( x-y ) \phi (y) \, ds(y) \,, \quad 
 \phi \in L^2(\p D), 
\eqn
then the scattering coefficients are defined as follows \cite{heteroscattering, homoscattering}:
\begin{Definition}
\label{defdefuse}
For $n,m \in \mathbb{Z}$, we define the scattering coefficient $W_{nm} (D, \varepsilon^*,k)$ as 
\begin{equation}
W_{nm} (D, \varepsilon^*,k) =
\int_{\partial \Omega} J_{n}( k r_x ) \, e^{-i n \theta_x} \phi_m (x) \, ds(x) \,,
\label{defint2}
\end{equation}
where $x= r_x (\cos \theta_x, \sin \theta_x)$ is in polar coordinates and the weight function $\phi_m \in L^2 (\partial D)$ is 
such that the pair $( \phi_m, \psi_m  )\in L^2 (\partial D) \times L^2 (\partial D)$ satisfies the following system of integral equations:
\beqn
\begin{cases}
S^{k \sqrt{\varepsilon^* +1}}_{\p D} [\phi_m] (x)- S^k_{\p D} [\psi_m] (x)= J_{m}(k r_x) e^{i m \theta_x} \,, \\
 \f{\p}{\p \nu}S^{k \sqrt{\varepsilon^* +1}}_{\p D}  [\phi_m] (x) \mid_-  - \f{\p}{\p \nu} S^k_{\p D} [\psi_m](x)\mid_+ = \f{\p}{\p \nu} (J_{m}(k r_x) e^{i m \theta_x} ).
\label{defint}
\end{cases}
\eqn

\end{Definition}
\noindent Here $+$ and $-$ in the subscripts indicate respectively the limits from outside $D$ and inside $D$ to $\partial D$ 
along the normal direction, and $\p / \p \nu$ denotes the outward normal derivative.

The scattering coefficients $W_{nm}(D, \varepsilon^*,k) $ are basically the Fourier coefficients of the far-field pattern (a.k.a. the scattering amplitude), 
which is a $2\pi$-periodic function in two dimensions \cite{superresolution, heteroscattering, homoscattering}. 
For the incident field $e^{i k
\widehat{d}\cdot x}$ with a unit vector $\widehat{d}$, we have 
$$
(u - u^i)(x) =  i e^{-\pi i/4} \frac{e^{i k |x|}}{\sqrt{8\pi k |x|}}  A_{\infty} ( \widehat{d}, \widehat{x}, k) + O(|x|^{-\frac{3}{2}}) \quad \mbox{as } |x|\rightarrow \infty,
$$
where $\widehat{x} = x/|x| =(\cos \theta_x, \sin \theta_x)$ and $\widehat{d}= (\cos \theta_d, \sin \theta_d)$ are in polar coordinates, and $A(\theta_d, \theta_x,k):= A_{\infty} (\widehat{d}, \widehat{x} ,k )$ is the so-called far-field pattern.
The following results can be found in 
\cite{superresolution, heteroscattering, homoscattering}.
\begin{Theorem}
\label{farfieldtheorem}
Let $\mathfrak{F}_{\theta_d, \theta_x} [A(\theta_d, \theta_x,k)] (m,n)$ be 
the $(m,n)$-th Fourier coefficient of the far-field pattern $A(\theta_d, \theta_x,k)$ 
with the background wave-number $k$, then it holds that 
\beqn
W_{nm} (D, \varepsilon^* ,k) = i^{(n-m)} \mathfrak{F}_{\theta_d, \theta_x} [A(\theta_d, \theta_x,k)] (-m,n),
\label{fourierfourier}
\eqn
or equivalently, 
\beqn
A(\theta_d, \theta_x,k)  =  \sum_{ m, n \in \mathbb{Z} } i^{(m-n)} e^{- i m \theta_d} e^{i n \theta_x} W_{nm} (D, \varepsilon^* ,k) \,.
\label{fouriersum}
\eqn
\end{Theorem}

%we directly get that the following theorem shall hold following the same arguments as in \cite{superresolution,heteroscattering,homoscattering}.
%\begin{Corollary}
%\label{coll_one}
%Let $u_0 = J_m (k r) e^{i m (\theta - \theta_d)}$, then
%\beqn
%(u- u_0 )^{\text{Far}} = e^{-i m \theta_d} \sum_{n\in\mathbb{Z}}  i^{n-m} W_{nm} (D, \varepsilon^*, k) e^{i n \theta}
%\label{SC_small_ball_purb}
%\eqn
%where the far-field $(u- u_0 )^{\text{Far}} $ is defined as
%\beqn
%(u- u_0 ) :=
%i e^{-\pi i/4} \frac{e^{i k |x|}}{\sqrt{8\pi k |x|}}
%(u- u_0 )^{\text{Far}}  + O(|x|^{-\frac{3}{2}}) \,.
%\eqn
%\end{Corollary}
The following result is a direct consequence of Corollary 7.1 in \cite{heteroscattering}.
\begin{Theorem}
When the contrast $\varepsilon^*$ is small, it holds that 
\beqn
W_{nm} (D,\varepsilon^* ,k)  = \varepsilon^* k^2 \int_{D} J_n(k r ) J_m(k r ) e^{i (n-m) \theta } \, dx + O(  {\varepsilon^*} ^2),
\label{SC_small}
\eqn
which can be simplified for the special case of domain $D$ being the circular shape $D= B_R(0)$:
\beqn
W_{nm} (B, \varepsilon^* ,k) = 2 \pi \varepsilon^* \delta_{nm} k^2 \int_{0}^R [J_n(k r )]^2 r dr  + O(  {\varepsilon^*} ^2).
\label{SC_small_ball}
\eqn
\end{Theorem}
%In particular, we immediately get to the following corollary.
%\begin{Corollary}
%Let $B := B_R(0)$, then
%\beqn
%W_{nm} (B, \varepsilon^* ,k) = 2 \pi \varepsilon^* \delta_{nm} k^2 \int_{0}^R [J_n(k r )]^2 r dr  + O(  {\varepsilon^*} ^2).
%\label{SC_small_ball}
%\eqn
%\end{Corollary}

We remark that the integral appearing in \eqref{SC_small_ball} 
can be calculated explicitly as a Lommel's integral, and this fact will become very helpful in section \ref{sec3}.

Before going to the discussion about 
the phased and phaseless reconstructions, we shall first provide an estimate of the scattering coefficient 
under a perturbation of an open ball $B_R(0)$, which is important for our subsequent analysis
about the resolution of both the phased and phaseless reconstructions.

Let $\nu(x)$ be the outward unit normal to $\p B$, and 
$D := B^\delta$ a $\delta$-perturbation of $B := B_R(0)$ along the variational direction $h \in C^1(\p D)$ 
with $||h||= 1 $:
\beqn
\p B^\delta := \{ \tilde{x} = x + \delta h(x) \nu(x) \, : \, x \in \p B \} \,, \label{eq:bdelta}
\eqn
then we can write  
the difference between the integrals over the domains $B$ and $B^\delta$ for an $L^1$ function $f$:
\beqnx
\int_{B^\delta} f(x) d x - \int_{B} f(x) d x = \delta  \int_{\p B} f(x) h(x) \,ds(x) + O(\delta^2) \,.
\eqnx
Now it follows from this and \eqref{SC_small} that 
\beqn
& & W_{nm} (B^\delta, \varepsilon^* ,k) - W_{nm} (B, \varepsilon^* ,k) \notag \\
&=& \varepsilon^* k^2 \int_{B^\delta \bigcup B \backslash B^\delta \bigcap B} \text{sgn}(h) J_n(k r ) J_m(k r ) e^{i (n-m) \theta } \, dx + O\left( \varepsilon^*\delta^2 + {\varepsilon^*} ^2\right) \notag \\
&=& \varepsilon^* \delta k^2 \int_{\p B} h(x) J_n(k r ) J_m(k r ) e^{i (n-m) \theta } \, dx + O\left( \varepsilon^*\delta^2 + {\varepsilon^*} ^2\right) \notag \\
&=& \varepsilon^* R \delta k^2 J_n(k R ) J_m(k R ) \int_{0}^{2 \pi} h(\theta)  e^{i (n-m) \theta } \, d \theta + O\left( \varepsilon^*\delta^2 + {\varepsilon^*} ^2\right) \notag \\
&  = & 2 \pi R k^2 {\varepsilon^*} \delta J_n(k R ) J_m(k R ) \mathfrak{F}[h](n-m) +
+ O\left( \varepsilon^*\delta^2 + {\varepsilon^*} ^2\right) \,, \notag
\eqn
where $\mathfrak{F}[h](n-m)$ stands for the $(n-m)$-th Fourier coefficient of the perturbation $h$ in the argument $\theta$.

If we further requires that the magnitude of $\delta$ is larger than $\varepsilon^*$ in a way such that $\delta = ({\varepsilon^*})^\alpha$ for some $0 < \alpha < 1$, then we arrive at the following result by some similar argument
to the one in \cite{superresolution}.
\begin{Theorem}
\label{goodtheorem}
Let $D := B^\delta$ be a $\delta$-perturbation of $B := B_R(0)$ as defined in 
\eqref{eq:bdelta}, 
then it holds for $\delta = ({\varepsilon^*})^\alpha$ with $0 < \alpha < 1$ that 
\beqn
W_{nm} (B^\delta, \varepsilon^* ,k) - W_{nm} (B, \varepsilon^* ,k)  = 2 \pi R k^2 ({\varepsilon^*})^{1+\alpha} J_n(k R ) J_m(k R ) \mathfrak{F}[h](n-m) +
O({\varepsilon^*}) ^2.
\label{SC_small_ball_purb}
\eqn
\end{Theorem}

\section{Sensitivity analysis of the phased measurement data}   \label{sec3}
In this section, we shall develop a sensitivity analysis of the phased measurement of the far-field data based on 
the result for the scattering coefficients in Theorem\,\ref{goodtheorem}.
This shall help us provide a crucial expression between the phaseless measurement of the far-field data (i.e., only its magnitude) and the shape $D$.

Suppose that $D := B^\delta$ is a $\delta$-perturbation of $B := B_R(0)$ along the variational direction $h \in C^1(\p D)$ with $||h||= 1 $ as described earlier.
Then it follows from \eqref{fouriersum}, \eqref{SC_small_ball} and \eqref{SC_small_ball_purb} that
\beqn
 &&A_{\infty} (\theta, \widetilde{\theta}, k) \notag \\
 &=&
 \sum_{n,m \in\mathbb{Z} } i^{n-m} e^{i n \theta} e^{ - i m \widetilde{\theta}} \, W_{nm} (B^\delta, \varepsilon^* ,k)
 \notag \\
% &=&
% 2 \pi \varepsilon^* \sum_{l\in\mathbb{Z}} e^{i l (\theta- \widetilde{\theta})} \int_{0}^R [J_l(k r )]^2 r dr + 2 \pi R {\varepsilon^*}^{1+\alpha}
% \sum_{n,m \in\mathbb{Z}} i^{n-m} e^{i n \theta} e^{ - i m \widetilde{\theta}} J_n(k R ) J_m(k R ) \mathfrak{F}(h)(n-m) + O(||\varepsilon^*||^2  ) \notag \\
  &=&
 2 \pi \varepsilon^* k^2 \sum_{l\in\mathbb{Z}} e^{i l (\theta- \widetilde{\theta})} \int_{0}^R [J_l(k r )]^2 r dr
  + 2 \pi R \,({\varepsilon^*})^{1+\alpha} k^2
 \sum_{n, l \in\mathbb{Z}} i^{l} e^{ i l \widetilde{\theta}} e^{i n (\theta- \widetilde{\theta})} J_n(k R ) J_{n-l}(k R ) \mathfrak{F}[h](l) + O(  {\varepsilon^*}) ^2 \notag \,.
%\label{ainfty}
\eqn
Although the above expression looks quite complicated, it can be greatly simplified by some well-known properties of the Bessel functions.
In fact, using the following form of the Graf's addition formula \cite{Watson2}:
\beqn
 \sum_{n = -\infty}^\infty J_n(x) J_{n-l}(y)   e^{i n \theta }
 =
 (-1)^l \left( \frac{x - y\exp(-i \theta) }{ x - y\exp(i \theta)} \right)^{l/2} J_{l} \left( \sqrt{x^2 + y^2 - 2xy \cos(\theta) }\right)
 \,
\eqn
for $x,y >0$ and $x \neq y$, and the well-known property for the second Lommel's integral:
\beqn
 \int_{0}^R [J_l(k r )]^2 r dr = \frac{R^2}{2} [ J_l(k R) ^2 - J_{l-1}(k R) J_{l+1}(k R)  ]
 \,,
\eqn
we can significantly simplify the above expression of the far-field pattern as 
\beqn
 %& &
 A_{\infty} (\theta, \widetilde{\theta}, k) %\notag \\
 &=& \pi R^2 \varepsilon^* \, k^2 \sum_{l\in\mathbb{Z}} e^{i l (\theta- \widetilde{\theta})} [ J_l(k R) ^2 - J_{l-1}(k R) J_{l+1}(k R)  ]  \notag \\
 & & + 2 \pi R \,({\varepsilon^*})^{1+\alpha} \, k^2
 \sum_{n, l \in\mathbb{Z}} i^{l} e^{ i l \widetilde{\theta}} e^{i n (\theta- \widetilde{\theta})} J_n(k R ) J_{n-l}(k R ) \mathfrak{F}[h](l) + O(  {\varepsilon^*}) ^2 \notag \\
&=& \pi R^2 \varepsilon^*k^2 [ J_{0}(2 k R \sin((\widetilde{\theta} - \theta) /2) )  - J_{2}(2 k R \sin((\widetilde{\theta} - \theta) /2) )  ] \notag \\
 & & + 2 \pi R \,({\varepsilon^*})^{1+\alpha} k^2
 \sum_{l \in\mathbb{Z}} (-i)^{l} e^{i l (\widetilde{\theta} + \theta)/2} J_{l}(2 k R \sin((\widetilde{\theta} - \theta) /2) )  \mathfrak{F}[h](l) + O(  {\varepsilon^*}) ^2 \,.
\label{ainfty}
\eqn
An interesting point to note is that the constants $\pi R^2$ and $2 \pi R$ in front of the two terms $\varepsilon^*$ and $({\varepsilon^*})^{1+\alpha}$ are respectively the volume and surface area of the open ball of radius $R$.

Summarizing the above discussions, we come directly to the following theorem.
\begin{Theorem} \label{proj}
If $\delta = ({\varepsilon^*})^\alpha$ for $0 < \alpha < 1$, then
\beqn
A_{\infty} (\theta, \widetilde{\theta}, k) =
 \pi R^2 \varepsilon^* k^2 P_R(\theta, \widetilde{\theta}, k) + 2 \pi R \,({\varepsilon^*})^{1+\alpha} k^2 \langle \mathfrak{F}[h], S_R(\theta, \widetilde{\theta}, k) \rangle_{l^2(\mathbb{C})}
+ O({\varepsilon^*})^2,
\label{A infty_small_ball_purb}
\eqn
where $P_R(\theta, \widetilde{\theta}, k)$ represents the quantity
\beqn
P_R(\theta, \widetilde{\theta}, k) &:=& J_{0}(2 k R \sin((\widetilde{\theta} - \theta) /2) ) - J_{2}(2 k R \sin((\widetilde{\theta} - \theta) /2) )
\label{scalar}
\eqn
and
$S_R(\theta, \widetilde{\theta}, k) \in l^2(\mathbb{C}) $ is a vector given by 
\beqn
S_R(\theta, \widetilde{\theta}, k)_l &:=& i^{l} e^{- i l (\widetilde{\theta} + \theta)/2} J_{l}(2 k R \sin((\widetilde{\theta} - \theta) /2) )  \,.
\label{vector}
\eqn
\end{Theorem}

With the above estimate of the far-field pattern, we can calculate the expression of the magnitude of the far-field pattern,
namely $| A_{\infty} (\theta, \widetilde{\theta}, k) |$, by 
%\beqn
%| A_{\infty} (\theta, \widetilde{\theta}, k) |=
%\left| \pi R^2 \varepsilon^* k^2 P_R(\theta, \widetilde{\theta}, k) + 2 \pi R {\varepsilon^*}^{1+\alpha} k^2 \langle \mathfrak{F}[h], S_R(\theta, \widetilde{\theta}, k) \rangle_{l^2}
%+ O({\varepsilon^*}^2) \right|
%\label{A infty_mag_small_ball_purb}
%\eqn
%handily as follows:
\beqn
\f{| A_{\infty} (\theta_i, \widetilde{\theta}_i, k_i)|^2 - \pi^2 R^4 ({\varepsilon^*})^2 k^4 \left(P_R(\theta_i, \widetilde{\theta}_i, k_i)\right)^2 }{4 \pi^2 R^3 ({\varepsilon^*})^{2+\alpha} k^4 P_R(\theta_i, \widetilde{\theta}_i, k_i) }
&=&  \text{Re} \langle \mathfrak{F}[h], S_R(\theta_i, \widetilde{\theta}_i, k_i)  \rangle_{l^2(\mathbb{C}) }
+ O({\varepsilon^*})^{1-\alpha} \notag\\
&=&  \langle \mathfrak{F}[h], S_R(\theta_i, \widetilde{\theta}_i, k_i)  \rangle_{l^2(\mathbb{R}^2)}
+ O({\varepsilon^*})^{1-\alpha} \,.
\eqn
Due to its great importance for both the subsequent phased and phaseless reconstructions, 
we state it in the following corollary. 
\begin{Corollary}
\label{proj2}
For $\delta = ({\varepsilon^*})^\alpha$ for $0 < \alpha < 1$ it holds that 
\beqn
\langle \mathfrak{F}[h], S_R(\theta_i, \widetilde{\theta}_i, k_i)  \rangle_{l^2(\mathbb{R}^2)} =
\f{| A_{\infty} (\theta_i, \widetilde{\theta}_i, k_i)|^2 - \pi^2 R^4( {\varepsilon^*})^2 k^4 \left(P_R(\theta_i, \widetilde{\theta}_i, k_i)\right)^2 }{4 \pi^2 R^3 ({\varepsilon^*})^{2+\alpha} k^4 P_R(\theta_i, \widetilde{\theta}_i, k_i) }
+ O({\varepsilon^*})^{1-\alpha} \,. \label{eq:new}
\eqn
\end{Corollary}
One interesting observation is that $P_R(\theta, \widetilde{\theta}, k)$ and $S_R(\theta, \widetilde\theta, k) $ 
become very simple for $\theta = \widetilde{\theta}$:
\beqn
P_R(\theta, \theta, k) =1 \,,\quad
S_R(\theta, \widetilde{\theta}, k)_l = \delta_{l 0}   \,.
\eqn
And the expression for the far-field pattern is simplified to
\beqn
A_{\infty} (\theta, \theta, k) =
 \pi R^2 \varepsilon^* k^2  + 2 \pi R \,({\varepsilon^*})^{1+\alpha} k^2 \mathfrak{F}[h](0)
+ O({\varepsilon^*})^2,
\label{A infty_small_ball_purb_2}
\eqn
which illustrates that the direct backscattering data $A_{\infty} (\theta, \theta, k)$ may only provide the information about the area and volume of the inclusions but not the first order perturbation.%, considering the fact that we can always choose $R$ such that $\mathfrak{F}[h](0)= 0$.

%The aforementioned theorems, especially Theorem \ref{goodtheorem}
%and Corollaries \ref{proj} and \ref{proj2} are very useful for the subsequent discussion of the phased and phaseless reconstructions.
%From \cite{Grafnew}, we have the following form of the Graf's addition formula
%\beqn
% \sum_{n = -\infty}^\infty J_n(x) J_{n-l}(y)   e^{i n \theta }
% =
% \left( \frac{x - y\exp(-i \theta) }{ x - y\exp(i \theta)} \right)^{l/2} J_{-l} \left( \sqrt{x^2 + y^2 - 2xy \cos(\theta) }\right)
% \,,
%\eqn
%which, together with the following well-known property for the second Lommel integral,
%\beqn
% \int_{0}^R [J_l(k r )]^2 r dr = \frac{R^2}{2} [ J_l(k R) ^2 - J_{l-1}(k R) J_{l+1}(k R)  ]
% \,,
%\eqn
%%\beqn
%% \int_{0}^R J_l(a r ) J_l(b r ) r dr = \frac{R}{a^2 -b^2} [ b J_l(a R) J_{l-1}(b R) - a J_{l-1}(a R) J_{l}(b R)  ]
%% \,,
%%\eqn
%helps us to significantly simplify the above vector $S(\theta, \widetilde{\theta}, k)_l$ as
%\beqn
%S(\theta, \widetilde{\theta}, k)_l =
%2 \pi R {\varepsilon^*}^{1+\alpha}  (-1)^{l} e^{i l (\widetilde{\theta} - \theta) /2 }   J_{l}(2 k R \sin((\widetilde{\theta} - \theta) /2) )
%+  \pi {\varepsilon^*} R^2  i^{l} e^{i l \theta  }[ J_l(k R) ^2 - J_{l-1}(k R) J_{l+1}(k R)  ] \,.
%\eqn

\smallskip
We end this section with an important remark about some similarities and differences between the phased and phaseless reconstructions in the linearized case. As we see 
from \eqref{A infty_small_ball_purb} that 
%\beqn
%\f{( A_{\infty} (\theta_i, \widetilde{\theta}_i, k_i))^2 - \pi^2 R^4 {\varepsilon^*}^2 k^4 \left(P_R(\theta_i, \widetilde{\theta}_i, k_i)\right)^2 }{4 \pi^2 R^3 {\varepsilon^*}^{2+\alpha} k^4 P_R(\theta_i, \widetilde{\theta}_i, k_i) }
%&=& \langle \mathfrak{F}[h], S_R(\theta_i, \widetilde{\theta}_i, k_i)  \rangle_{l^2(\mathbb{C}) }
%+ O({\varepsilon^*}^{1-\alpha}) \notag\\
%\eqn
%and therefore
\beqn
\langle \mathfrak{F}[h], S_R(\theta, \widetilde{\theta}, k) \rangle_{l^2(\mathbb{C})} =
\f{ A_{\infty} (\theta_i, \widetilde{\theta}_i, k_i) - \pi R^2 \varepsilon^* k^2 P_R(\theta_i, \widetilde{\theta}_i, k_i) }{2 \pi R^2 ({\varepsilon^*})^{1+\alpha} k^2} + O({\varepsilon^*})^{1-\alpha} \,, \label{phased}
%\langle \mathfrak{F}[h], S_R(\theta_i, \widetilde{\theta}_i, k_i)  \rangle_{l^2(\mathbb{C}) } =
%\f{( A_{\infty} (\theta_i, \widetilde{\theta}_i, k_i))^2 - \pi^2 R^4 {\varepsilon^*}^2 k^4 \left(P_R(\theta_i, \widetilde{\theta}_i, k_i)\right)^2 }{4 \pi^2 R^3 {\varepsilon^*}^{2+\alpha} k^4 P_R(\theta_i, \widetilde{\theta}_i, k_i) }
%+ O({\varepsilon^*}^{1-\alpha}) \,,
\eqn
which might be comparable  to Corollary \ref{proj2} above.  However, we do see several differences here.
First, we obtain an approximate value of $\langle \mathfrak{F}[h], S_R(\theta_i, \widetilde{\theta}_i, k_i)  \rangle_{l^2(\mathbb{C}) }$
with the phased measurements in the linearized case, while 
an approximate value of  $\langle \mathfrak{F}[h], S_R(\theta_i, \widetilde{\theta}_i, k_i)  \rangle_{l^2(\mathbb{R}^2)}$
with the phaseless measurements, which is the projection of the original complex inner product to the real part.  Therefore, we can regard the linearized phaseless reconstruction as an "half-dimension" analogy of the linearized phased reconstruction. Second, in the phased reconstruction, the denominator of the right hand side of the equation 
\eqref{phased}
does not involve the division of the term $ P_R(\theta_i, \widetilde{\theta}_i, k_i)$, whereas in the phaseless reconstruction 
the division of the term is involved (cf.\,\eqref{eq:new}). Both differences make the phaseless reconstruction more ill-posed than the phased one.
These differences will be clearly elaborated in section \ref{sec74}.
As the last point, 
it is well-known that the phaseless reconstruction is not unique in a sense that any translation of the inclusion yields the same phaseless measurement. But this is not reflected from the above equation, as we have assumed 
the inclusion is in the center for the sake of exposition.

\section{A phased reconstruction algorithm in the linearized case} \label{sec4}

In this section, we provide a reconstruction algorithm for the phased measurement in the linearized case using the concept of the scattering coefficients, and then a resolution analysis of this algorithm.

\subsection{An algorithm for phased reconstruction}
We recall that $\varepsilon^*$ is the contrast of the inclusion $D$
(cf.\,\eqref{eq:q}) and the perturbation parameter $\delta$ of $D$ is of 
the order $\delta = ({\varepsilon^*})^{\alpha}$ for $ 0 < \alpha < 1$.
Then motivated by the results in Theorems \ref{farfieldtheorem} and \ref{goodtheorem}, 
we come to the following reconstruction algorithm in the linearized case.

\ms
\textbf{Algorithm 1}. Given the measurement $A^{\text{meas}}_{\infty} (\theta, \widetilde{\theta}, k)$.
\begin{enumerate}
\item
Compute  $W_{nm}^{\text{meas}}$ from the Fourier transform as in \eqref{fourierfourier} for $ -N < n,m < N$.
\item
Find $R$, $\varepsilon^*$ from the following minimization problem
\beqn
%\underset{R, \,\varepsilon^*}{{\rm argmin}}
\min_{R, \,\varepsilon^*} 
\sum_{-N < n < N} \left|W^{\text{meas}}_{nn} - \pi R^2 \varepsilon^* k^2 [ J_l(k R) ^2 - J_{l-1}(k R) J_{l+1}(k R)  ]  \right|^2 \,.
\eqn
%Then find $\alpha$.

\item
Compute from \eqref{SC_small_ball_purb} the estimator $(\delta \,\mathfrak{F}[h])^{est}$ of the product of magnitude $\delta$ and Fourier coefficients $\mathfrak{F}[h]$ of the perturbation $h$ for $l \neq 0$:
\beqn
 (\delta \, \mathfrak{F}[h])^{est} := \frac{1}{2N - l} \sum_{m-n = l, \, -N < n,m < N}  \frac{W_{nm}^{\text{meas}} - W_{nm}(B,\varepsilon^*,k) }{ 2 \pi R \varepsilon^* k^2J_n(k R ) J_m(k R ) } \,. \label{recon_sim}
\eqn
\end{enumerate}

We remark that the reconstruction formula \eqref{recon_sim} is similar to the one (5.3) in \cite{superresolution}.
Indeed, considering equations (3.26) in \cite{superresolution}, with any contrast $\varepsilon^*$, the Fourier coefficients of any perturbation $h$ of $B = B_R(0)$ can always be recovered by an inversion of the operator $A(\varepsilon^*)$ as defined in (4.62) in \cite{superresolution} (after a normalization of its wave number $k$ to $k = 1$). However, the coefficients of the matrix $A(\varepsilon^*)$, i.e. $C(\varepsilon^*,n,m)$ defined in (3.27) in \cite{superresolution}, is only given by an expression of resolvent operators, and therefore their explicit expressions are unknown.  The inversion formula (5.3) in \cite{superresolution} is hence inconvenient to be used.
Nonetheless, for a small contrast $\varepsilon^*$, we know now from Theorem \ref{goodtheorem} 
an explicit approximation of coefficients $C(\varepsilon^*,n,m)$ as $C(\varepsilon^*,n,m)\approx  2 \pi R k^2 ({\varepsilon^*})^{1+\alpha} J_n(k R ) J_m(k R ) $.
Therefore \eqref{recon_sim} can be regarded as an easy-to-use approximation of the inversion formula (5.3) of the operator $A(\varepsilon^*)$ described in \cite{superresolution} when the contrast $\varepsilon^*$ is small.

\subsection{Resolution analysis with respect to signal-to-noise ratio}
In this subsection, we perform a resolution analysis of Algorithm 1 in the previous section, which applies also to other reconstruction process derived from \eqref{phased}, since the above algorithm is just a Fourier-transformed version of \eqref{phased}.
Resolution analysis of the above reconstruction with respect to the signal-to-noise ration (SNR) can be conducted following the spirit of the work \cite{GPT}.

In what follows, we assume the following noise model for the far-field pattern:
\begin{equation}
A_{\infty}^{\text{meas}}  (\theta_i, \widetilde{\theta}_j, k) :=  A_{\infty} (\theta_i, \widetilde{\theta}_j, k) +   N (\theta_i, \widetilde{\theta}_j, k)
\end{equation}
where the pairs $\{(\theta_i, \widetilde{\theta}_i)\}_{i,j = 1}^M$ represent the $M$ incident and receiving angles of the measurement evenly distributed on the circle (where $M$ is very large) and $( N (\theta_i, \widetilde{\theta}_j, k) )_{i,j = 1}^M $ is modeled as, for any fixed value of $k$, a complex circular symmetric Gaussian white noise vector with variance:
\begin{equation}
\mathbb{E}[|N (\theta_i, \widetilde{\theta}_j, k) |^2] = \sigma^2 k^4 \,.
\label{noisenoise21}
\end{equation}
Here $\sigma$ represents the noise magnitude and the noise term is assumed to have a variance of quadratic growth with respect to $k$,
as it is direct from \eqref{fouriersum} and \eqref{SC_small} to see that the magnitude of $A_{\infty}^{\text{meas}}  (\theta_i, \widetilde{\theta}_j, k)$ grows in the order of $k^2$ as $k$ grows.

>From the well-known fact that any orthogonal transformation of a Gaussian white random vector will result in another Gaussian white random vector, we arrive at, after taking the discrete Fourier transform in the variables $\theta$ and $\widetilde{\theta}$, that the following model for the scattering coefficient should be in force:
\begin{equation}
W_{nm}^{\text{meas}} (B^\delta,\varepsilon^*,k) = W_{nm} (B^\delta,\varepsilon^*, k) + \hat{N}_{n,m,\varepsilon^*} \, ,
\label{noisenoise}
\end{equation}
where the noise term  $\hat{N}_{n,m,\varepsilon^*}$ is another complex circular symmetric Gaussian random variable such that its variance $\mathbb{E}[|\hat{N}_{n,m,\varepsilon^*}|^2]$ (i.e. the power spectrum of the original random variable $N (\theta_i, \widetilde{\theta}_j, k)$) behaves like 
\begin{equation}
\mathbb{E}[|\hat{N}_{n,m,\varepsilon^*}|^2] = \sigma^2 k^4 \,.
\label{noisenoise2}
\end{equation}
%with $\kappa \in \mathbb{R}$ being a constant. Then $\sigma$ represents the noise magnitude and the noise term is assumed to have a variance of exponential growth as a function of $n,m$. %and linear with respect to $k$.
%This noise model is actually similar to \cite{GPT}, which is a legitimate model to be considered and includes many practical cases.
%Indeed, if $\kappa$ is smaller than zero, then this noise model assumes a variance which has an exponential decay comparable to the magnitudes of $W_{nm}$; see e.g. (3.13) in \cite{heteroscattering} and (2.10) in \cite{homoscattering}. On the other hand, if $\kappa = 0$, this corresponds to a white noise model in the Far-field pattern.
%Worse still, if $\kappa$ is larger than zero, this model actually assumes a variance of an exponentially growth, which shall seriously contaminate the signals from the higher Fourier modes of the perturbation $h$ in the far-field measurements according to Theorem \ref{goodtheorem}.
%We should like to emphasize that our subsequent analysis holds true for all the three cases.
%\footnote{Is this remark good to explain or do we need any modification of the claim to justify the use of this noise model?}

Assume a pair of generic $(k,R)$ such that $kR$ is not a zero of $J_n$ for all $n$.
Then for $l \neq 0$, we obtain from a direct subtraction of \eqref{SC_small_ball_purb} from \eqref{recon_sim}, 
together with \eqref{noisenoise}, that
\beqnx
\mathfrak{F}[h](l)=
(\mathfrak{F}[h])^{est} (l)
+ \frac{\sigma}{ ({\varepsilon^*})^{1+\alpha}} N_{l}  + ({\varepsilon^*})^{1-\alpha} V_l \,,
\eqnx
where $V_l$ represents the approximation error and $N_{l}$ a noise term satisfying the following estimate for its variance for a small $R<1$ and large $N<M$ using \eqref{noisenoise2}:
\beqn
\mathbb{E} [ | N_{l} |^2 ] &=&
 C \frac{1}{(2N - l)^2} \sum_{m-n = l, \, -N < n,m < N}   R^{-2} [J_n(k R ) J_m(k R )]^{-2}  \notag \\
 &\leq& C \frac{1}{(2N - l)^2R^2}\sum_{m-n = l, \, -N < n,m < N}  \frac{  m^m n^n }{ R^{2(m+n)}  }  \notag\\
 &\leq& C \frac{N^{4N}}{R^{2+4N}} \,.
\eqn
The last second inequality above comes from the asymptotic behavior of the Bessel function 
\cite{handbook}:
\beqn
     J_n (t) \bigg/ \f{1}{\sqrt{2 \pi |n|}}\left(\f{e t}{2 |n|}\right)^{|n|} \rightarrow 1
     \label{decayhaha}
\eqn
as $|n| \rightarrow \infty$.
Now assume further that $\varepsilon^* << \sqrt{\sigma}$ and
\beqnx
\text{SNR} := \left( \frac{\varepsilon^*}{\sigma}\right)^2 \,,
\eqnx
then we get
\beqn
\mathbb{E} [ | (\mathfrak{F}[h])^{est} (l)  | ] = \mathfrak{F}[h](l) \,, \quad \mathbb{E} [ | (\mathfrak{F}[h])^{est} (l) -   \mathfrak{F}[h](l)   |^2 ] \leq C \frac{N^{4N}}{R^{2+4N}} \text{(SNR)}^{-(1 + \alpha/2)} \,,
\eqn
which enables us to conclude the following result. 
\begin{Theorem}
\label{resolution}
Suppose that $\delta = ({\varepsilon^*})^\alpha$ for $0 < \alpha < 1$, and $M >> 1$ is the number of measurement points.
If $N<M$ is selected such that
\beqn
C \frac{N^{4N}}{R^{2+4N}} < \text{(SNR)}^{1 + \alpha/2}
\eqn
and that $ \mathfrak{F}[h](l)$ ($|l|\leq N$) are of order $1$, then the $l$-th mode of $h$ can be resolved for $|l|\leq N$, i.e. $\mathbb{E} [ | (\mathfrak{F}[h])^{est} (l) -   \mathfrak{F}[h](l)   |^2  ] < 1$.
\end{Theorem}

%\noindent \textbf{Remark}: As a remark, we may also consider the case of one incidence measurement, which is very similar to what is done in this section.
%For instance, we may propose a reconstruction algorithm with one measurement event given by the incident field $u_0 = J_0(kr)$, i.e. an incident field with is angularly constant.
%From the expression of $(u-u_0)^{\text{Far}}$ given in Corollary \ref{coll_one} and similar ideas as in the previous algorithm, we can following obtain the perturbation $h$ form the following method:
%\begin{enumerate}
%\item
%Find a pair of $(R, \varepsilon^*)$ such that
%\beqn
%\overline{(u-u_0)^{\text{Far}}} = \pi R^2 \varepsilon^* k^2 [ J_0(k R) ^2 + J_{1}(k R)^2  ] \,.
%\eqn
%(Notice the non-uniqueness of the problem. This may help explain the instability of reconstruction from single frequency $k$. Whereas, with multiple frequencies $k$, this process of determining $(R, \varepsilon^*)$ becomes more stable.)
%
%\item
%Let $K = \exp(i k d_\theta \cdot x) = \sum_{n \in \mathbb{Z}} i^n J_n(k R ) e^{i n \theta}$, where the last equality comes from the well-known Jacobi-Anger Identity.
%Then from \eqref{SC_small_ball_purb}, we can obtain an estimator of $(\delta h)^{est}$, where $h$ is the perturbation and $\delta$ is the magnitude, as follows:
%\beqn
%(\delta h) * K = 2 \pi R {\varepsilon^*}^{1+\alpha}  k^2 J_0(k R )  (u-u_0)^{\text{Far}}
%\eqn
%A stable deconvolution technique, e.g. Richardson-Lucy Deconvolution, shall solve for this estimator in a robust way.
%\end{enumerate}

\section{Introduction to phaseless reconstruction}  \label{sec5}

Phaseless reconstruction originates from the physical background that we can usually only measure the magnitude of some data, for example, the magnitude of the far-field pattern.
As briefly explained in section \ref{sec1}, it is quite difficult and expensive to obtain the phased data
in many physical and engineering applications, and the phase of a measurement is easily 
contaminated by noise. On the other hand, 
the phaseless data is much easier to obtain and less contaminated in many practical situations.  Due to 
these facts, the phaseless reconstruction has attracted wide attention.

\subsection{Brief history of a general phaseless reconstruction problem}\label{sec:history}
Let us first give a brief introduction and history of a general phaseless reconstruction. As in \cite{Candes}, for a given set of $m$ sampling vectors, $\textbf{z}_1, \cdots, \textbf{z}_m$, we intend to recover a vector $\textbf{x}$ from some phaseless data. This may be formulated as 
\beqn
\text{Find } \textbf{x} \text{ such that } A(\textbf{x}) = b
\label{eq:vector}
\eqn
where $A: \mathbb{C}^N \rightarrow R^m $ is given by $A(\textbf{x})_i = | \langle \textbf{x},\textbf{z}_i\rangle|^2$.
One may consider a convexification of the problem \eqref{eq:vector} \cite{Candes}:
\beqn
\text{Find } \textbf{X} \geq 0 \text{ such that } \mathcal{A}(\textbf{X}) =
b\, ,
\eqn
where $\mathcal{A}: \mathcal{H}^{N \times N} \rightarrow R^m $ is given by 
$\mathcal{A}(\textbf{X})=\textbf{z}_i^* \textbf{X} \textbf{z}_i $, which helps reduce the complexity of solving the problem, 
as well as provide uniqueness results under some practical conditions.
For instance, this problem is proven to have a high probability that it is uniquely solvable up to a unit complex number stably from $O(N \log N )$ random measurements \cite{Demanet}.
We remark that a stablized version of convexification is given by
\begin{equation}
\text{Find } \textbf{X} \geq 0 \text{ such that } || \mathcal{A}(\textbf{X}) - b || \leq \epsilon ||X_0||_2 \,.
\end{equation}

Another more general form of phaseless reconstruction (which generalizes the above) comes from recovering the phase of a function/vector from the modulus of its evaluation by a family of functionals.
In a more precise way, let $E$ be a complex vector space and
$\{L_i\}_{i \in I}$ be a family of functionals. Then this phaseless reconstruction reads:
\beqn
\text{Find } \textbf{f} \in E \text{ such that } |L_i(\textbf{f})| = b\,.
\eqn
In the case where $\{L_i\}_{i \in I}$ represents the wavelet transform by the Cauchy wavelets, it was shown in \cite{Irene} 
that the modulus of the wavelet
transform uniquely determines the function up to a global phase, and the reconstruction operator is 
continuous but not uniformly continuous.

\subsection{Introduction to our phaseless reconstruction problem}
The convexification discussed in section\,\ref{sec:history}  
is a very interesting approach, but the purpose,  framework and analysis 
of our phaseless reconstruction here are very different.
We aim to achieve numerical reconstructions of inhomogeneous domains in the linearized case.
We will provide an algorithm for the domain reconstruction from some phaseless far-field data, 
and estimate the condition number of this reconstruction process, and establish an upper bound of 
its infimum over all phaseless measurement strategies. This casts light on how we can obtain an optimal strategy to perform effective phaseless measurements such that the phaseless inversion process shall be well-posed.
For comparison purpose, a similar analysis technique is also performed on its phased counterpart, and 
a comparison between the phased and phaseless reconstructions shall be made.

\section{Phaseless domain reconstruction algorithm in linearized cases}  \label{sec6}

In this section, we provide a new method for the domain reconstruction from the phaseless far-field data based on 
our analyses and results in sections \ref{sec2} and \ref{sec3}. We first recall from Theorem \ref{proj2} 
the relation  
\beqn
\f{| A_{\infty} (\theta_i, \widetilde{\theta}_i, k_i)|^2 - \pi^2 R^4 ({\varepsilon^*})^2 k^4 \left(P_R(\theta_i, \widetilde{\theta}_i, k_i)\right)^2 }{4 \pi^2 R^3 ({\varepsilon^*})^{2+\alpha} k^4 P_R(\theta_i, \widetilde{\theta}_i, k_i) }
=  \langle \mathfrak{F}[h], S_R(\theta_i, \widetilde{\theta}_i, k_i)  \rangle_{l^2(\mathbb{R}^2)}
+ O({\varepsilon^*})^{1-\alpha}\,,
\eqn
where $P_R(\theta_i, \widetilde{\theta}_i, k_i) \in \mathbb{R}$ and $S_R(\theta_i, \widetilde{\theta}_i, k_i) \in l^2{(\mathbb{C})}$ are given in \eqref{scalar} and \eqref{vector}.
Therefore, from a finite number of $M$ measurements $| A_{\infty} (\theta_i, \widetilde{\theta}_i, k_i) |$ 
($ 1 \leq i \leq M$),
%\beqnx
%| A_{\infty} (\theta_i, \widetilde{\theta}_i, k_i)| =
% \left| \pi R^2 \varepsilon^* k^2 P_R(\theta_i, \widetilde{\theta}_i, k_i)  + 2 \pi R {\varepsilon^*}^{1+\alpha} k^2 \langle \mathfrak{F}[h], S_R(\theta_i, \widetilde{\theta}_i, k_i)  \rangle_{l^2{(\mathbb{C})}}
%+ O({\varepsilon^*}^2) \right|
%\eqnx
we obtain the following linear approximation of $\langle \mathfrak{F}[h], S_R(\theta_i, \widetilde{\theta}_i, k_i)\rangle$
as the measurement quantities from the phaseless measurements:
\beqn
\langle \mathfrak{F}[h], S_R(\theta_i, \widetilde{\theta}_i, k_i)  \rangle_{l^2(\mathbb{R}^2)} \approx
\f{| A_{\infty} (\theta_i, \widetilde{\theta}_i, k_i)|^2 - \pi^2 R^4( {\varepsilon^*})^2 k^4 | P_R(\theta_i, \widetilde{\theta}_i, k_i)|^2 }{4 \pi^2 R^3 ({\varepsilon^*})^{2+\alpha} k^4 P_R(\theta_i, \widetilde{\theta}_i, k_i) } \,.
\eqn
This is of crucial importance for us to derive an algorithm for the domain reconstruction from the phaseless far-field measurements.

\subsection{Phaseless reconstruction algorithm}  \label{phaseless_algo}

We are now ready to introduce our phaseless reconstruction algorithm.
Following Theorem \ref{resolution} from the resolution analysis for the phased reconstruction in section \ref{sec3}, we can directly infer that the resolution with respect to SNR in the phaseless reconstruction should not surpass the $N$-th Fourier mode, where $N$ satisfies the inequality $C \frac{N^{4N}}{R^{2+4N}} < \text{(SNR)}^{1 + \alpha/2}$ for some $C$ and $\alpha$.
%Therefore the inversion of the Fourier mode over the $N$-th mode 
%is unstable numerically.
Hence in our reconstruction algorithm,
we may always assume that $\mathfrak{F}[h](l) = 0$ for $|l|> N$ for some $N$ and consider only the inversion of finite dimensional operators, and the contribution of $\mathfrak{F}[h](l)$ for $|l|> N$ to the measurement data can be regarded as noise.
Now since $h(\theta) \in \mathbb{R}$ for all $\theta$, we have the following additional constraints on the Fourier coefficients:
\beqn
\mathfrak{F}[h]( - l) = \overline{ \mathfrak{F}[h]( l) } \,.
\eqn
This set of constraints is very important in our subsequent analysis.
We assume again that the magnitude of the perturbation $\delta$ is of the form 
$\delta = ({\varepsilon^*})^{\alpha}$ for $0 < \alpha < 1$, where $\varepsilon^*$ is the contrast of the inclusion.
>From Theorem \ref{proj2}, we can now suggest the following phaseless reconstruction algorithm.

\ms 
\textbf{Algorithm 2}. 
Given a positive integer $N$ and $M$ measurements of the magnitude $| A_{\infty}^{\text{meas}} (\theta_i, \widetilde{\theta}_i, k_i) |$ ($ 1 \leq i \leq M$) of the far-field. 
\begin{enumerate}
\item
Find the pair $(R , \varepsilon^*)$ that minimizes the following functional: 
\beqn
%\text{argmin} 
\sum_{1 \leq i \leq M } \left| |A_{\infty}^{\text{meas}} (\theta_i, \widetilde{\theta}_i, k_i)|^2 - \pi^2 R^4 k^4 {\varepsilon^*}^2 \left(P_R(\theta_i, \widetilde{\theta}_i, k_i)\right)^2 \right|^2 \,,
\eqn
where the values $P_R(\theta_i, \widetilde{\theta}_i, k_i)$ are computed from \eqref{scalar}.

\item
%From each of the values $| A_{\infty}^{\text{meas}} (\theta_i, \widetilde{\theta}_i, k_i)|$, calculate $| A_{\infty}^{\text{meas}} (\theta_i, \widetilde{\theta}_i, k_i)|^2 - \pi^2 R^4 {\varepsilon^*}^2 k^4 | P_R(\theta_i, \widetilde{\theta}_i, k_i)|^2$; then divides each value by $4 \pi^2 R^3 {\varepsilon^*}^{2} k^4P_R(\theta_i, \widetilde{\theta}_i, k_i)$ to obtain 
%
Compute the following quantities for $ 1 \leq i \leq M$:
\beqn
\f{| A_{\infty}^{\text{meas}} (\theta_i, \widetilde{\theta}_i, k_i)|^2 - \pi^2 R^4 {\varepsilon^*}^2 k^4 | P_R(\theta_i, \widetilde{\theta}_i, k_i)|^2 }{4 \pi^2 R^3 {\varepsilon^*}^{2} k^4P_R(\theta_i, \widetilde{\theta}_i, k_i) } \,.
\label{step2}
\eqn

\item
Calculate the estimator $(\delta \,\mathfrak{F}[h])^{est}(l)$ of the product of magnitude $\delta$ and Fourier coefficient 
$\mathfrak{F}[h]$ of the perturbation $h$ for $|l| \leq N$ by the inversion of the following system of linear equations:
\beqn
\langle (\delta \,\mathfrak{F}[h])^{est} , S_R(\theta_i, \widetilde{\theta}_i, k_i)  \rangle_{l^2(\mathbb{R}^2)} =
\f{| A_{\infty}^{\text{meas}} (\theta_i, \widetilde{\theta}_i, k_i)|^2 - \pi^2 R^4 {\varepsilon^*}^2 k^4 | P_R(\theta_i, \widetilde{\theta}_i, k_i)|^2 }{4 \pi^2 R^3 {\varepsilon^*}^{2} k^4P_R(\theta_i, \widetilde{\theta}_i, k_i) } \,,
\label{step3}
\eqn
%\item
%Recover $\mathfrak{F}(h)(l)$ from $X$ as $X = \mathfrak{F}(h) \mathfrak{F}(h)^*$
under the constraints
\beqn
(\mathfrak{F}[h])^{est}( - l) = \overline{ (\mathfrak{F}[h])^{est}( l) } \,.
\eqn
%by direct inverting linear operator $\langle \, \cdot \, , S_R(\theta_i, \widetilde{\theta}_i, k_i)  \rangle_{l^2(\mathbb{R}^2)}$ in the subspace satisfying the set of constraints.
%with a similar algorithm as
%\beqn
%\text{Find } \textbf{X} \text{ such that } || \mathcal{A}(\textbf{X}) - b || \leq \epsilon ||X_0||_2 \,.
%\eqn
%with $\mathcal{A}: \mathcal{H}^{M \times M} \rightarrow R^N $ with $\textbf{X} \rightarrow S(\theta_i, \widetilde{\theta}_i, k_i) ^* \textbf{X} S(\theta_i, \widetilde{\theta}_i, k_i) $, i.e. choose $\textbf{z}_i = (\theta_i, \widetilde{\theta}_i, k_i)$.
\end{enumerate}

We should be reminded that the above algorithm can only provide stable inversion and reasonable resolution of the perturbation $h$ up to at most the $N$-th Fourier mode, 
where $N$ shall satisfy the inequality $C \frac{N^{4N}}{R^{2+4N}} < \text{(SNR)}^{1 + \alpha/2}$ for some $C$ and $\alpha$, as shown in Theorem \ref{resolution} in section \ref{sec3}.

\section{Stability of the phaseless domain reconstruction} \label{sec_sta}

We are now ready to discuss the stability of the phaseless reconstruction by estimating the condition number of this inversion process.
Before going into detailed estimates, we shall state our inversion problem in a more concise manner, which will provide a clear framework for our subsequent analysis.
For this purpose, we first define three operators for a given pair $(N,M) \in \mathbb{N}$, where
two of them are linear in nature while the other one is nonlinear:
\begin{enumerate}
\item
the component-wise squaring of a vector followed by a subtraction of another known vector, i.e. the action
$v_i \mapsto v_i^2 - \pi^2 R^4 {\varepsilon^*}^2 k^4 | P_R(\theta_i, \widetilde{\theta}_i, k_i)|^2$,
which appears in Step $2$ of Algorithm $2$;
and we write this nonlinear operator as $F:\mathbb{R}^M \rightarrow \mathbb{R}^M$;
\item
the component-wise multiplication of a vector $v_i \mapsto 4 \pi^2 R^3 {\varepsilon^*}^{2} k^4P_R(\theta_i, \widetilde{\theta}_i, k_i)\, v_i$, which appears in Step $2$ of Algorithm $2$;
and we will write this linear operator as $L : \mathbb{R}^{M} \rightarrow \mathbb{R}^{M}$;
\item
the linear operator $ v \mapsto \left( \langle \, v \, , S_R(\theta_i, \widetilde{\theta}_i, k_i) \rangle_{l^2(\mathbb{R}^2)} \right)_{i = 1}^M$, which appears in Step $3$ of Algorithm $2$.
We shall write this linear operator as $T : \mathbb{C}^{2N} \oplus \{ 0\} \cong \mathbb{R}^{4N}  \rightarrow \mathbb{R}^{M}$.
\end{enumerate}
Without loss of generality, we may always choose a radius $R$ such that the zeroth Fourier coefficient $ \mathfrak{F}[h](0) $ is zero.
With these preparations, we can write \eqref{step3} as
\beqn
(L \circ T)  \left[ (\delta \,\mathfrak{F}[h])^{est} \right] = F(| A_{\infty}^{\text{meas}} (\theta_i, \widetilde{\theta}_i, k_i)| ) \,.
\label{operatoreq}
\eqn
Then our phaseless inversion problem can be precisely stated as follows:
given a value of $SNR$, with a number $N$ such that
$
C {N^{4N}}/{R^{2+4N}} < \text{(SNR)}^{1 + \alpha/2}
$
for some $C$ and $\alpha$, we recover the Fourier coefficients $ \delta \left( \mathfrak{F}[h]^{est} (l) \right)_{l = -N}^N \in \mathbb{C}^{2N} \oplus \{ 0\} \cong \mathbb{R}^{4N} $ from \eqref{operatoreq} with $M$ measurements:
 \beqnx
 (b_i)_{i = 1}^M  := F(| A_{\infty}^{\text{meas}} (\theta_i, \widetilde{\theta}_i, k_i)|) = \left(| A_{\infty}^{\text{meas}} (\theta_i, \widetilde{\theta}_i, k_i)|^2 - \pi^2 R^4 {\varepsilon^*}^2 k^4 | P_R(\theta_i, \widetilde{\theta}_i, k_i)|^2  \right)_{i = 1}^M   \in \mathbb{R}^{M} \,
\eqnx
subjected to the following extra set of constraints in $ \mathbb{R}^{4N} $ as
\beqnx
\text{Re}(\mathfrak{F}[h])^{est}( - l) - \text{Re} (\mathfrak{F}[h])^{est}( l) = 0 \,, \qquad  \text{Im}(\mathfrak{F}[h])^{est}( - l) + \text{Im} (\mathfrak{F}[h])^{est}( l) = 0\,.
\eqnx
>From now on, we shall denote this set of linear constraints as
\beqn
 C \left[ (\delta \,\mathfrak{F}[h])^{est} \right] = 0 \,. \label{constraint}
\eqn
After this restatement of the phaseless reconstruction problem, we can directly infer that the stability of the inversion lies in the stability of inversion of the linear operators $L$ and $T$ in the subspace $\text{ker} (C)$  under a certain noise level.
Therefore, the aim of this section is to estimate the condition numbers of the operators $T$ and $L$ in this subspace.
To the best of our knowledge, the stability estimates on condition numbers are novel for inverse problems, and are very important for us to understand the degree of ill-posedness and stability of the reconstruction problem, as well as 
to provide optimal methods to minimize these two condition numbers by making wise measurements or regularizations.

\subsection{Estimation of the condition number of $T$} \label{condT}
We come now to the estimate of the condition number of operator $T$. 
For notational sake, we first introduce two more operators,
$\iota_0 : \mathbb{C} \rightarrow \mathbb{R}^{2} \,, z \mapsto (\text{Re}(z) ,  \text{Im}(z))$
and their liftings on the linear operators over the corresponding spaces
$\iota : \mathfrak{L}(\mathbb{C}) \cong  \mathbb{C} \rightarrow \mathfrak{L}(\mathbb{R}^{2}) \cong M_{2 \times 2} \,, z \mapsto
\begin{pmatrix}
\text{Re}(z)  & \text{Im}(z) \\
-\text{Im}(z) & \text{Re}(z)
\end{pmatrix}
$. And we shall also use the projection map
$\pi_{\text{Re}} : M_{2 \times 2} \rightarrow M_{1 \times 2}\,, \begin{pmatrix}
a & b \\
c & d
\end{pmatrix} \mapsto \begin{pmatrix}
a & b
\end{pmatrix} $. It is easy to check
\beqn
 [ \pi_{\text{Re}} \circ \iota (\bar{z}) ] \left( \iota_0(w) \right)=  \langle \iota_0(z), \iota_0(w) \rangle_{\mathbb{R}^{2}} = \text{Re}(\bar{z} w)\,. \label{special}
\eqn

Before we go on to study the stability of the reconstruction problem, we shall also provide a clear concept and define the condition number of $T$ subjected to the constraint $Cx = 0$, denoted as $\kappa(T,\text{ker}(C) )$, where 
$C :\mathbb{R}^P \rightarrow \mathbb{R}^Q$ for $Q \leq P$ is another linear operator.
First, for the sake of exposition, we shall denote $C^{\perp}$ as the set of all matrices $E$ such that its column vectors are linearly independent and span the orthogonal complement of the row space of $C = (C_1, C_2, \dots, C_n)^T$,  i.e.
\beqnx
C^{\perp} := \{
(E_1, E_2, \dots, E_n): \langle E_1, ... , E_n, C_1, ... , C_n \rangle = \mathbb{R}^P ,  \langle E_1, ..., E_{j-1} , E_{j+1},..., E_n, C_1,... , C_n \rangle \neq \mathbb{R}^P ~~\forall j \}  \,.
\eqnx
Now, if we solve the following constraint problem for a given triple $(T, b, C)$:
\beqnx
\text{Find $x \in \mathbb{R}^P$ such that } \quad  Tx = b \quad \text{and} \quad C x = 0 \,,
\eqnx
or its least-squares formulation:
\beqnx
\min_{ \underset{s.t.\,Cx=0}{x\in \mathbb{R}^P}} || Tx - b ||_{2}^2  \,,
\eqnx
%
%\beqnx
%\text{Find $x\in \mathbb{R}^P$ which minimizes } ~~\min_{x\in \mathbb{R}^P s.t. \,C x = 0} || Tx - b ||_{2}^2 \quad \text{such that} \quad C x = 0 \,,
%\eqnx
we are actually parameterizing the kernel $\text{ker}(C)$ by an orthogonal complement of the row space of $C$, and then solve the equation $T x = b$ under this parametrization (either in the strict sense, or the least-squares sense), i.e., 
solve for $y$ the equation (with $E \in C^{\perp}$):
\beqnx
(T \circ E ) \, y = b \,,
\eqnx
or the least-squares minimization (with $E \in C^{\perp}$):
\beqnx
\min_{ y \in \mathbb{R}^Q}   || (T \circ E )\,  y - b ||_{2}^2\,.
\eqnx
%
%\beqnx
%\text{Find $y \in \mathbb{R}^Q$ which minimizes} ~~|| (T \circ E )\,  y - b ||_{2}^2 \quad \text{ for } E \in C^{\perp} \,.
%\eqnx
>From this definition, one can easily get that the operator $T$ is invertible with its solution in the subspace $C x = 0$ if and only if $T \circ E$ is invertible. One can also directly get that if $T x = b$ and $T \tilde{x} = \tilde{b}$ where $x , \tilde{x} \in \text{ker} C$, then we have the following estimate
\beqnx
\frac{|| y - \tilde{y} || }{ || y ||} \leq \kappa(T \circ E ) \frac{|| b - \tilde{b} || }{ || b ||} \,,
\eqnx
for any $E \in C^{\perp}$, where $y, \tilde{y}$ are defined such that $E y = x , E \tilde{y} = \tilde{x}$.
Hence, in order to study the stability of the inversion process of $T$ in the subspace, we are motivated to define the condition number of $T$ under the constraint $Cx = 0$ as
\beqn
\kappa(T,\text{ker}(C) ) : = \inf \{ \kappa(T \circ E ) \,:\, E \in C^{\perp} \} \,.
\eqn

\subsubsection{A measurement strategy for phaseless reconstruction} \label{sec71}
%In the previous subsection, we demonstrates that with random measurement events, we can have a quite promising well-conditioning of the inversion process of $T$.

In this subsection, we proceed to develop a good measurement strategy which can minimize the condition number of $T$ and ensure the well-posedness of the inversion concerned.

Indeed, we shall intuitively expect to have a good strategy in choosing the measurement set $(\theta, \widetilde{\theta}, k)$ by gazing at the vector $S_R$ in \eqref{vector}: for a given target resolution $N$, one may choose $2 k_i R \sin((\widetilde{\theta}_i - \theta_i) /2)$ such that they attain the $m_0$-th local extremum of $J_l$ for $1<l<N$, i.e., the values of $b_{l,m_0}$ where $J_l' (b_{l,m_0}) = 0$ for a given $m_0$.
With this particular choice of $M$ sets of measurement data $\{ (\theta_i, \widetilde{\theta_i}, k_i) \}_{i=1}^M$ such that $2 k_i R \sin((\widetilde{\theta}_i - \theta_i) /2) \in \{b_{1,m_0}, b_{2,m_0}, ... , b_{N,m_0}\}$, the operator $T$ is expected to be well-conditioned and therefore provide a good set of information for the geometry of the inclusion. This shall be indeed verified in this subsection.

In what follows, we aim to estimate, for a given resolution $N$, the infimum over the condition numbers of all the operators $T:\mathbb{R}^{4N} \rightarrow \mathbb{R}^M$ subjected to the constraint $C x = 0$, i.e.,
\beqnx
\kappa_{\text{inf},N} := \inf \bigg\{ \kappa(T,\text{ker}(C) ) \,:\, \{ (\theta_i, \widetilde{\theta_i}, k_i) \}_{i=1}^M \in [0.2\pi]^2 \times (0,\infty) \, , \, M \in \mathbb{Z} \bigg\}
\label{minumum}
\eqnx
by appropriately choosing the vectors $\{(\theta_i, \widetilde{\theta_i}, k_i)\}_{i=1}^M$.
Indeed, from the following well-known asymptotic of $J_l$ \cite{handbook} for all $l$:
\beqn
J_l \left( z  \right) = \sqrt{\frac{2}{ \pi z}} \cos \left( z - \frac{2 l + 1}{ 4} \pi \right) + O( {z}^{-3/2}) \,,
\label{order0}
\eqn
we directly have for a fixed $l$ that 
\beqn
b_{l,m_0} \bigg\slash  \frac{(4 m_0 + 2l + 1)\pi}{4}  \rightarrow 1
\label{asymroot}
\eqn
as $m_0$ goes to infinity.
Therefore, for a given large $m_0$, if we choose $(\theta_i, \widetilde{\theta_i}, k_i)$ as the form $(\theta_i, \theta_i + \pi,  \frac{4 m_0 + 2 J_i + 1}{8 R})$, where $(J_i)_{i=1}^M \in \mathbb{Z}$ are some integral indices to be specified later, then we have directly from \eqref{order0} and \eqref{asymroot} that 
\beqn
S_R(\theta_i, \widetilde{\theta_i}, k_i)_l = e^{i l \theta_i}
\sqrt{\frac{16 R}{ \pi (4 m_0 + 2 J_i + 1) }} \cos \left( \left( m_0 + \frac{ J_i - l  }{2} \right)  \pi \right)
+ O( m_0^{-3/2})\,.
\eqn
Let $T_{m_0}$ be the linear operator $T$ with this specific arrangement of measurements for a given $m_0$.
If we further denote $L := T_{m_0}^T T_{m_0}$ in the form of a block matrix $\left( L_{lm} \right)_{-N\leq l,m\leq N,\, l,m \neq 0}$, then from \eqref{special}, each of the block $L_{lm}$ will be the following $2 \times 2$ matrix:
\beqn
L_{lm}
&=& \iota_0 \left[S_R(\theta_i, \widetilde{\theta_i}, k_i)_l\right]^T  \iota_0\left[ S_R(\theta_i, \widetilde{\theta_i}, k_i)_m \right]  \notag \\
&=& \frac{16 R}{ \pi} \, \iota \left[
\sum_{i=1}^M
L_{l,m,\theta_i}
\frac{1}{ 4 m_0 + 2 J_i + 1 } \cos \left( \left( m_0 + \frac{ J_i - m  }{2} \right) \pi \right) \cos \left( \left( m_0 + \frac{ J_i - l  }{2} \right)  \pi \right)
\right] + O( m_0^{-2}) \,,  \notag
\eqn
where the matrix $L_{l,m,\theta}$  has the form \beqnx
L_{l,m,\theta} = \begin{pmatrix}
\cos\left(m   \theta  \right) \cos\left(l  \theta \right)  & \cos\left(m  \theta \right) \sin\left(l  \theta \right)\\
\cos\left(m  \theta \right) \sin\left(l  \theta \right) & \sin\left(m   \theta\right) \sin\left(l  \theta \right)  \end{pmatrix} \,.
\eqnx
For the sake of exposition, we further denote $\theta_i =  2 \pi I_i   / N $ where $(I_i)_{i=1}^M \in \mathbb{Z}$ are
some indices to be chosen later.
%, then
%$L_{l,m,i} = \begin{pmatrix}
%\cos^2\left((l-m) ( \frac{ 2 \pi I_k}{N} + \frac{\pi}{2})\right)  & \cos\left((l-m) ( \frac{ 2 \pi I_k}{N} + \frac{\pi}{2}) \right) \sin\left((l-m) ( \theta_i + \frac{\pi}{2}) \right)\\
%\cos\left((l-m) ( \frac{ 2 \pi I_k}{N} + \frac{\pi}{2}) \right) \sin\left((l-m)( \frac{ 2 \pi I_k}{N} + \frac{\pi}{2}) \right) & \sin^2\left((l-m)( \theta_i + \frac{\pi}{2})\right)
%\end{pmatrix} \,.$
%\beqn
%L_{lm}
%= \frac{16 R}{ \pi} \,
%\sum_{k=1}^M  \iota \left[ e^{i (m-l) \frac{ 2 \pi I_k}{N} } \right]
%\frac{1}{ 4 m_0 + 2 J_k + 1 } \cos \left( \left( m_0 + \frac{ J_k - m  }{2} \right) \pi \right) \cos \left( \left( m_0 + \frac{ J_k - l  }{2} \right)  \pi \right)
% + O( m_0^{-2}) \,.  \notag
%\eqn

We are now ready to specify our choice of indices $\{ ( I_i, J_i) \}_{i=1}^M$.  In particular we let the array $\{ ( I_i, J_i) \}_{i=1}^M$ be such that it enumerate the index set $\{ (I,J): 1\leq I \leq N, 1\leq J \leq N\}$, i.e., we have $M = N^2$.
With the above definition, we readily see 
\beqn
& & L_{lm} \notag\\
%&=& \frac{16 R}{ \pi} \, \sum_{J=1}^N \sum_{I=1}^N \iota \left[ e^{i (m-l) \frac{ 2 \pi I}{N} } \right]
%\frac{1}{ 4 m_0 + 2 J + 1 } \cos \left( \left( m_0 + \frac{ J - m  }{2} \right) \pi \right) \cos \left( \left( m_0 + \frac{ J - l  }{2} \right)  \pi \right)
% + O( m_0^{-2}) \,.  \notag \\
&=& \frac{16 R}{ \pi} \,
\sum_{J=1}^N  \left[ \sum_{I=1}^N  L_{l,m, \frac{ 2 \pi I}{N} } \right]
\frac{1}{ 4 m_0 + 2 J + 1 } \cos \left( \left( m_0 + \frac{ J - m  }{2} \right) \pi \right) \cos \left( \left( m_0 + \frac{ J - l  }{2} \right)  \pi \right)
 + O( m_0^{-2})  \notag \\
%&=& \frac{16 R N}{ \pi} \delta_{lm}  \iota \left[  1 \right] \,
%\sum_{J=1}^N
%\frac{1}{ 4 m_0 + 2 J + 1 } \cos^2 \left( \left( m_0 + \frac{ J - m  }{2} \right) \pi \right)
% + O( m_0^{-2})   \notag \\
&=& \frac{16 R N}{ \pi}  \delta_{|l|,|m|}  \begin{pmatrix}
1  & 0 \\
0 & \text{sgn}(l) \, \text{sgn}(m)
\end{pmatrix} \,
\sum_{J=1}^N
\frac{1}{ 4 m_0 + 2 J + 1 } \cos^2 \left(  \frac{ (J - m) \pi  }{2} \right)
 + O( m_0^{-2}) \,, \notag
\eqn
where $\delta_{a,b}$ is the kronecker delta for any $a,b \in \mathbb{N}$.

>From the above summation, we can directly infer that
\beqn
L_{lm} = \frac{16 R N}{ \pi}  \delta_{|l|,|m|}  \begin{pmatrix}
1  & 0 \\
0 & \text{sgn}(l) \, \text{sgn}(m)
\end{pmatrix} \,
\sum_{J=0}^{[\frac{N-1}{2}]}  \frac{1}{ 4 m_0 + 4 J + 2\, \texttt{mod}_2(m) +1}
 + O( m_0^{-2}) \,,
\eqn
where $\texttt{mod}_2$ is the standard mod-$2$ function and $[\cdot]$ is the floor function.  
Now, for the sake of exposition, we denote for given $C, m_0, \widetilde{M}$ 
a coefficient $K_{C,m_0,\widetilde{M}}$ as
\beqn
K_{C,m_0,\widetilde{M}} := \sum_{J=1}^{\widetilde{M}} \frac{1}{ 4 m_0 +1 + 4 J + 2 C } \,.
\eqn
With this definition, we now hope to approximate $K_{C,m_0,\widetilde{M}}$.
%With a fixed $m_0$ and for a large $\widetilde{M}$, for any given smooth function with compact support $\varphi$, we have that \footnote{I will double-check.}
%\beqn
%\int_{\infty}^\infty K_{C,m_0,\widetilde{M}}(r) \varphi(r) dr &= & \frac{1}{4} \, \textit{p.v.}\int_{-\pi}^\pi \left( \int_{\infty}^\infty \frac{e^{i r x } }{ x  } \varphi(r) dr \right) dx  + O ( | \text{supp}(\varphi)| || \varphi ||_{\infty} \widetilde{M}^{-2}) \notag \\
%& = &\frac{1}{4} \, \textit{p.v.}\int_{-\pi}^\pi \widehat{( \int \varphi ) }[x] dx  + O ( | \text{supp}(\varphi)| || \varphi ||_{\infty} \widetilde{M}^{-2}) \,.
%\eqn
In fact, from the comparison test, we directly arrive at, for any fixed $m_0, \tilde{M}$ and any $C= 0,1$, that the following holds:
\beqn
\frac{1}{4} \log\left(1 + \frac{ \tilde{M} }{m_0 + 2}\right)
%\leq \int_{1}^{[\frac{M}{2}]+1} \frac{1}{ 4 m_0 +1 + 4 x + 2 C } dx
\leq  K_{C,m_0,\tilde{M}}
%\leq \int_{\frac{3}{4}}^{[\frac{M}{2}]+\frac{3}{4}} \frac{1}{ 4 m_0 +1 + 4 x + 2 C } dx
\leq
\frac{1}{4} \log\left(1 + \frac{ \tilde{M} }{m_0 + 1}\right)
\,.
\label{ineq_K}
\eqn
Then we can write 
\beqn
L_{lm} = \frac{16 R N}{ \pi} \delta_{|l|,|m|} K_{\texttt{mod}_2(m) ,m_0,[\frac{N-1}{2}] }  \begin{pmatrix}
1  & 0 \\
0 & \text{sgn}(l) \, \text{sgn}(m)
\end{pmatrix}
 + O( m_0^{-2}) \,,
\eqn
with $K_{\texttt{mod}_2(m) ,m_0,[\frac{N-1}{2}] } $ satisfying estimate \eqref{ineq_K}.
We may now observe a seemingly pathological situation: the matrix $L$ is actually not invertible in $\mathbb{R}^{4N}$. % since we see that the blocks $L_{lm} = \begin{pmatrix}
%1  & 0 \\
%0 & -1
%\end{pmatrix} L_{l \, -m} $ for all $N \leq l,m \leq N$.
However, this is actually not as pathological as we think it is, because the constraint $C x = 0$ shall come in to play a fundamental role.
To proceed, we can take a matrix $E \in C^{\perp}$ in the block form $\left( E_{lm} \right)_{-N\leq l\leq N,  l \neq 0, 1\leq m\leq N,\,}$ as follows:
\beqn
E_{lm} =   \delta_{|l|,m}  \begin{pmatrix}
1  & 0 \\
0 & \text{sgn}(l)
\end{pmatrix} \,.
\eqn
One easily check that the above block matrix $E$ is indeed in $C^{\perp}$.
Then one directly calculate that, for all $1\leq l,m\leq N$,
\beqn
{\left(E^T \circ L\circ E \right)}_{lm} = \frac{64 R N}{ \pi} \delta_{l, m} K_{\texttt{mod}_2(m) ,m_0,[\frac{N-1}{2}] }  \begin{pmatrix}
1  & 0 \\
0 & 1
\end{pmatrix}
 + O( m_0^{-2}) \,,
\eqn
which is now invertible.
Hence for a fixed $N$ and the choice 
$(\theta_i, \widetilde{\theta_i}, k_i)$ of the form $( 2 \pi I_i   / N, 2 \pi I_i   / N + \pi,  \frac{4 m_0 + 2 J_i + 1}{8 R})$, where
$\{ ( I_i, J_i) \}_{i=1}^M$ enumerates through $\{ (I,J): 1\leq I \leq N, 1\leq J \leq N\}$, 
we can directly derive the following estimate for the singular values of $T_{m_0}$:
{\small
\beqnx
\frac{4 \sqrt{R N}}{ \sqrt{\pi}}  \sqrt{\log\left(1 + \frac{ [\frac{N-1}{2}]  }{m_0 + 2}\right) } - \frac{C_N}{m_0^{2}} \leq s_{\min} (T_{m_0} \circ E)  \leq  s_{\max} (T_{m_0} \circ E ) \leq \frac{4 \sqrt{R N}}{ \sqrt{\pi}}  \sqrt{ \log\left(1 + \frac{ [\frac{N-1}{2}]  }{m_0 + 1}\right) } + \frac{C_N}{m_0^{2}}
\eqnx
}where $C_N$ is a constant only depending on $N$. %In here, we would like to remind that this big-$O$ term $O(m_0^{-2})$ can be absorbed when $m_0$ is large (by replacing the respective multiplicative constants with a slightly smaller/larger one.)
Therefore, if we write $s_{max}$ and $s_{min}$ respectively as the largest and smallest singular values, then 
it follows that 
\beqn
\kappa(T_{m_0} \circ E) = \frac{s_{\max} (T_{m_0} \circ E)}{s_{\min} (T_{m_0} \circ E)}  \leq  \sqrt{\log\left(1 + \frac{ [\frac{N-1}{2}]  }{m_0 + 1}\right) \bigg \slash  \log\left(1 + \frac{ [\frac{N-1}{2}]  }{m_0 + 2}\right)  } + O({m_0}^{-2}) .\,
\label{log_growth}
\eqn
The Taylor series of $\log(1+x)$ and $\sqrt{a+x}$ then give rise to the following estimate for large $m_0$ that
\beqn
 \kappa_{\text{inf},N}  \leq \kappa(T_{m_0}, \text{ker}(C) ) \leq  \kappa(T_{m_0} \circ E)  \leq  \sqrt{\frac{ m_0+2 }{m_0 + 1} } + O(m_0^{-2}) \leq 1 + O(m_0^{-1}) \,,
\eqn
where we should remind ourselves that the big-$O$ terms are bounded by a constant only depending on $N$. Since $m_0$ is arbitrary, we get for any given $N$ that the infimum of the condition number $\kappa(T, \text{ker}(C) )$ is given by
$
 \kappa_{\text{inf},N}  =1\,,
$
and a minimizing sequence to attain this infimum can be actualized by measurements
$(\theta_i, \widetilde{\theta_i}, k_i)$ as previously specified as $m_0$ goes to infinity. This implies that we can always make an appropriate choice of the target resolution $N$ such that the inversion process of $T$ is well-posed.
%Moreover, in view of the fact that the constant in front of the big-$O$ term in $\label{log_growth}$ come from \eqref{order0}, we can see that with a same choice of large $m_0$, the growth of condition number is bounded by order $\log(N)$.
The above analysis can be summarized into the following theorem.
\begin{Theorem} \label{theorem_imp}
For a given target resolution $N$, the infimum $\kappa_{\text{inf},N} $ of the condition number $\kappa(T, \text{ker}(C) )$ defined as in \eqref{minumum} over the set of linear operators $T$ is given by
\beqn
 \kappa_{\text{inf},N}  = 1.
\eqn
A minimizing sequence $\kappa(T_{m_0}, \text{ker}(C) ), m_0 \in \mathbb{Z} $ of this infimum acquires the following bound
\beqn
 \kappa(T_{m_0}, \text{ker}(C) )  \leq 1 + O(m_0^{-1})
\label{drastic}
\eqn
if we make the arrangement of phaseless measurements in the way that the following equality holds:
\beqn
(\theta_i, \widetilde{\theta_i}, k_i) = \left( 2 \pi I_i   / N, 2 \pi I_i   / N + \pi,  \frac{4 m_0 + 2 J_i + 1}{8 R} \right) \,,
\label{measurement}
\eqn
where
$\{ ( I_i, J_i) \}_{i=1}^M$ enumerates through $\{ (I,J): 1\leq I \leq N, 1\leq J \leq N\}$ and $m_0$ is large, 
hence $N^2$ phaseless measurements shall be made.
\end{Theorem}
This theorem gives us a very effective strategy of data measurement such that the phaseless reconstruction process shall be well-posed. In particular, an increase of $m_0$ in the aforementioned measurement method reduces 
the condition number of the inversion process with an order of $O(m_0^{-1}) $ according to \eqref{drastic}.

\subsection{Estimation of the condition number of $L$} \label{cond_L}

>From the previous analysis, we can see that the inversion process of $T$ can be made impressively stable and one can suppress its condition number appropriately.  However this does not ensure 
a very stable phaseless inversion process, owing to the fact from \eqref{operatoreq} that 
the total inversion process is given by $T^{-1} \circ L^{-1}$.
%and that the inversion $L^{-1}$ might be indeed tragic.

Although the action of $L$ is simple and explicit, the inversion process may not be as simple as one might think.
The condition number of $L$ can be directly calculated as ${ \max_i | P_R(\theta_i, \widetilde{\theta}_i, k_i)|}/{ \min_i |P_R(\theta_i, \widetilde{\theta}_i, k_i)|}$. Therefore the inversion process becomes severely ill-posed when some measurement data has a very small value $| P_R(\theta_i, \widetilde{\theta}_i, k_i)|$, which in turn pushes up the condition number to an arbitrary magnitude. This causes the reconstruction process to be very unstable in practice.

However, a very simple regularization technique can get rid of this instability.
Thanks to the fact that $P_R(\theta, \widetilde{\theta}, k)$ is analytic, its value cannot be zero on an open neighborhood, and therefore a simple regularization can be performed on the inversion of $L$ by the operator $L_\alpha^{-1}$ defined as follows:
{\small
\beqn
L_\alpha^{-1} = \text{diag} \Big( \chi_{_{x > \alpha} } ( | P_R(\theta_i, \widetilde{\theta}_i, k_i) | ) [P_R(\theta_i, \widetilde{\theta}_i, k_i)]^{-1} + {\alpha^{-1}}  \chi_{_{x \leq \alpha} } ( | P_R(\theta_i, \widetilde{\theta}_i, k_i) | )\lim_{(\theta, \widetilde{\theta}, k) \rightarrow (\theta_i, \widetilde{\theta}_i, k_i)}\frac{P_R(\theta, \widetilde{\theta}, k)}{|P_R(\theta, \widetilde{\theta}, k)|} 
\Big)\label{labL}
\eqn
}where $ \chi_{_{x > \alpha} } $ and $ \chi_{_{x \leq \alpha} } $ are the respective characteristic functions on the intervals $\{x > \alpha \}$ and  $\{x \leq \alpha  \}$.
With this definition, we come readily to the following simple but important lemma.
\begin{Lemma}
Let $L_\alpha^{-1}$ be defined as in \eqref{labL}, then we have
\beqn
\kappa (L^{-1}) = \frac{ \max_i | P_R(\theta_i, \widetilde{\theta}_i, k_i)|}{ \min_i |P_R(\theta_i, \widetilde{\theta}_i, k_i)|} \,, \quad \kappa (L_\alpha^{-1}) \leq \frac{ 2 }{ \alpha } \,.
\label{linverse}
\eqn
\end{Lemma}
We can see from above that $\kappa (L^{-1})$ cannot be controlled but $\kappa (L_\alpha^{-1})$ has an upper bound, 
therefore it provides a stable inversion process if $\alpha$ is appropriately chosen.

>From \eqref{operatoreq}, a stable shape reconstruction process is therefore provided by $ T^{-1} \circ L_{\alpha}^{-1} $.
Indeed, the stability estimates \eqref{drastic} and \eqref{linverse} for the condition numbers of $ T^{-1} $ and $ L_{\alpha}^{-1} $ subjected to $C x = 0$ ensure us the stability of this reconstruction method and provide us optimal strategies to lower the degree of ill-posedness for the phaseless reconstruction problem under the corresponding measurement cases.
The stability of our proposed method will be verified in the numerical experiments.
To the best of our knowledge, these estimates of condition numbers are completely new to our inverse problems. 

%\noindent \textbf{Remark}: Again we may also consider the case of finite incident measurement. The inversion process is similar and estimates can be done in a similar fashion. However the results of the estimates may be much worse owing to a lacking of enough measurements.

\subsection{A comparison with the phased reconstruction} \label{sec74}
As we have remarked in section \ref{sec3},
together with the fact that any translation of the inclusion yields the same phaseless measurement, 
the phaseless reconstruction is not unique in this sense. And 
the linearized phased and phaseless reconstructions share some fundamental differences.
>From \eqref{phased} or its Fourier-transformed version (Algorithm 1), we see that any algorithm derived from \eqref{phased} for the phased reconstruction is equivalent to solving $ \delta \left( \mathfrak{F}[h]^{est} (l) \right)_{l = -N}^N \in \mathbb{C}^{2N} \oplus \{ 0\} \cong \mathbb{R}^{4N} $ such that
\beqn
 \tilde{T}  \left[ (\delta \,\mathfrak{F}[h])^{est} \right] = G( A_{\infty}^{\text{meas}} (\theta_i, \widetilde{\theta}_i, k_i) ) \,,
\eqn
where $N$ satisfies
\beqnx
C \frac{N^{4N}}{R^{2+4N}} < \text{(SNR)}^{1 + \alpha/2}
\eqnx
for some $C$ and $\alpha$,
and $(\mathfrak{F}[h])^{est} $ is again subjected to the constraints
\beqnx
\text{Re}(\mathfrak{F}[h])^{est}( - l) - \text{Re} (\mathfrak{F}[h])^{est}( l) = 0 \,, \qquad  \text{Im}(\mathfrak{F}[h])^{est}( - l) + \text{Im} (\mathfrak{F}[h])^{est}( l) = 0\,.
\eqnx
The operators $G$ and $\tilde{T}$ above are respectively given by 
\begin{enumerate}
\item
$G:\mathbb{C}^M \rightarrow \mathbb{C}^M$, the component-wise affine map of a vector, i.e. the action
$v_i \mapsto \frac{v_i - \pi R^2 {\varepsilon^*} k^2  P_R(\theta_i, \widetilde{\theta}_i, k_i)}{2 \pi R^2 {\varepsilon^*}k^2 }$.
\item
$\tilde{T}  : \mathbb{C}^{2N}  \rightarrow \mathbb{C}^{M}$, the linear operator $ v \mapsto \left( \langle \, v \, , S_R(\theta_i, \widetilde{\theta}_i, k_i) \rangle_{l^2(\mathbb{C})} \right)_{i = 1}^M$.
\end{enumerate}

A similar stability analysis for the operator $\tilde{T}$ induced by the phased measurements can be performed 
to the one for the operator $\tilde{T}$ corresponding to the phaseless reconstruction as in section \ref{sec71}.
Since most of the steps are similar to the previous analysis for the phaseless reconstruction, we only provide a sketch of the argument.
Again we choose $(\theta_i, \widetilde{\theta_i}, k_i)$  of the form $( 2 \pi I_i   / N ,  2 \pi I_i   / N  + \pi,  \frac{4 m_0 + 2 J_i + 1}{8 R})$ where $\{ ( I_i, J_i) \}_{i=1}^M \in \mathbb{Z}$ are some integral indices to be specified and let $\tilde{T}_{m_0}$ be the linear operator $\tilde{T}$ with this specific arrangement of measurement with a given $m_0$.
Denoting $\tilde{L} := \iota [ \tilde{T}_{m_0}^* ] \iota  [\tilde{T}_{m_0}]$, then a similar argument, 
along with the fact that $\iota$ is an algebra homomorphism,  shows for $-N\leq l,m\leq N$ that 
$$
\tilde{L}_{lm}
= \frac{16 R}{ \pi} \, \iota \left[
\sum_{i=1}^M
\tilde{L}_{l,m,\theta_i}
\frac{1}{ 4 m_0 + 2 J_i + 1 } \cos \left( \left( m_0 + \frac{ J_i - m  }{2} \right) \pi \right) \cos \left( \left( m_0 + \frac{ J_i - l  }{2} \right)  \pi \right)
\right] + O( m_0^{-2}) \,,  \notag
$$
where each $\tilde{L}_{l,m,\theta_i} := e^{i (l-m) \theta_i} $ is invertible.
Again, letting the array $\{ ( I_i, J_i) \}_{i=1}^M$ enumerate the index set $\{ (I,J): 1\leq I \leq N, 1\leq J \leq N\}$, i.e., $M = N^2$ complex (phased) measurements, we have
\beqn
& & \tilde{L}_{lm} \notag \\
&=& \frac{16 R}{ \pi} \,
\sum_{J=1}^N  \iota \left[ \sum_{I=1}^N   e^{2 \pi i (l-m)  I / N } \right]
\frac{\cos \left( \left( m_0 + \frac{ J - m  }{2} \right) \pi \right) \cos \left( \left( m_0 + \frac{ J - l  }{2} \right)  \pi \right)}
{ 4 m_0 + 2 J + 1 } 
 + O( m_0^{-2})   \notag\\
&=& \frac{16 R N}{ \pi}  \delta_{l, m} \begin{pmatrix}
1  & 0 \\
0 & 1
\end{pmatrix}
\sum_{J=1}^N
\frac {\cos^2 \left(  \frac{ (J - m) \pi  }{2} \right)}{ 4 m_0 + 2 J + 1 } 
 + O( m_0^{-2}) \,. \notag
\eqn
>From here onward, the analysis is the same as in section \ref{sec71} to get the same block matrix $E$ such that for all $1\leq l,m\leq N$, 
\beqn
{\left(E^T \circ \tilde{L} \circ E \right)}_{lm} = \frac{32 R N}{ \pi} \delta_{l, m}  K_{\texttt{mod}_2(m) ,m_0,[\frac{N-1}{2}] }
\begin{pmatrix}
1  & 0 \\
0 & 1
\end{pmatrix}
 + O( m_0^{-2}) \,.
 \label{eqneqn}
\eqn
Now a similar argument as in section \ref{sec71} applied to get an identical result for the phased reconstruction:
\beqn
 \kappa(\tilde{T}_{m_0}, \text{ker}(C) )  \leq  \sqrt{\frac{ m_0+2 }{m_0 + 1} } + O(m_0^{-2}) \leq  1 + O(m_0^{-1}) \,,
\eqn
and by tracing all the constants, we can see the constants represented by big-$O$'s is of the same magnitude as in the phaseless reconstruction.
Therefore the ill-posedness in inverting $T$ and $\tilde{T}$ are actually of the same order of magnitude using a same set of measurement angles, and the following result holds.

\begin{Theorem} \label{theorem_imp2}
For a given target resolution $N$, the condition number $\kappa(\tilde{T}_{m_0}, \text{ker}(C) )$ of the operator $\tilde{T}_{m_0}$ for $m_0 \in \mathbb{Z}$ can be controlled by
\beqn
\kappa(\tilde{T}_{m_0}, \text{ker}(C) ) \leq 1 + O(m_0^{-1})
\label{drastic2}
\eqn
if we make an $N^2$ complex (phased) measurement arrangement:
\beqn
(\theta_i, \widetilde{\theta_i}, k_i) = \left( 2 \pi I_i   / N, 2 \pi I_i   / N + \pi,  \frac{4 m_0 + 2 J_i + 1}{8 R} \right) \,,
\label{measurement12}
\eqn
where
$\{ ( I_i, J_i) \}_{i=1}^M$ enumerates through $\{ (I,J): 1\leq I \leq N, 1\leq J \leq N\}$.
\end{Theorem}

Nonetheless, we notice a fundamental difference here between the phased and phaseless reconstructions. For the phased reconstruction, the matrix $\tilde{L}$ is invertible itself, therefore the constraint $C x = 0$ is redundant. However, in the phaseless reconstruction, this set of constraints is necessary for us to get to a solution in the inversion process.
Therefore, to fully exploit the constraints $Cx = 0$, it shall be possible to obtain the same stability estimate for $\tilde{T}$ even if the number of equations represented by the matrix are cut off by half.
There are different ways to realize this, and we suggest one of them below.
We shall not repeat all the details in the argument again but give only a sketch.

%\footnote{Please check this paragraph very carefully to see if it is right.}
Suppose we choose the set of measurement points
$(\theta_i, \widetilde{\theta_i}, k_i)$ as $\left( 2 \pi I_i   / N ,  2 \pi I_i   / N  + \pi,  \frac{4 m_0 + 2 J_i + 1}{8 R}\right)$ where $\{ ( I_i, J_i) \}_{i=1}^M $ enumerate the index set $\{ (I,J): 1\leq I \leq N/2, 1\leq J \leq N\}$, but 
we only measure the real part of the far-field pattern $A_{\infty}^{\text{meas}} (\theta_i, \widetilde{\theta}_i, k_i)$.
Clearly, we have $N^2$ real (phased) measurements.
>From the fact that $P_R(\theta_i, \widetilde{\theta}_i, k_i)$ is real, we have
\beqn
\text{Re} \left( \tilde{T}  \left[ (\delta \,\mathfrak{F}[h])^{est} \right] \right)= \text{Re} \left( G( A_{\infty}^{\text{meas}} (\theta_i, \widetilde{\theta}_i, k_i) ) \right) =  G\left( \text{Re} ( A_{\infty}^{\text{meas}} (\theta_i, \widetilde{\theta}_i, k_i) ) \right) \,.
\eqn
Therefore, by taking only $N^2$ real (phased) measurements, we are actually dropping half of the equations representing measurements from the imaginary part.
Now, in order to distinguish from the previous measurement setting, we denote the operator with these new measurement events as $\tilde{\tilde{T}}_{m_0}$ for a given $m_0$.

With this very particular choice of real (phased) measurements, we know from \eqref{special}
that the matrix $\tilde{\tilde{T}}_{m_0}$ is coincidentally the same as $T_{m_0}$.
Hence, if we write $\tilde{\tilde{L}} := \tilde{\tilde{T}}_{m_0}^T  \tilde{\tilde{T}}_{m_0}$, then 
$\tilde{\tilde{L}} = L$. Therefore, 
with the same $E$ as previously chosen, the same argument applies for us to get for all $1\leq l,m\leq N$ that 
\beqn
{\left(E^T \circ \tilde{\tilde{L}} \circ E \right)}_{lm} = \frac{64 R N}{ \pi}   \delta_{l, m}  K_{\texttt{mod}_2(m) ,m_0,[\frac{N-1}{2}] }
 + O( m_0^{-2}) \,.
 \label{eqneqneqn}
\eqn
This gives the following result.
\begin{Theorem} \label{theorem_imp3}
An effective choice of only $N^2$ real phased measurement ensures 
the following bound for the condition number:
\beqn
 \kappa(\tilde{\tilde{T}}_{m_0}, \text{ker}(C) )  \leq 1 + O(m_0^{-1}) \,.
\label{drastic3}
\eqn
\end{Theorem}
Other ways to fully exploit the constraints $Cx = 0$ by dropping at most half of the equations represented by $\tilde{T}$, such as measuring the projection of complex number by another phase angle other than taking the real part, or taking only a special set of undersampling measurements, shall be possible, but for the sake of simplicity, we shall not proceed
further.

>From the above analysis, we can see that although the structures of $T$ and $\tilde{T}$ are fundamentally different, they 
have similar behavior on their condition numbers.
Yet the phaseless reconstruction is still much more ill-posed than the phased counterpart, owing to the following very simple yet important point.
In the phaseless reconstruction, we shall also need to invert $L$ by a regularized inversion process $L_\alpha^{-1}$; however, in a phased reconstruction, such an inversion of $L$ is unnecessary. Therefore, instability imposed by $L$ exists only in the phaseless reconstruction.  Considering this fact, the total regularized inversion of the phaseless reconstruction is still much more ill-posed than the phased counterpart, having its condition number being $1/\alpha$ times of the phased reconstruction.

\section{Numerical experiments} \label{numerical}

In this section, we will first present numerical results illustrating some behaviors of the condition number $\kappa(T_{m_0}, \text{ker}(C) )$ using our measurement strategy described in the section \ref{condT}, then focus on the inverse problem of shape reconstruction from the observed magnitude of far-field data.

\subsection{Condition number of $T$ subjected to $C x = 0$}
In what follows, %we shall first compare the behaviours of condition number $\kappa(T)$ in several random measurement case, and numerically verify the probabilistic asymptotic estimate provided in Theorem \ref{theorem_imp1}.
we shall observe the behaviors of the condition number $\kappa(T_{m_0}, \text{ker}(C) )$ using our measurement strategy given in Theorem \ref{theorem_imp2} and check the asymptotic estimate of $\kappa(T)$ in the theorem as $m_0$ grows.
With a given $m_0$, we now fix the resolution $N=51$ and choose the wave-numbers $k$ such that $k =  \frac{ 4 m_0 + 2 J + 1 }{8 R}$ with $R = 0.2$ and $J = 5, ..., 10$. The measurement points are the same as stated in Theorem \ref{theorem_imp2}.
We compute the condition number of the operator $T$ with $m_0 = 1, ... ,20$. The values of the corresponding condition numbers are plot in Figure \ref{cond_stra}.

\begin{figurehere} \centering
\includegraphics[width=5cm,height=4cm]{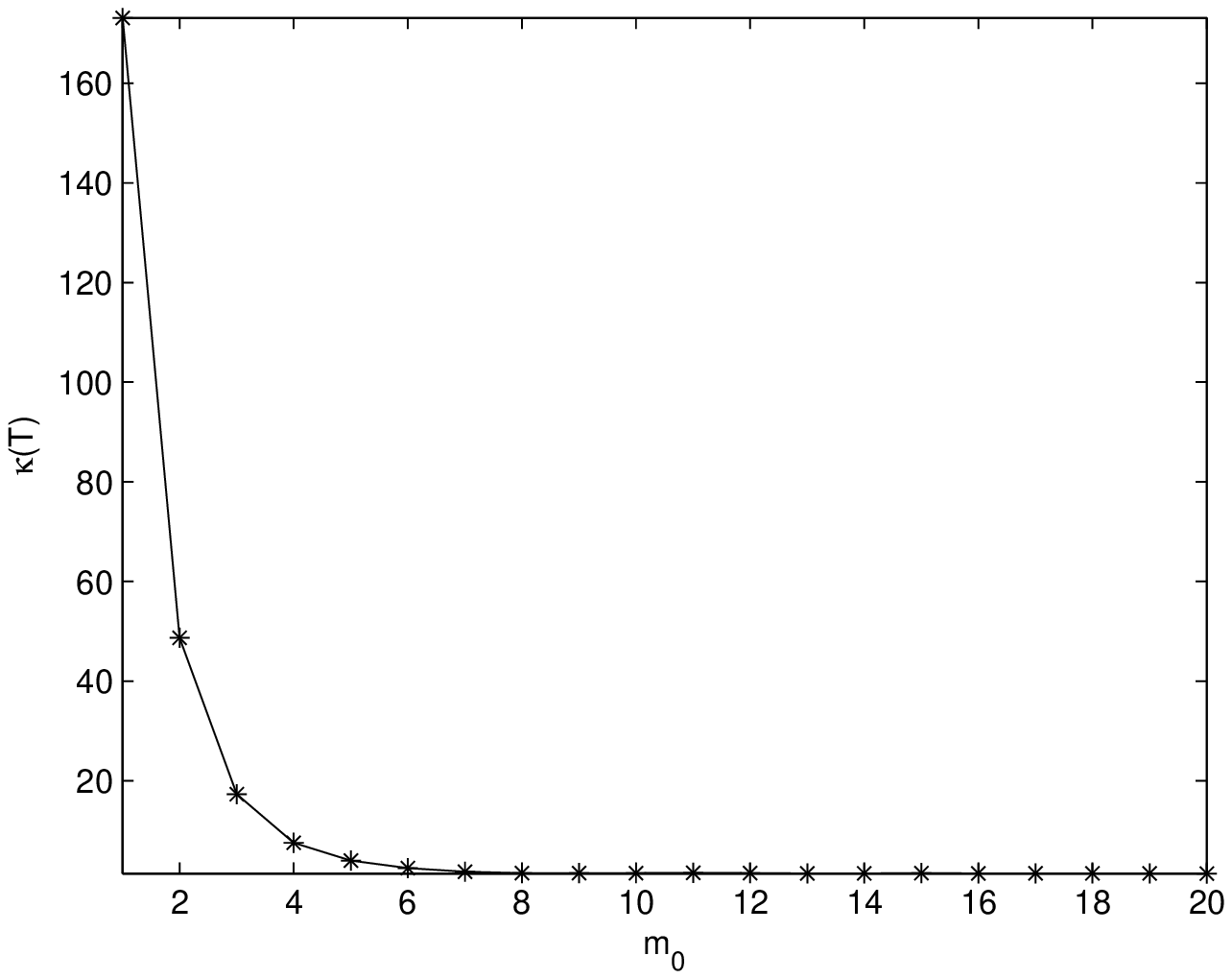}
\caption{Decay of $\kappa(T, \text{ker}(C))$ with respect to $m_0$.} \label{cond_stra}
\end{figurehere}
We see clearly the drastic decay of the condition number as $m_0$ grows, showing the effectiveness of increasing stability by the increment of $m_0$.
This agrees with the result we obtained in Theorem \ref{theorem_imp2}.

\subsection{Phaseless reconstruction}

We shall now proceed to present several numerical examples to show the performance
of the newly proposed reconstruction algorithm, i.e., Algorithm $2$ in Section\,\ref{phaseless_algo}, 
from phaseless far-field data.  In order to attain robustness and stability of our algorithm, we approximate 
the inversion of $L$ in Step 2 by $L_{\alpha}$ as in section \ref{cond_L} for some regularization $\alpha$ described below.

In the following $3$ examples, we consider an infinite homogeneous background medium with its material coefficient being $1$. In each example, an inhomogeneous inclusion $D = B^\delta$ is then introduced as a perturbation of a circular domain $B = B(0,R)$ for some $\delta>0$ and its radius $R=0.2$ sitting inside the homogeneous background medium, with its contrast always set to be $\varepsilon^* = 0.05$.

Given a domain $B^\delta$, we first obtain the observed data of the forward problem, namely the magnitude of far-field data.
In order to generate the far-field data for the forward problem and the observed scattering coefficients, we use the SIES-master package developed by H.~Wang \cite{hansite}.
For a fixed wave-number $k$, we first solve for the solutions $(\phi_m, \psi_m)$ of \eqref{defint} for $|m|\leq 50$ using the rectangular quadrature rule with mesh-size $s/1024$ along the boundary of the target, where $s$ denotes the length of the inclusion boundary. The scattering coefficients of $B^{\delta}$ of orders $(n,m)$ for $|n|,|m| \leq 50 $ are then calculated, 
and the far-field data $A(\theta_d, \theta_x,k)$ is evaluated 
using \eqref{fouriersum} with $\theta_d, \theta_x \in (0,2\pi]$ on a uniform mesh of size $N = 50$.
Then the magnitude of the far-field pattern $| A(\theta_d, \theta_x,k)|$ is taken for our reconstruction process. In order to test the robustness of our reconstruction algorithm against the noise, we introduce some
multiplicative random noise in the magnitude of far-field pattern $| A(\theta_d, \theta_x,k)|$ point-wisely in the form:
\begin{equation}
| A^{\text{meas}}(\theta_d, \theta_x,k)|^{\sigma}
= | A^{\text{meas}}(\theta_d, \theta_x,k)| (1 + \sigma\, \xi ) \label{noise} \end{equation}
where $\xi $ is uniformly distributed between $[-1,1]$ and $\sigma$ refers to the relative noise level.
In the following $4$ examples, we always set the noise level to be $\sigma = 5 \%$.

Then we apply our reconstruction algorithm for shape reconstruction with the noisy phaseless data as $ T^{-1} \circ L_{\alpha}^{-1} \circ F $ following the notation introduced in section \ref{sec_sta}, where the regularization parameter is chosen as $\alpha = 10^{-3}$.
In view of \eqref{measurement}, we make the choice of measurements such that the measured wave-numbers $k$ satisfy $k =  \frac{ 4 m_0 + 2 J + 1 }{8 R}$
with $m_0 = 10$ in all the examples, and $J = 5, ..., 5+\tilde{C}$ for some $\tilde{C}$ to be chosen in each of the example.
%As the phaseless data are invariant under scaling, rotation and translation,
%a postprocessing of scaling, rotating and translating is performed for better comparisons between the reconstructed
%and exact shapes.
The relative error of the reconstruction is defined by
\beqn
\text{Relative Error} :=\frac{\text{Area}\left( (D^{\text{approx}} \bigcup D )\backslash (D^{\text{approx}} \bigcap D ) \right)}{\text{Area}(D)} \,,
\label{error}
\eqn
where $D^{\text{approx}}$ is the reconstructed domain of the exact one $D$.
To demonstrate the effectiveness of our algorithm and illustrates the necessity of a certain number of measurements angles in the phaseless reconstruction (i.e. to test its resolution limit), we shall try $3$ different sets of measurements angles:

\begin{enumerate}
\item[\textbf{Set 1}]
Full measurement angles (over-abundant number of measurements):
\beqn
\left((\theta_d)_i, (\theta_x)_i\right) = \left( 2 \pi I_i   / N_0, 2 \pi K_i   / N_0 \right)\,,\quad 1 \leq I_i, K_i \leq N_0 \,,
\eqn
when $N_0$ is always chosen as $50$ in all the examples;
\item[\textbf{Set 2}]
Transmission measurement angles (critical number of measurements):
\beqn
\left((\theta_d)_i, (\theta_x)_i\right) = \left( 2 \pi I_i   / N, 2 \pi I_i   / N + \pi + U \right) \,,\quad 1 \leq  I_i \leq N  \,,
\eqn
where $N := \min\left\{N: [\mathfrak{F}(h)](k) = 0\,, \quad |k| > N \right\}$ and $h$ is the perturbation in the corresponding example and $U = (-\frac{\pi}{5},\frac{\pi}{5})$;
\item[\textbf{Set 3}]
Half of transmission measurement angles (insufficient number of measurements):
\beqn
\left((\theta_d)_i, (\theta_x)_i\right) = \left( 2 \pi I_i   / [[N/2]], 2 \pi I_i   / [[N/2]] + \pi + U  \right) \,,\quad 1 \leq  I_i \leq [[N/2]]  \,,
\eqn
which consists of $[[N/2]]$ measurement angles, where $N$ is the same as previously mentioned and $[[\cdot]]$ is the ceiling function.
\end{enumerate}
The purpose of introducing an interval $U$ instead of one single point is to increase numerical stability in reconstruction.
We emphasize that \textbf{Sets 2} and \textbf{3} are set up only to test the resolution limit of our phaseless reconstruction algorithm. We are not suggesting the necessity to determine $\min\left\{N: [\mathfrak{F}(h)](k) = 0\,, |k| > N \right\}$ from $h$ before utilizing our algorithm.  Such information is unnecessary and unavailable in a practical phaseless reconstruction.

In order to further increase numerical stability using a critical number of measurements (\textbf{Set 2}) and an insufficient number of measurement (\textbf{Set 3}), we further regularize our inversion process by a $L^1$ regularizer to enforce sparsity in the Fourier modes of our reconstructed perturbation, i.e., we solve 
\beqn
\min_{C \textbf{X}=0 }|| (L_\alpha \circ T)   \textbf{X}  - F(| A_{\infty}^{\text{meas}} (\theta, \widetilde{\theta}, k)| ) ||_2^2 + \beta || \textbf{X}  ||_1
\eqn
where $\beta$ is a regularization parameter that is always chosen as $\beta = 0.05$.  We perform the $L^1$ minimization by a standard Bregman iteration \cite{Bregman}.

\textbf{Example 1}.
In this example, we consider an inhomogeneous domain of a flori-form shape $D = B^\delta$ described by the following parametric form (with $\delta = 0.1$ and $n = 3$ ):
    \begin{equation}
        r = 0.2 ( 1 + \delta \cos (n \theta) )\, , \quad \theta \in (0,2\pi]\,,
        \label{circle_perturb_form}
    \end{equation}
which is a perturbation of the domain $B = B(0, 0.2)$; see Figure \ref{test1A} (left) and (right) 
respectively for the shape of the domain and the contrast of the inhomogeneous medium. 
% and Figure \ref{test1A} (right) for a comparison between the domains $B^\delta$ and $B$ using the values of a sum of characteristic functions $\chi_{B} + \chi_{B^\delta}$.

\begin{figurehere} \centering
\includegraphics[width=5cm,height=4cm]{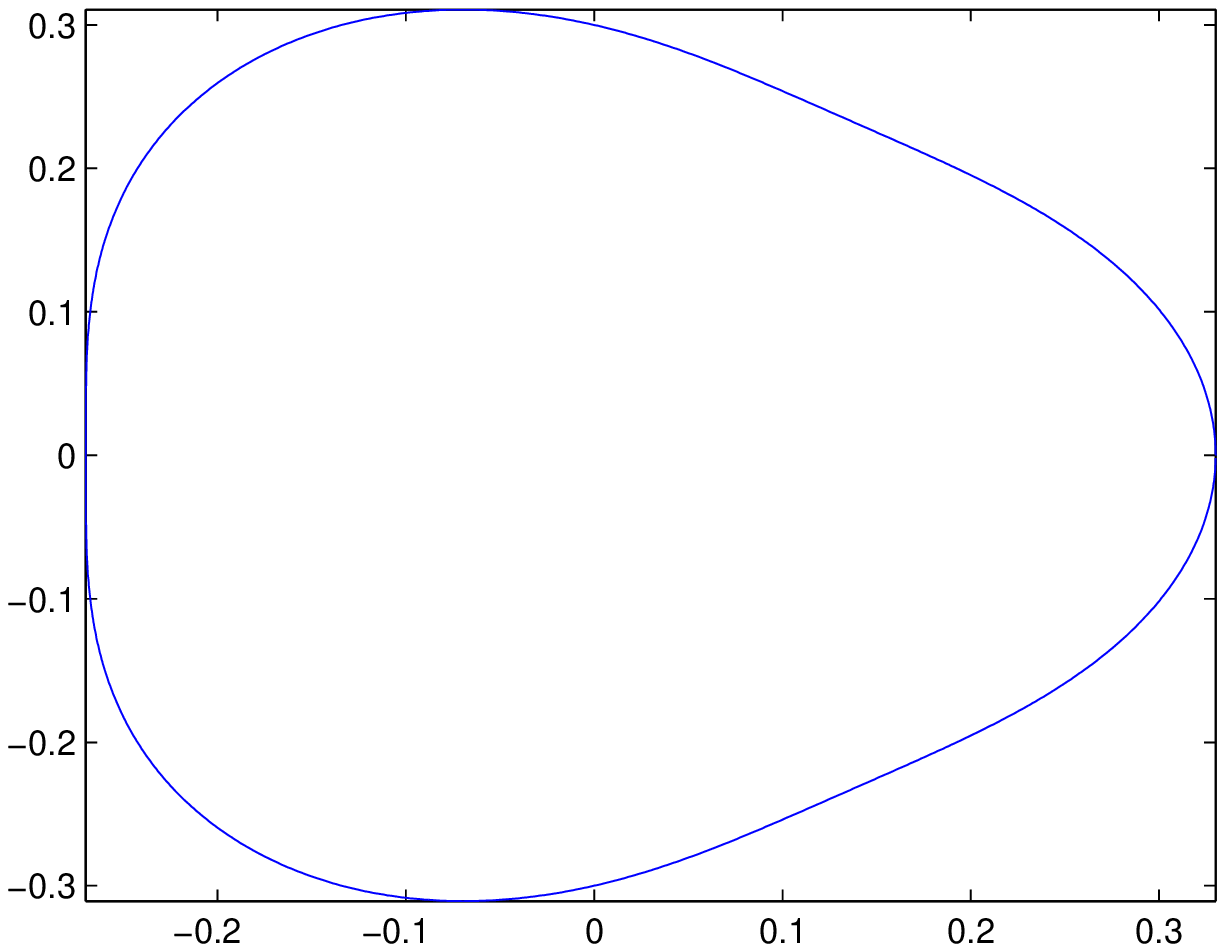}
\includegraphics[width=5cm,height=4cm]{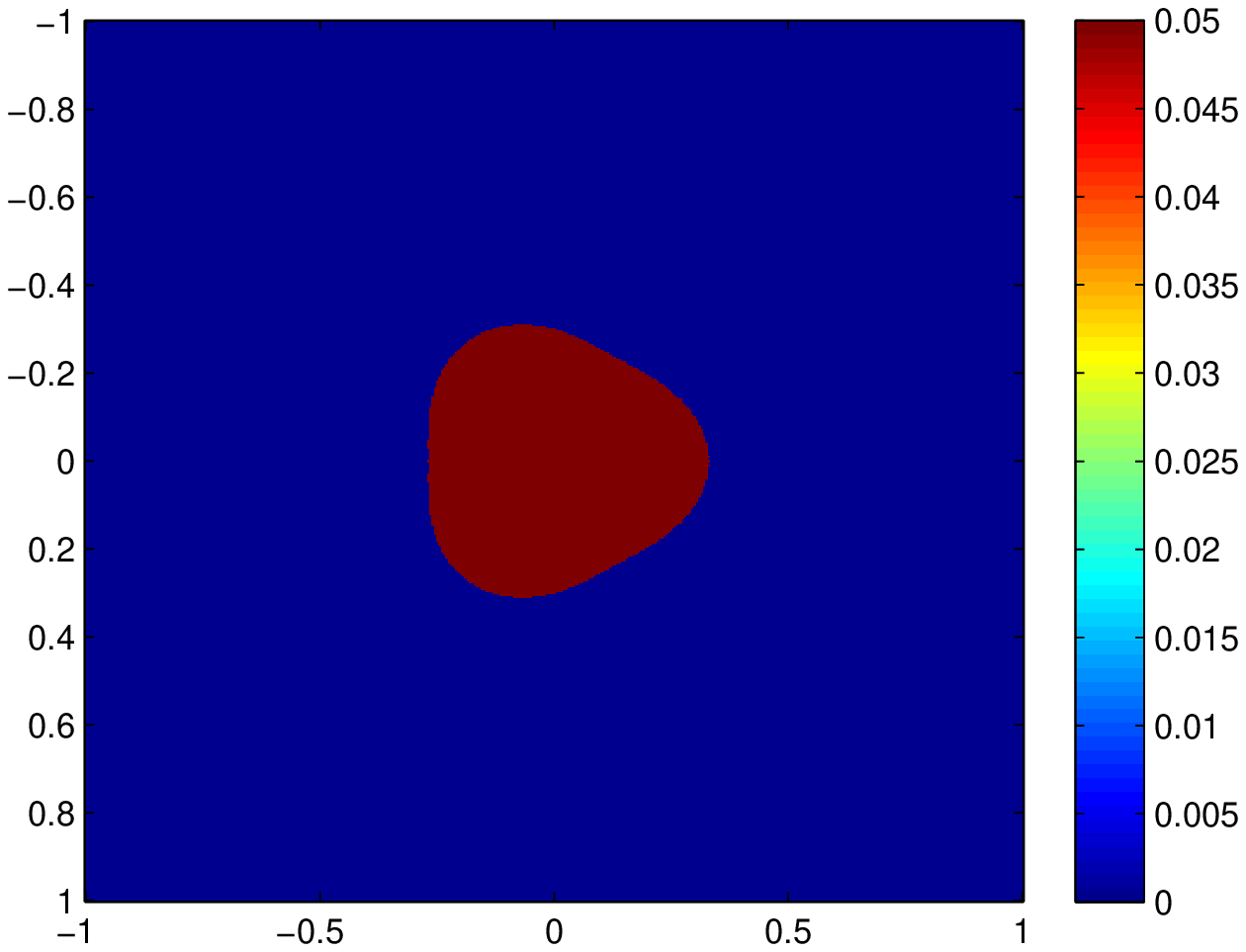}
\caption{Exact inhomogeneous domain (left) and the contrast of the inclusion (right) in Example 1.}
% left: Exact shape; right: exact inclusion.} 
\label{test1A}
\end{figurehere}

The magnitude of the far-field pattern for $6$ wave-numbers are used for shape reconstruction, i.e. $\tilde{C} = 5$, 
and the Fourier coefficients of the reconstructed perturbations using the respective measurement sets are 
shown in Figure \ref{test1B}.

\begin{figurehere} \centering
\includegraphics[width=5cm,height=4cm]{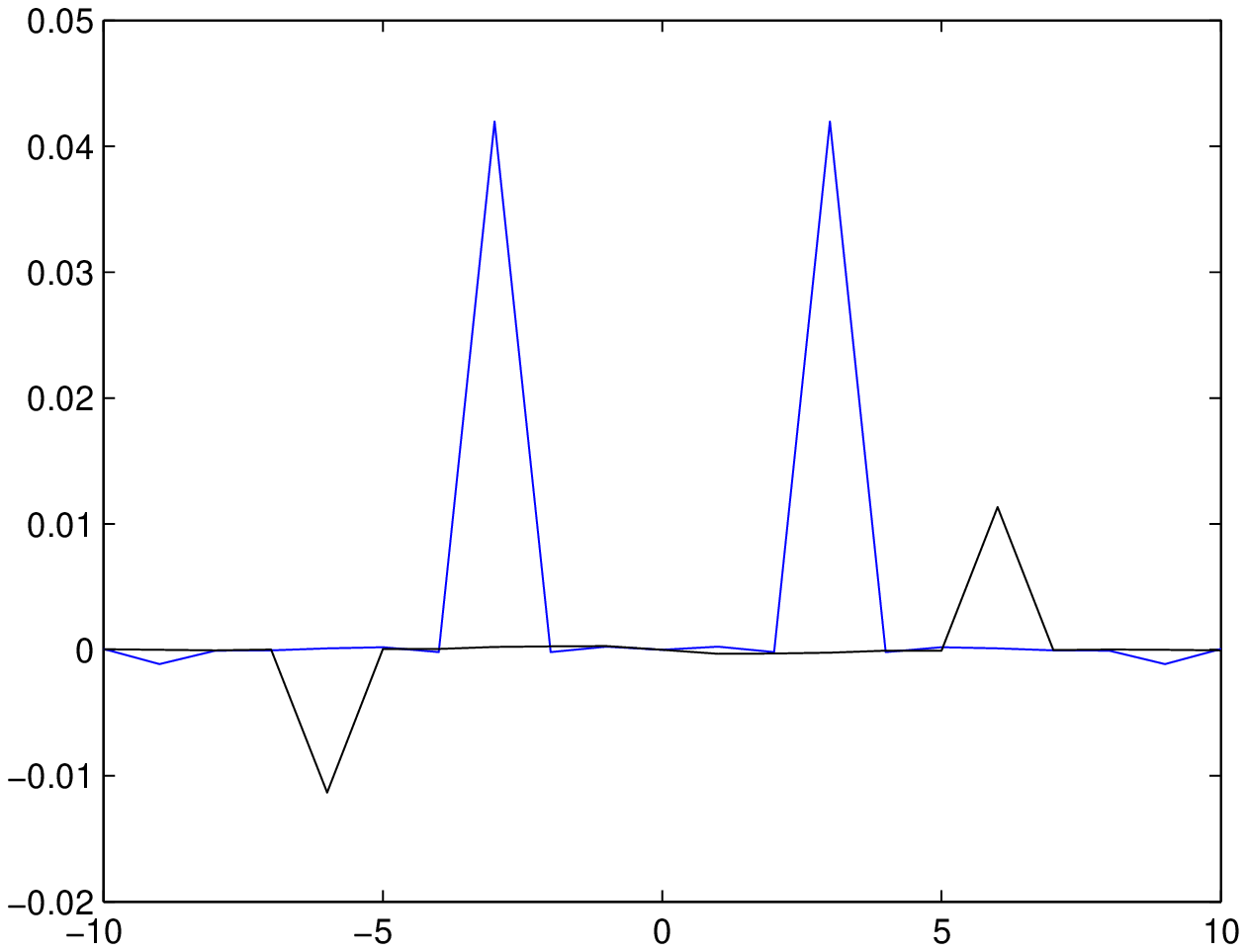}
\includegraphics[width=5cm,height=4cm]{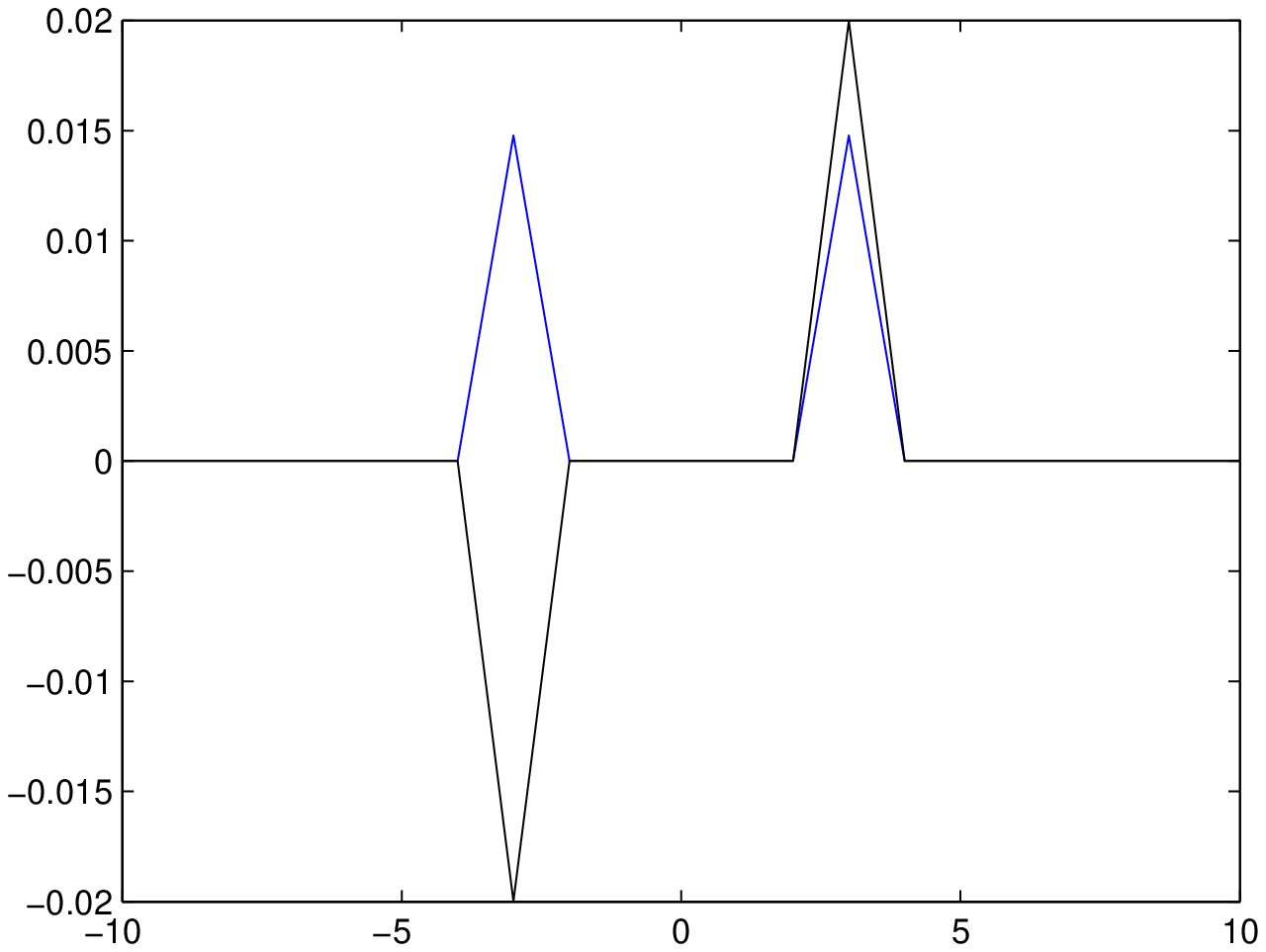}
\includegraphics[width=5cm,height=4cm]{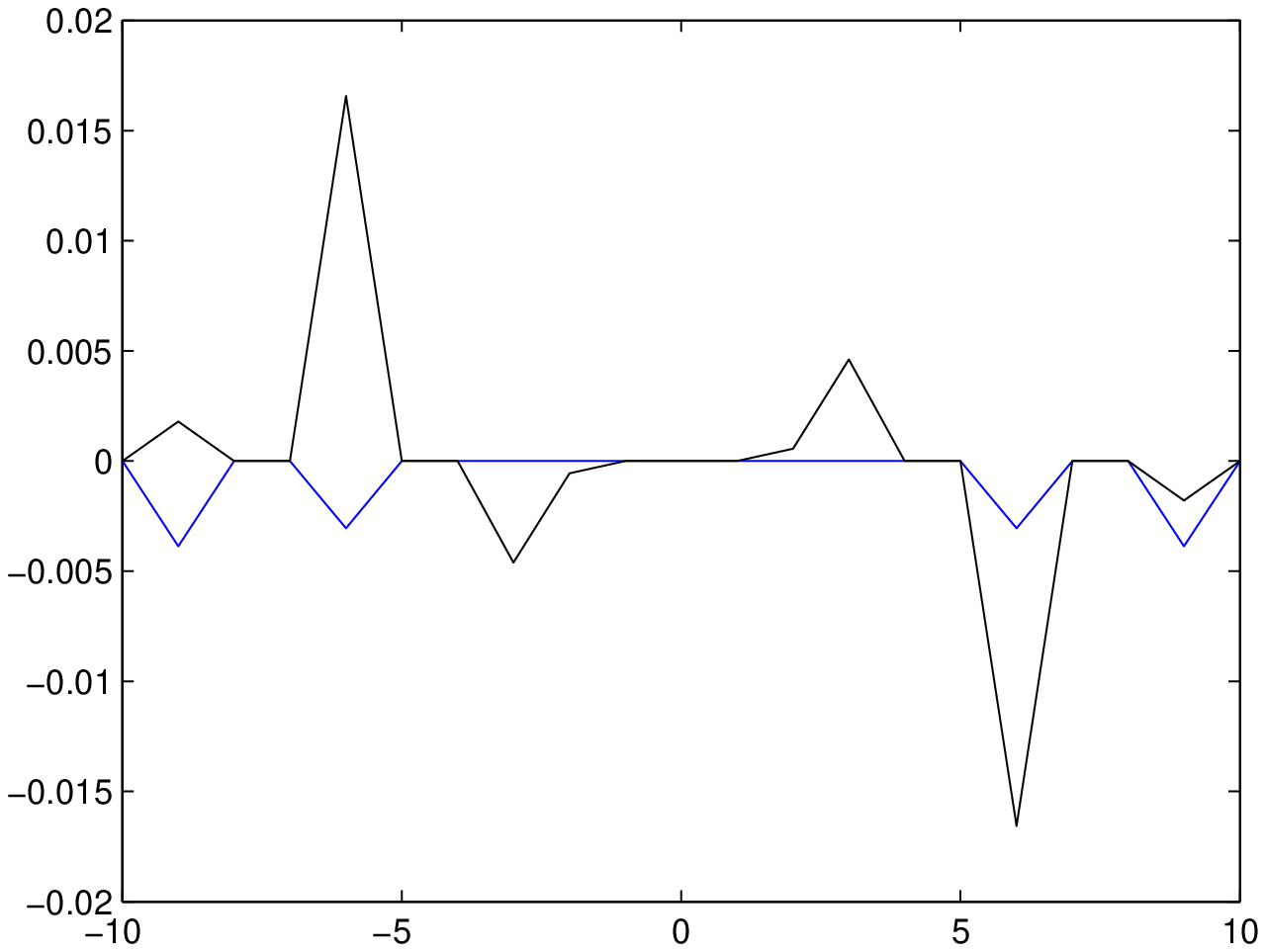}
\caption{Fourier coefficients of reconstructed perturbations in Example 1; \textbf{Set 1} to \textbf{Set 3}
from left to right;  blue: real part; black: imaginary part.} \label{test1B}
\end{figurehere}
Although there are some deficiencies in the reconstruction of Fourier modes, we can see from these figures the Fourier coefficients reconstructed from \textbf{Set 1} is largest at $|n| = 3$ with its magnitude almost between $0.04$ and $0.05$, which clearly indicates a strong dominance of $\delta \cos (3 \theta)$ with magnitude $\delta$ between $0.08$ and $0.1$ and corresponds to the signal from the exact inclusion.
The reconstruction from \textbf{Set 2} has more deficiency, that the Fourier mode is somehow shifted to $\delta ( \cos (3 \theta) + \sin (3 \theta) ) $ with magnitude $\delta$ between $0.02$ and $0.025$. However, the location of the peak Fourier mode is still correct.
Nonetheless, the reconstruction from \textbf{Set 3} deviates totally from the exact solution, indicating its insufficiency in number of measurements to reconstruct the perturbation. This goes with the theoretical analysis in section \ref{sec71}.

Now we show in Figure \ref{test1C} (top) the shapes of reconstructed domains, Figure \ref{test1C} (middle)  the contrast of the reconstructed media and Figure \ref{test1C} (right)
a comparison between the reconstructed domains $D^{\text{approx}}$ and exact domain $D$ using the values of a sum of characteristic functions $\chi_{D} + \chi_{D^{\text{approx}}}$.
The relative $L^2$ errors of the reconstructions for \textbf{Set 1} to \textbf{Set 3} are respectively $3.29 \%$, $6.34 \%$ and $11.72 \%$.  In view of the severe ill-posedness of the phaseless reconstruction problem and $5\%$ percent of measurement noise, the reconstructions from \textbf{Set 1} and \textbf{Set 2} measurements are quite reasonable

\begin{figurehere} \centering
\includegraphics[width=5cm,height=4cm]{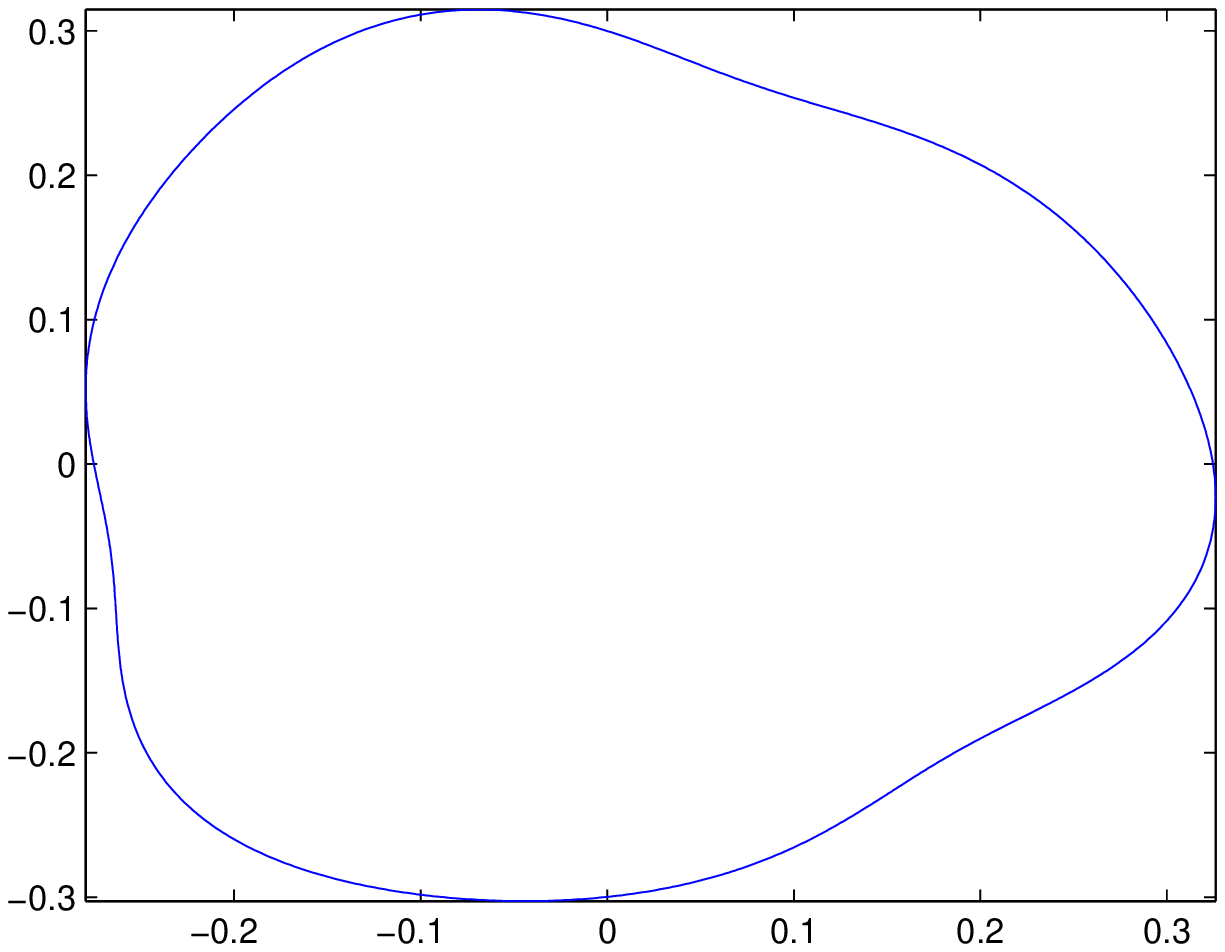}
\includegraphics[width=5cm,height=4cm]{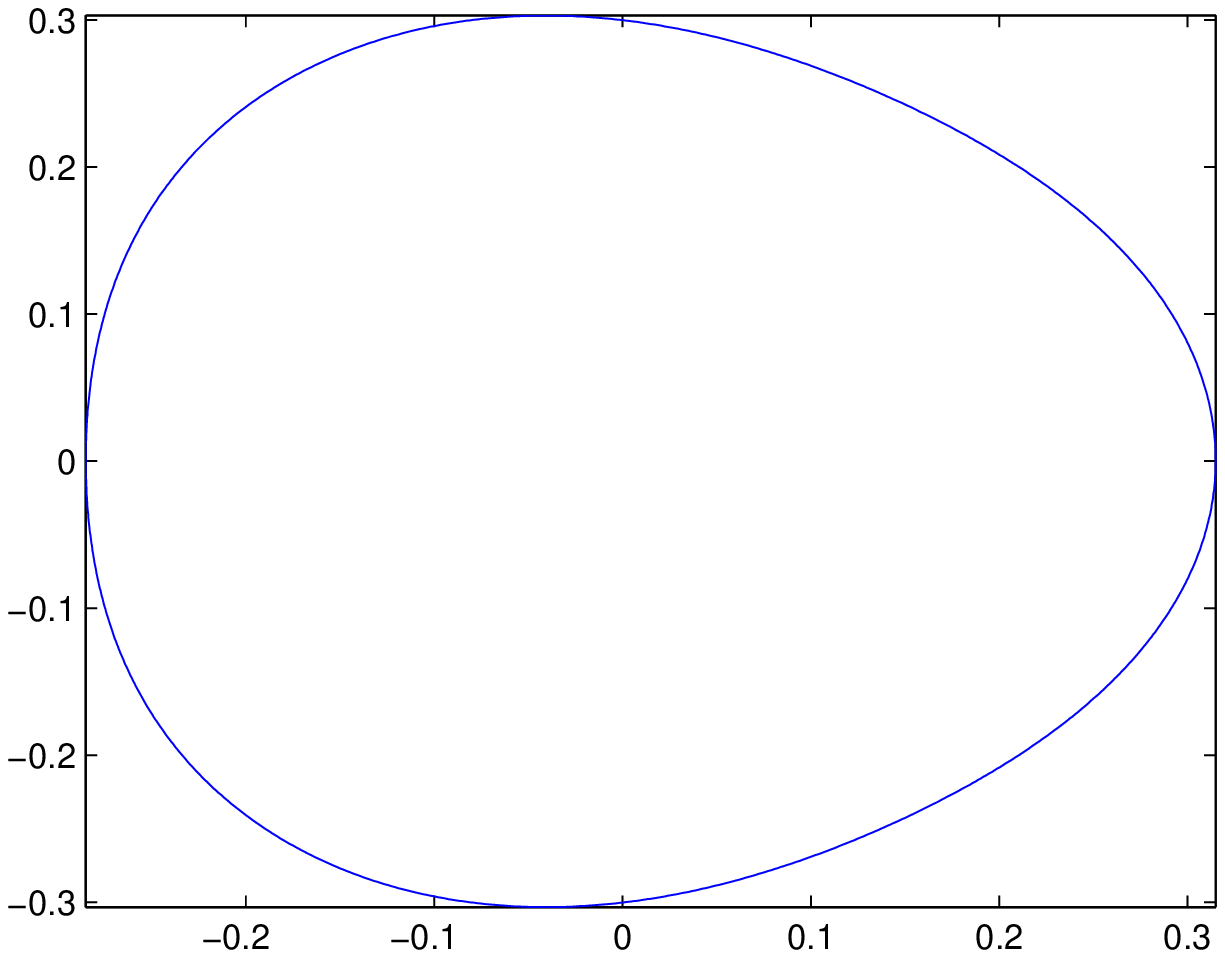}
\includegraphics[width=5cm,height=4cm]{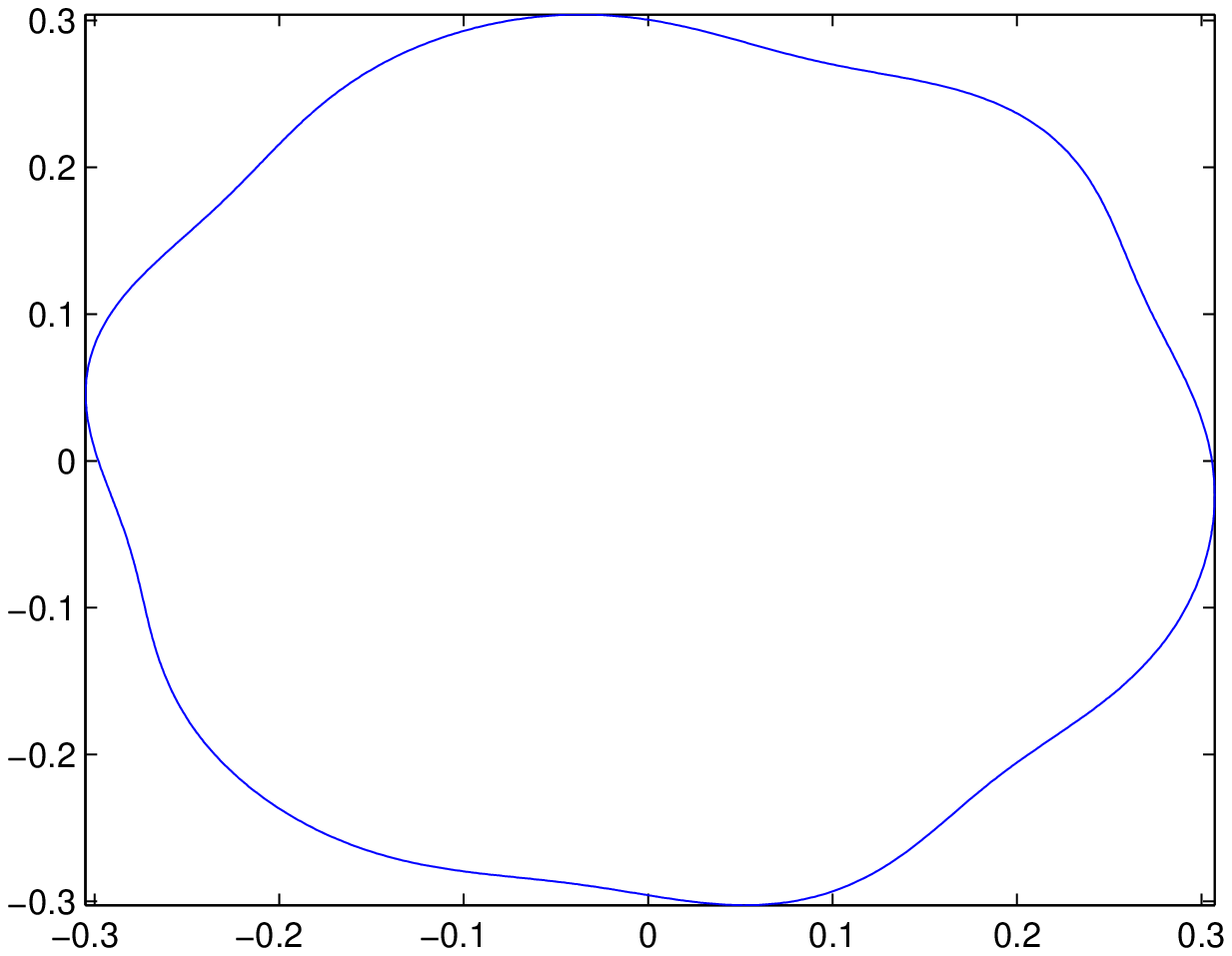} \\
\includegraphics[width=5cm,height=4cm]{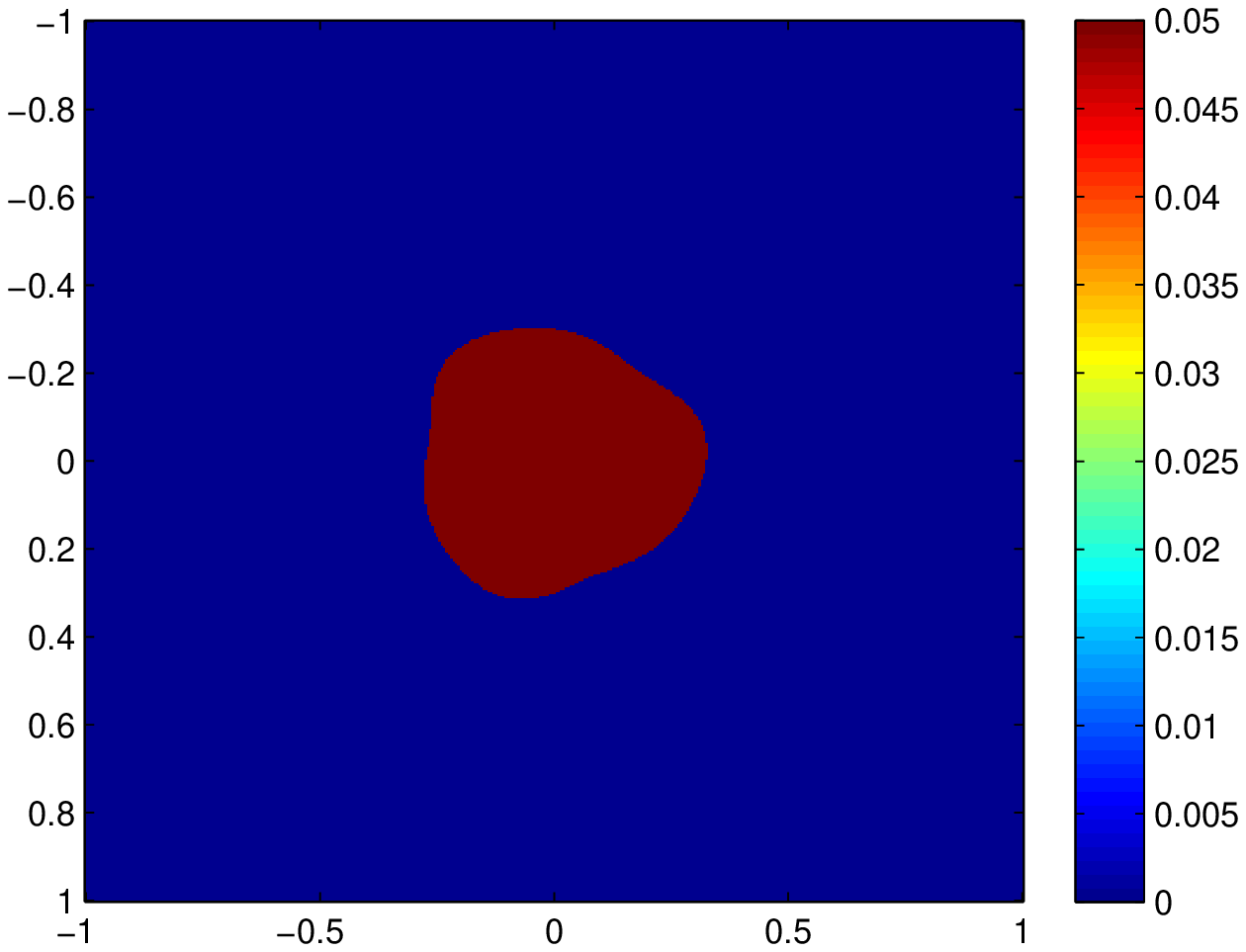}
\includegraphics[width=5cm,height=4cm]{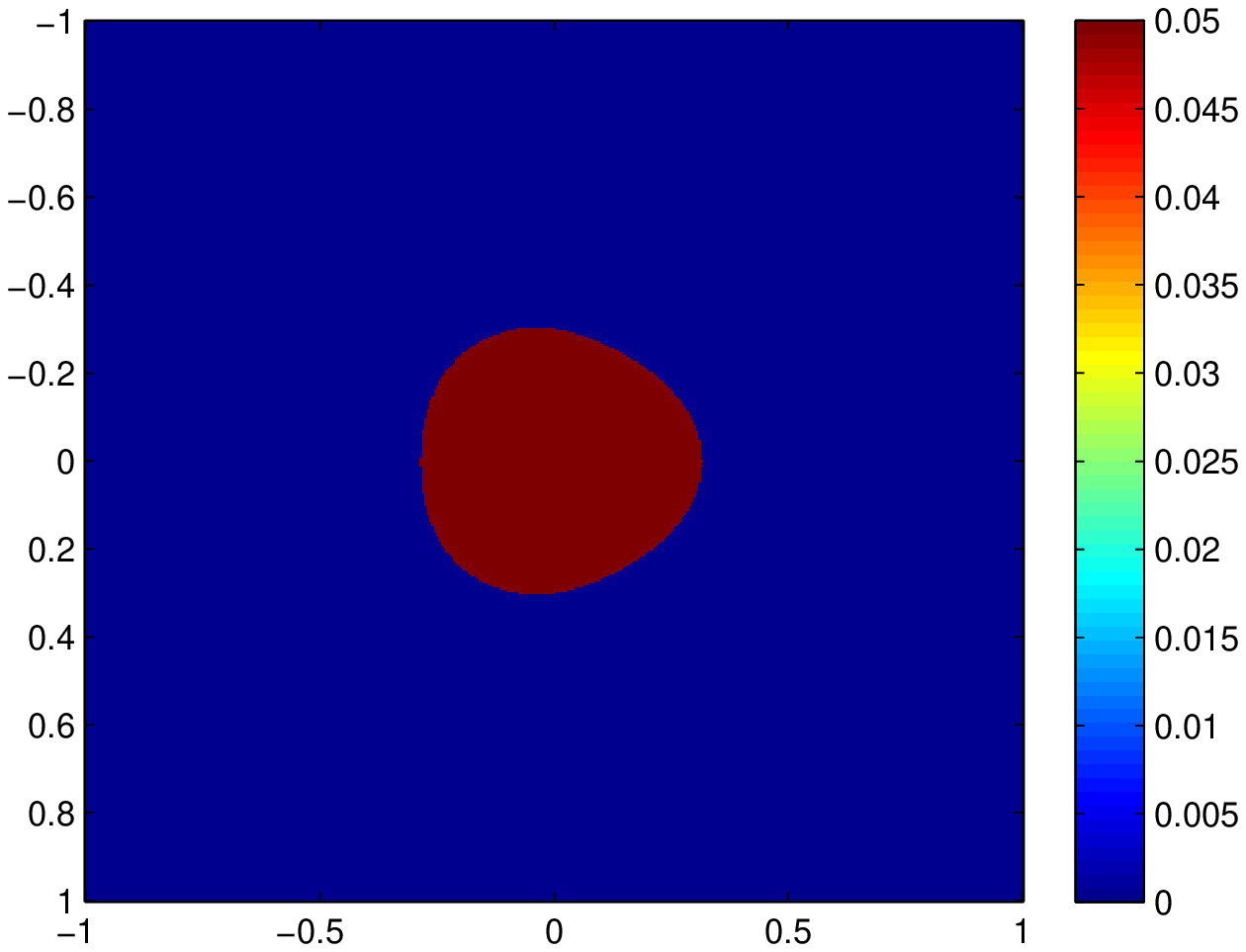}
\includegraphics[width=5cm,height=4cm]{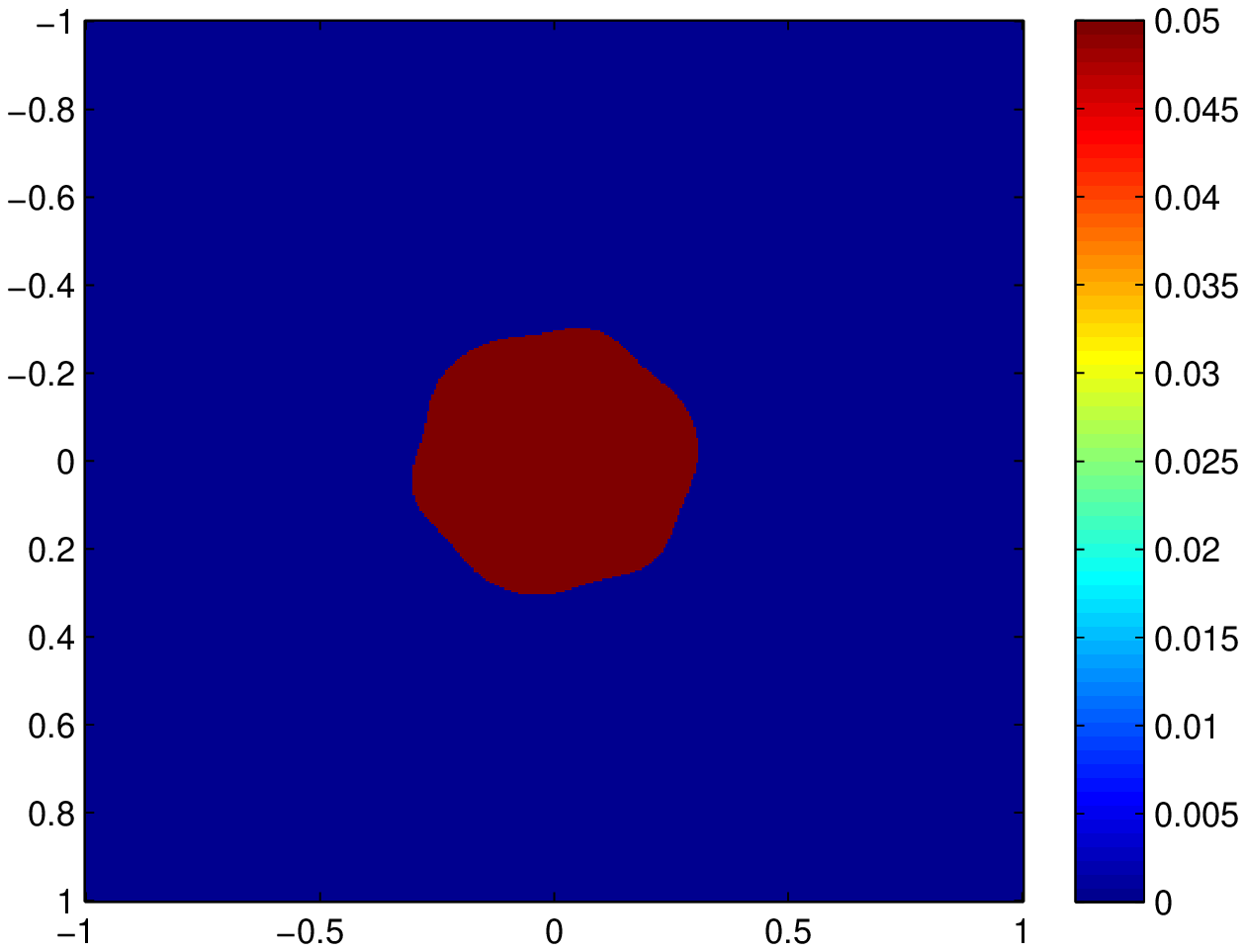} \\
\includegraphics[width=5cm,height=4cm]{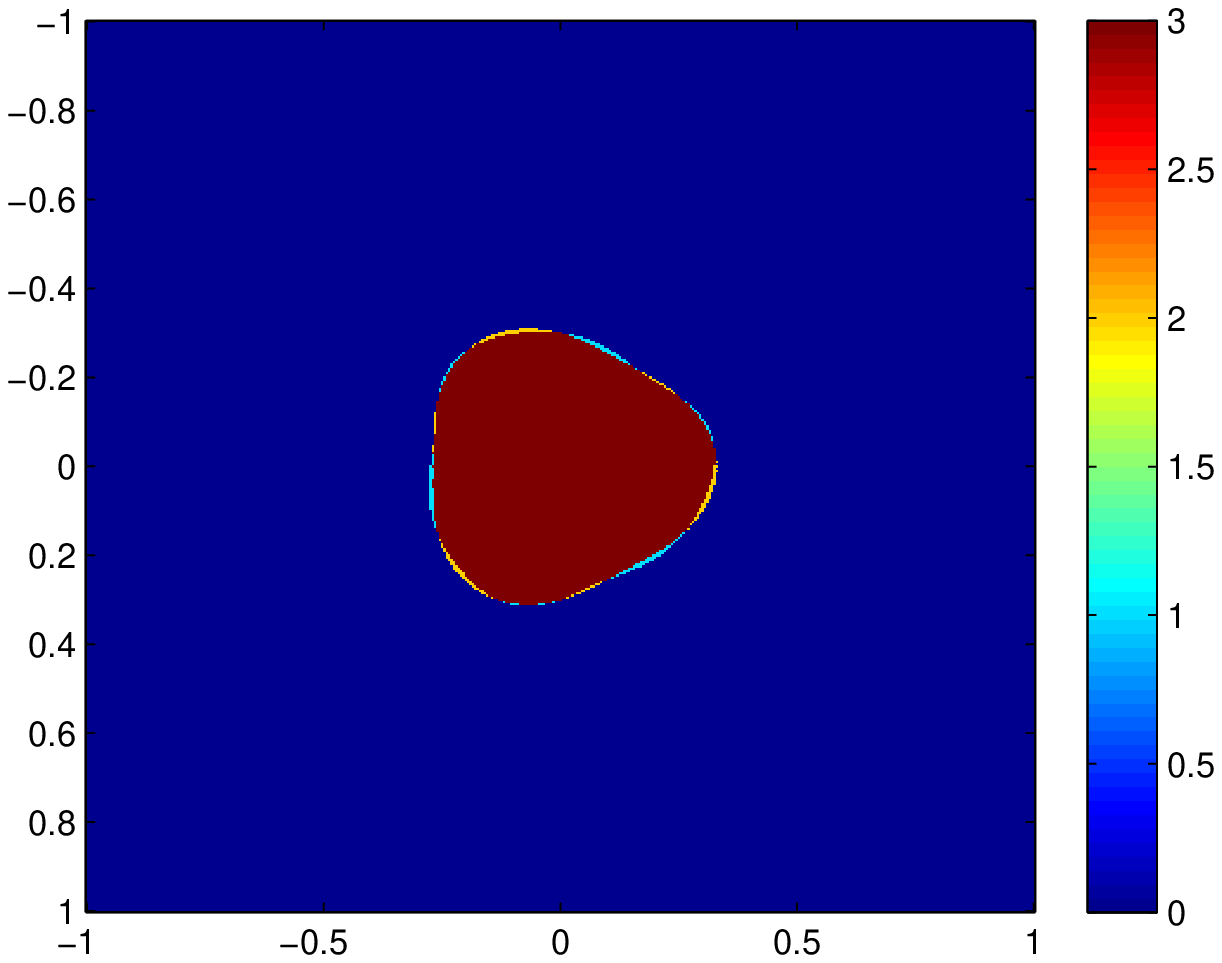}
\includegraphics[width=5cm,height=4cm]{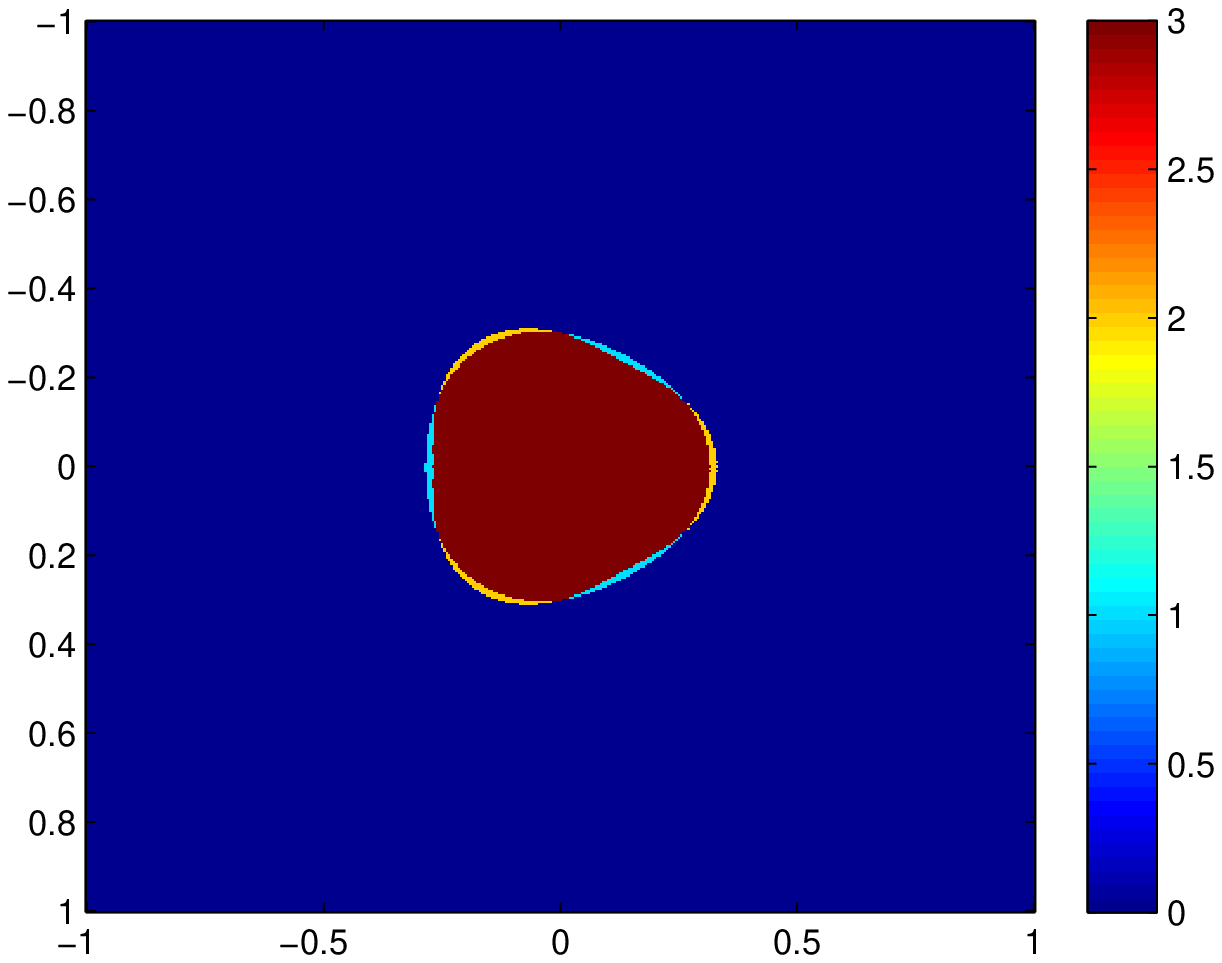}
\includegraphics[width=5cm,height=4cm]{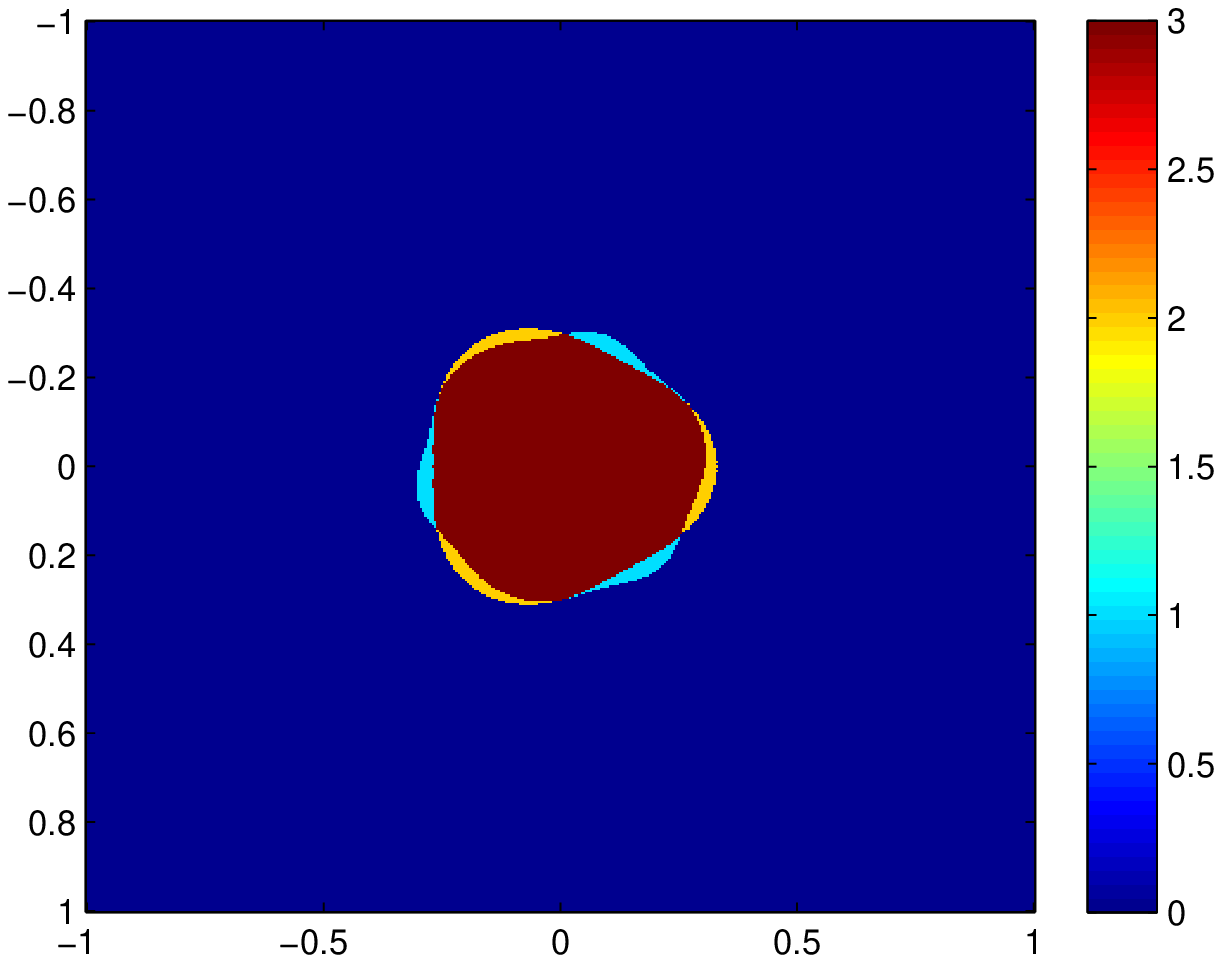} \\
\caption{Reconstructed domain and medium in Example 1 and comparison between the exact and reconstructed domains; 
\textbf{Set 1} to \textbf{Set 3} from left to right; from top to bottom: 
reconstructed shape, reconstructed inclusion, and comparison between reconstructed and exact domains.}  \label{test1C}
\end{figurehere}

\textbf{Example 2}.
We test another domain of the flori-form shape described by \eqref{circle_perturb_form} with $\delta = 0.1$ and $n = 5$.
Figure \ref{test2A} (left) and (right) show the shape of the domain and the contrast of the inhomogeneous medium
respectively.

\begin{figurehere} \centering
\includegraphics[width=5cm,height=4cm]{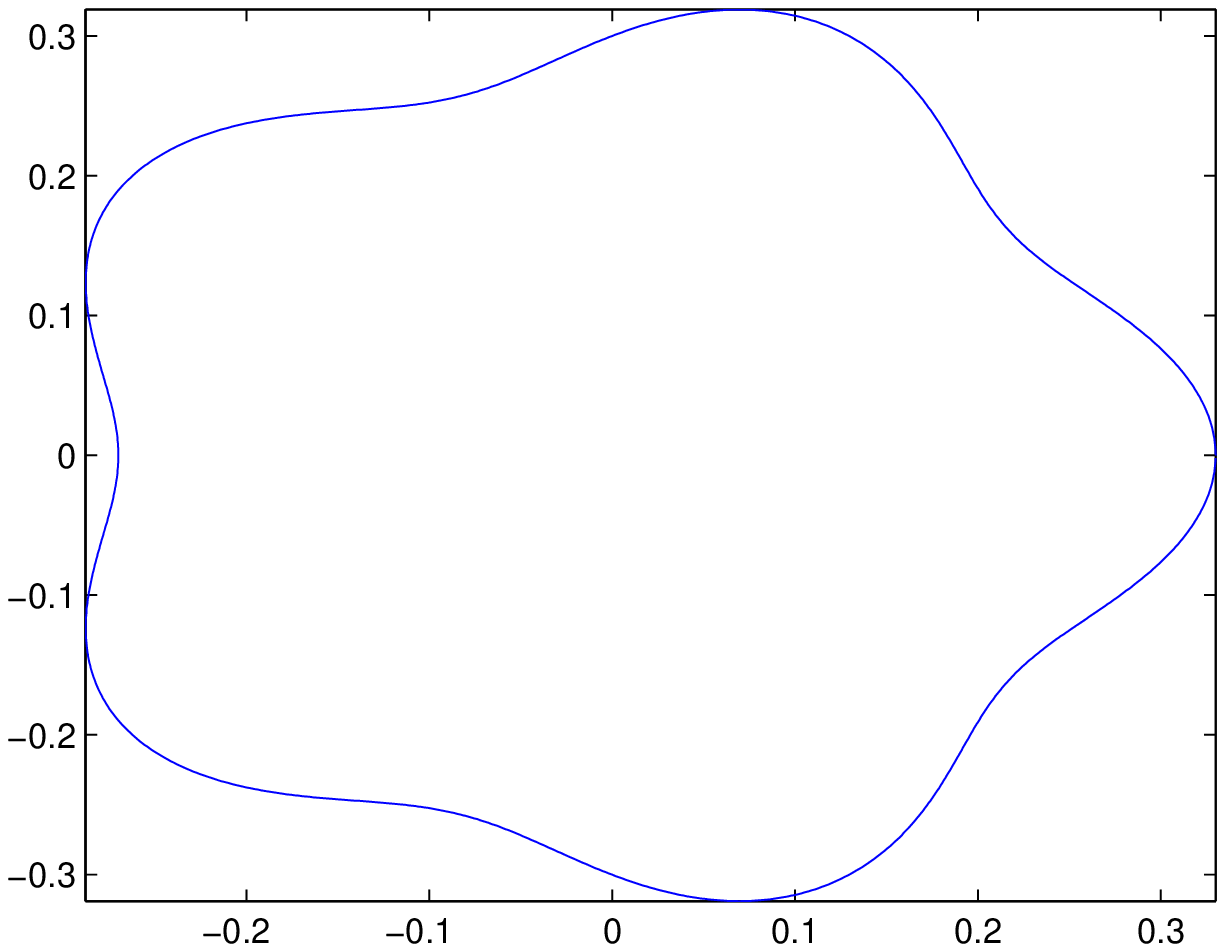}
\includegraphics[width=5cm,height=4cm]{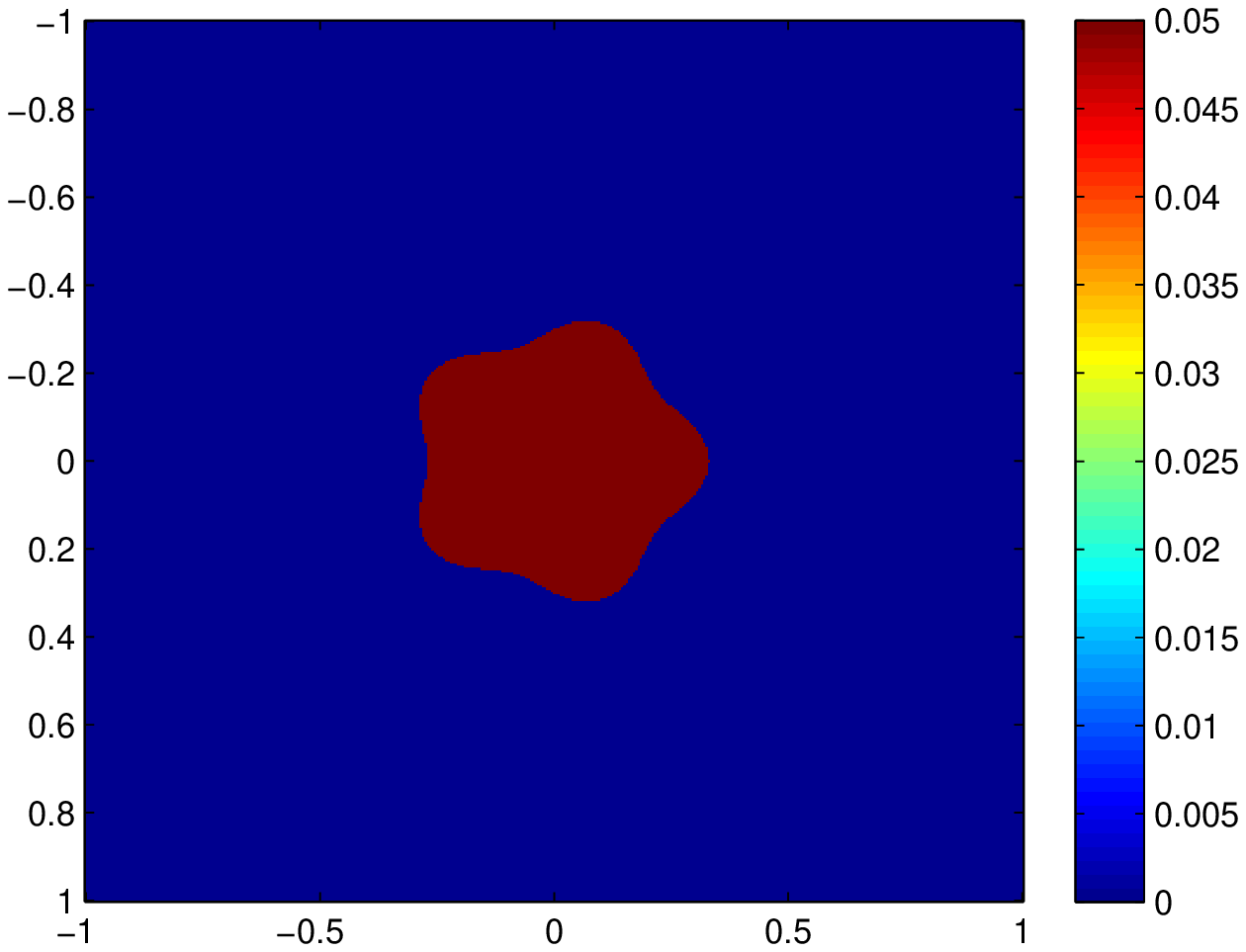}
\caption{Exact inhomogeneous domain (left) and the contrast of the inclusion in Example 2.} \label{test2A}
\end{figurehere}

In this example, the magnitude of the far-field pattern for $16$ wave-numbers are used for shape reconstruction, i.e. $\tilde{C} = 15$.
The Fourier coefficients of the reconstructed perturbations using the respective measurement sets are 
shown in Figure \ref{test2B}.

\begin{figurehere} \centering
\includegraphics[width=5cm,height=4cm]{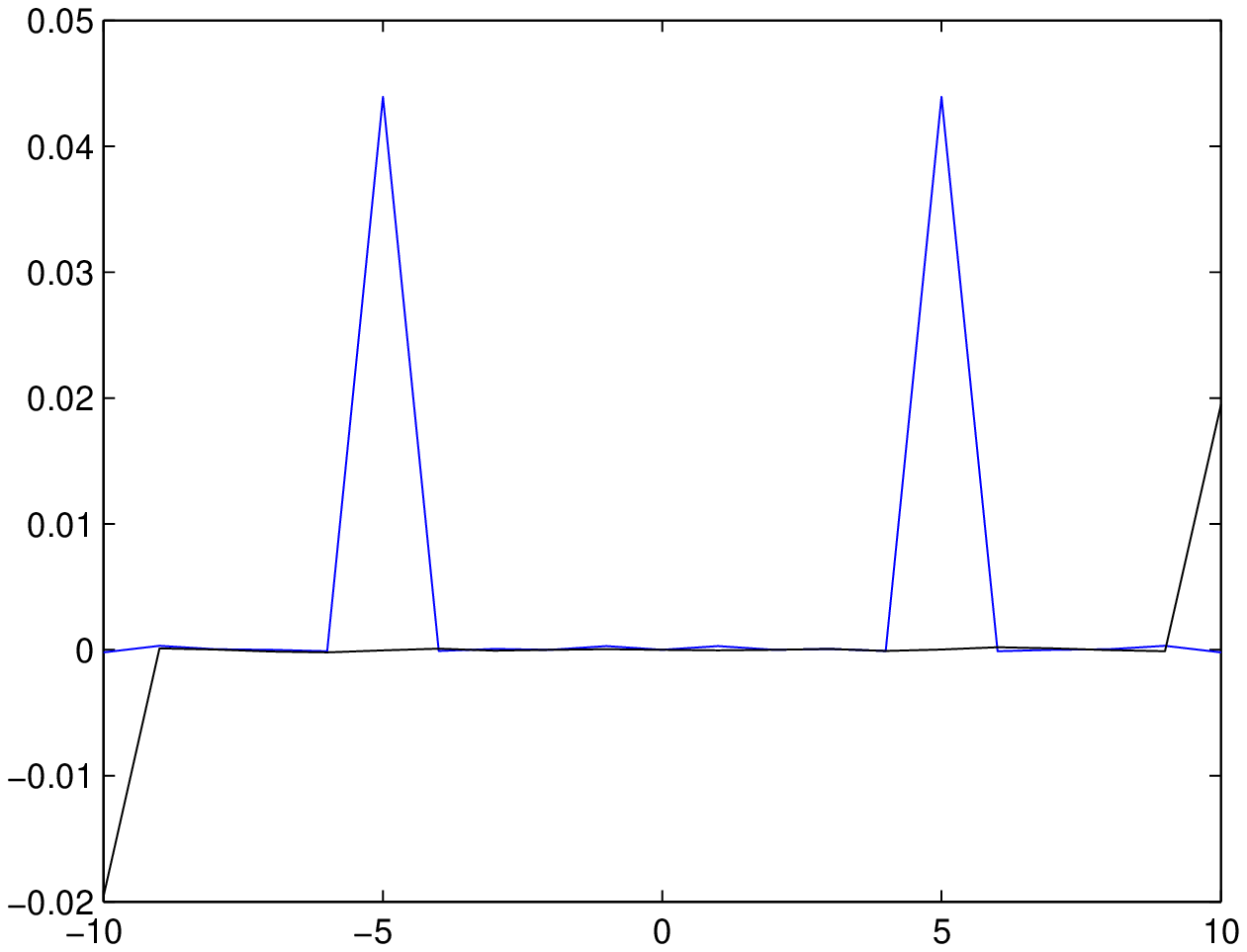}
\includegraphics[width=5cm,height=4cm]{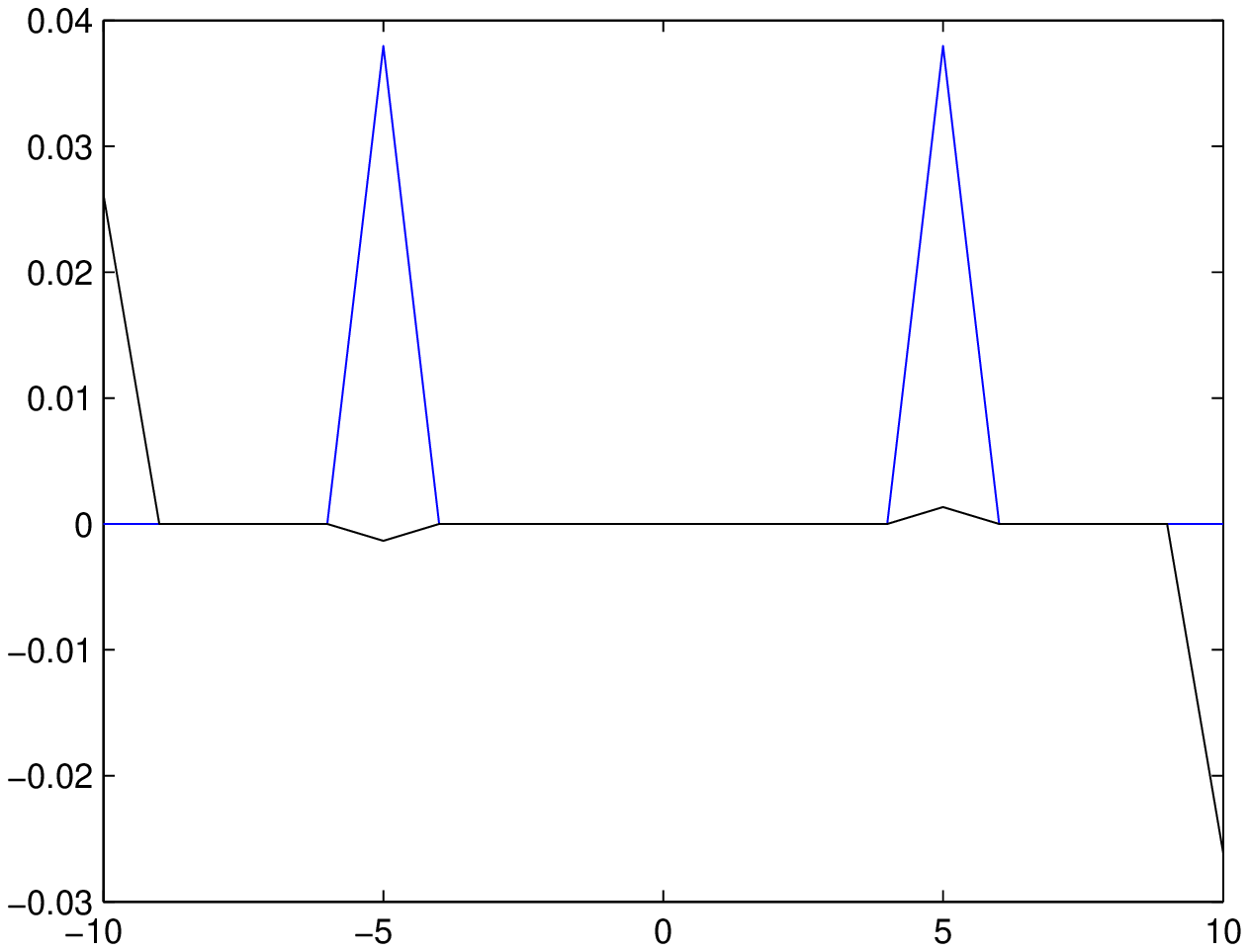}
\includegraphics[width=5cm,height=4cm]{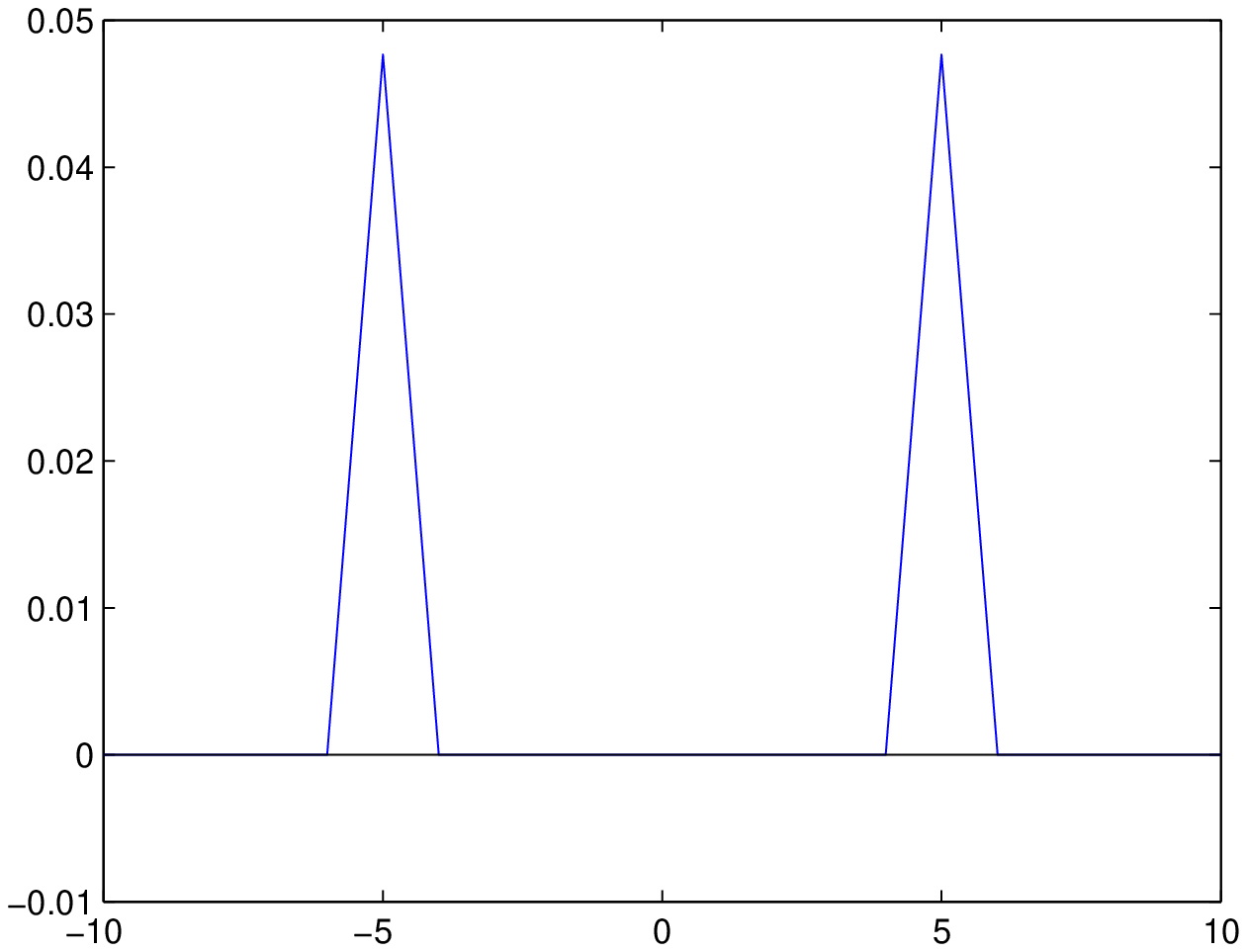}
\caption{Fourier coefficients of reconstructed perturbations in Example 2; 
\textbf{Set 1} to \textbf{Set 3} from left to right;   blue: real part; black: imaginary part.} \label{test2B}
\end{figurehere}

We can now see that both the reconstructions from \textbf{Set 1} and \textbf{Set 2} are reasonable and indicate the correct peak Fourier modes and its magnitude. It is no surprise to see 
the reconstruction for \textbf{Set 2} is worse than that for \textbf{Set 1}.
However we can see that in this particular case, the reconstruction for \textbf{Set 3} coincidentally collides with the exact solution after regularization.

In Figure \ref{test2C} (top), (middle) and (right), the shapes of reconstructed domains, 
the contrast of the reconstructed media and 
the comparison between the reconstructed domains $D^{\text{approx}}$ and exact domain $D$ 
(by showing a sum of characteristic functions $\chi_{D} + \chi_{D^{\text{approx}}}$) are presented respectively.
The relative $L^2$ errors of the reconstructions for \textbf{Set 1} to \textbf{Set 3} are respectively $4.91 \%$,
$6.94 \%$
and
$0.57\%$.  As we mentioned above, quite surprisingly, the $L^1$ regularizer coincidentally provides a very good estimate 
for \textbf{Set 3}.

\begin{figurehere} \centering
\includegraphics[width=5cm,height=4cm]{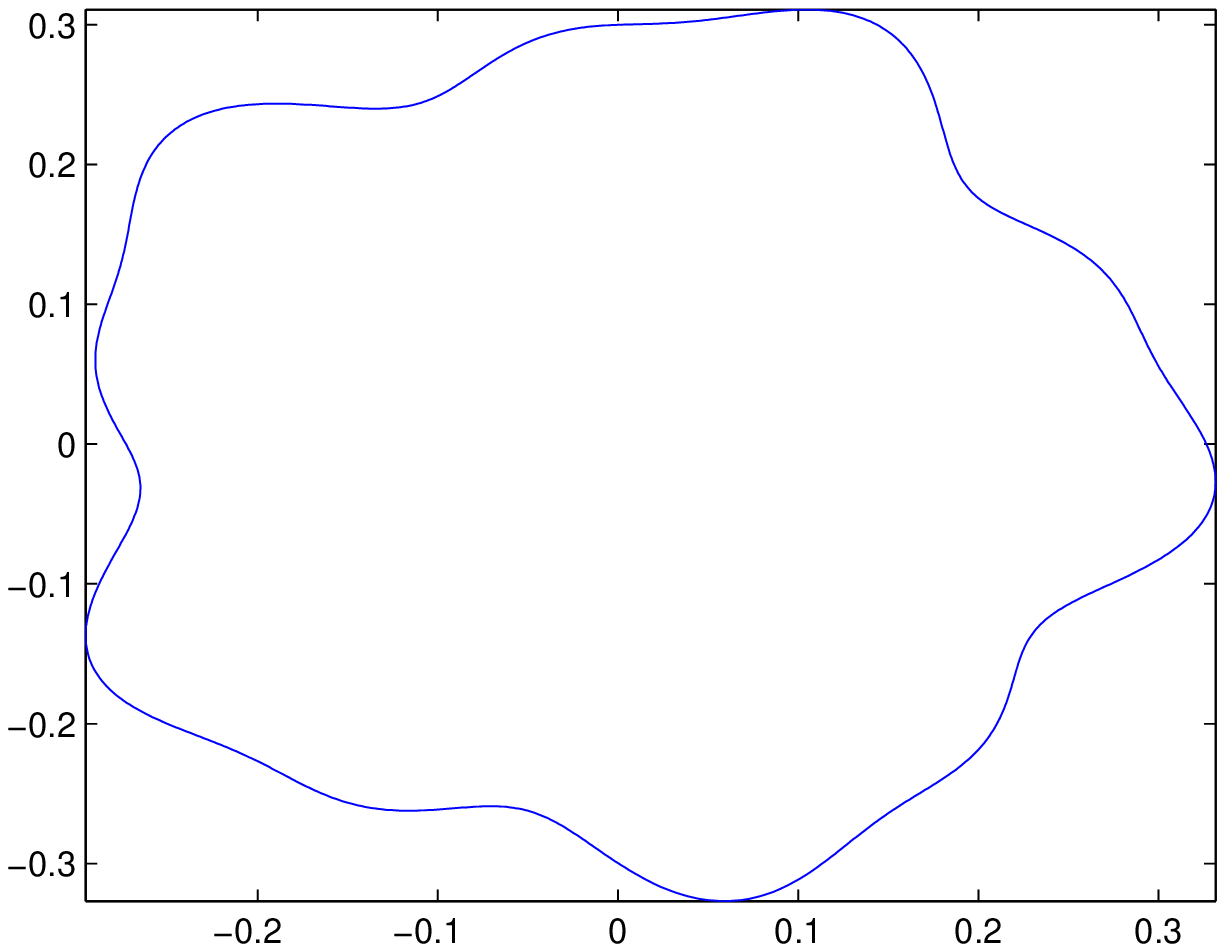}
\includegraphics[width=5cm,height=4cm]{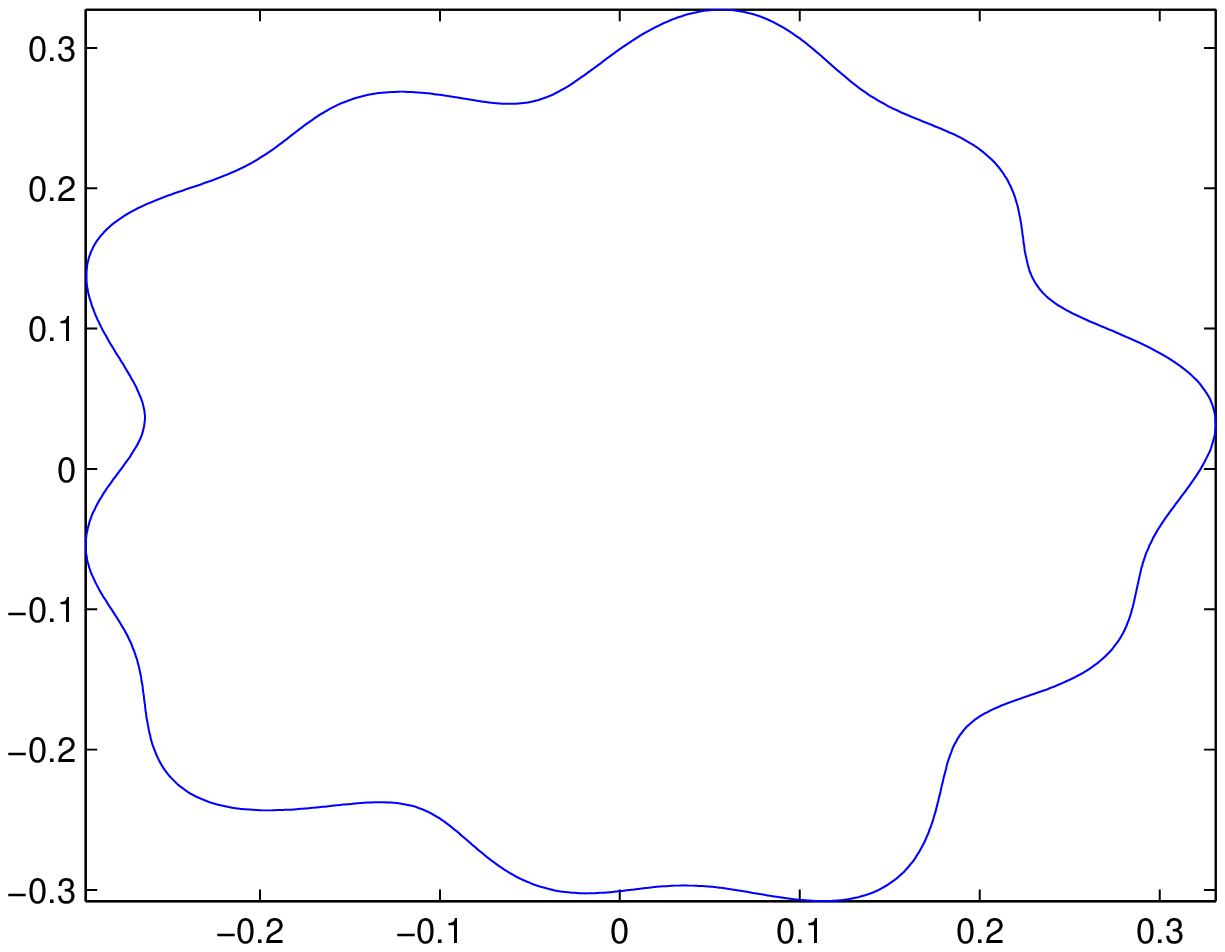}
\includegraphics[width=5cm,height=4cm]{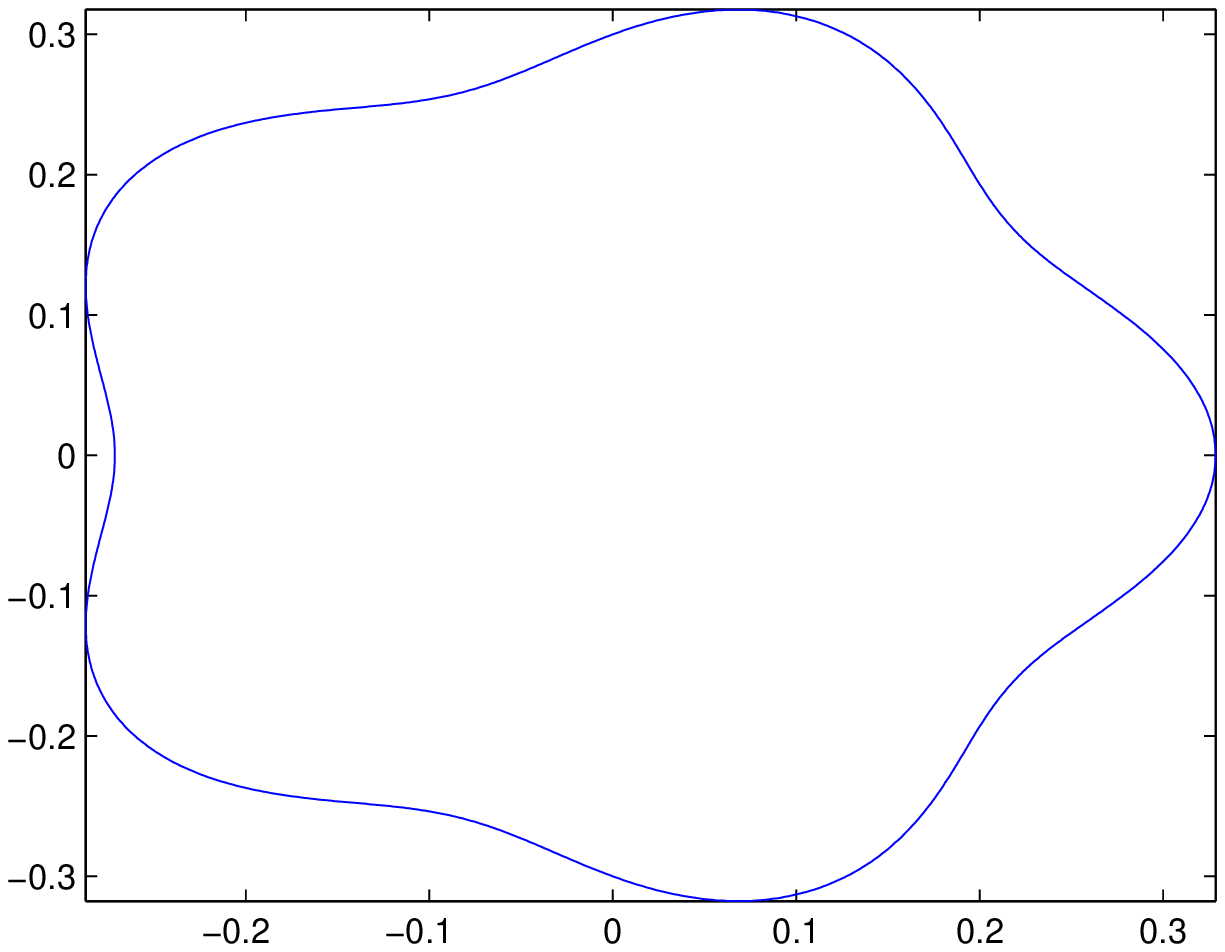} \\
\includegraphics[width=5cm,height=4cm]{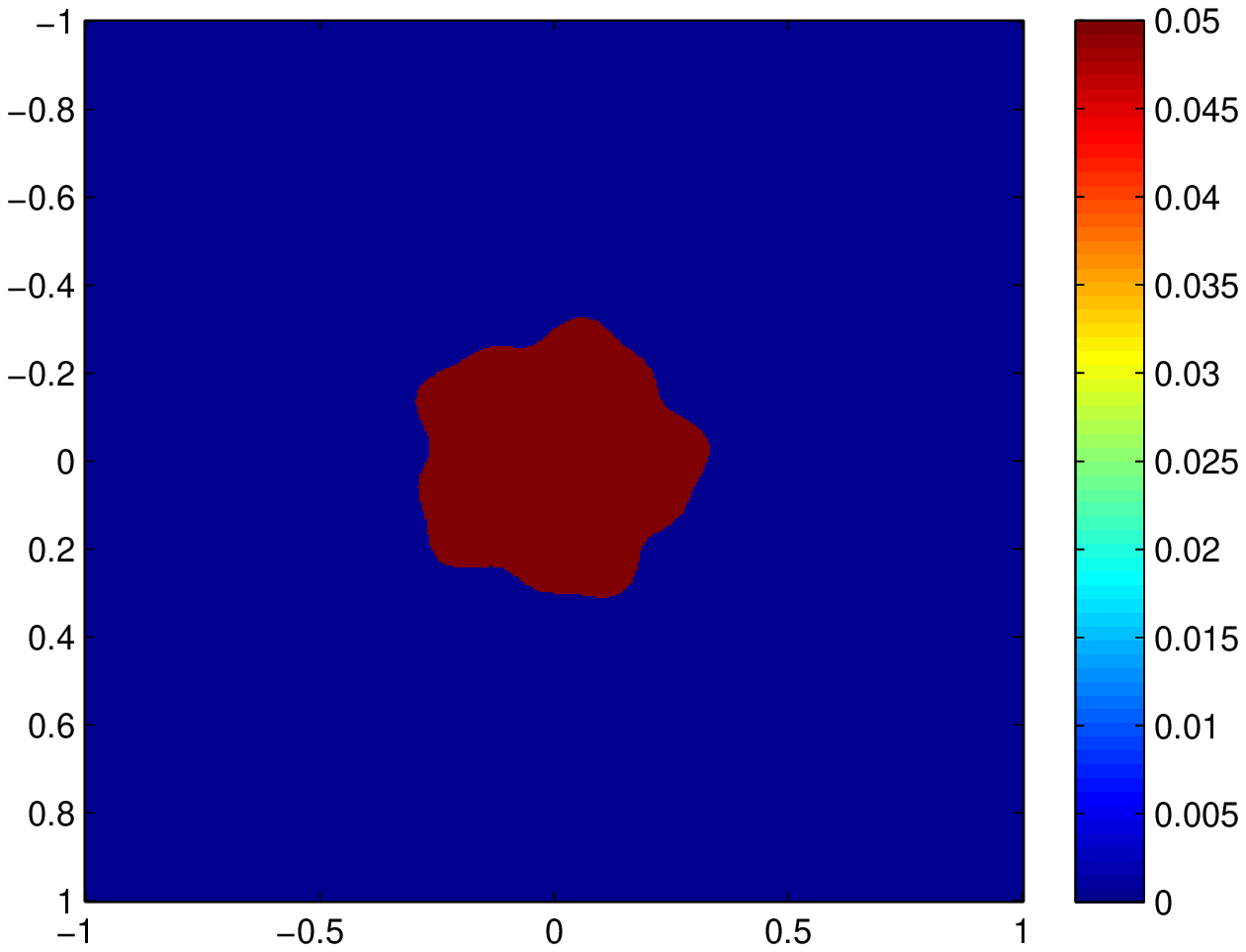}
\includegraphics[width=5cm,height=4cm]{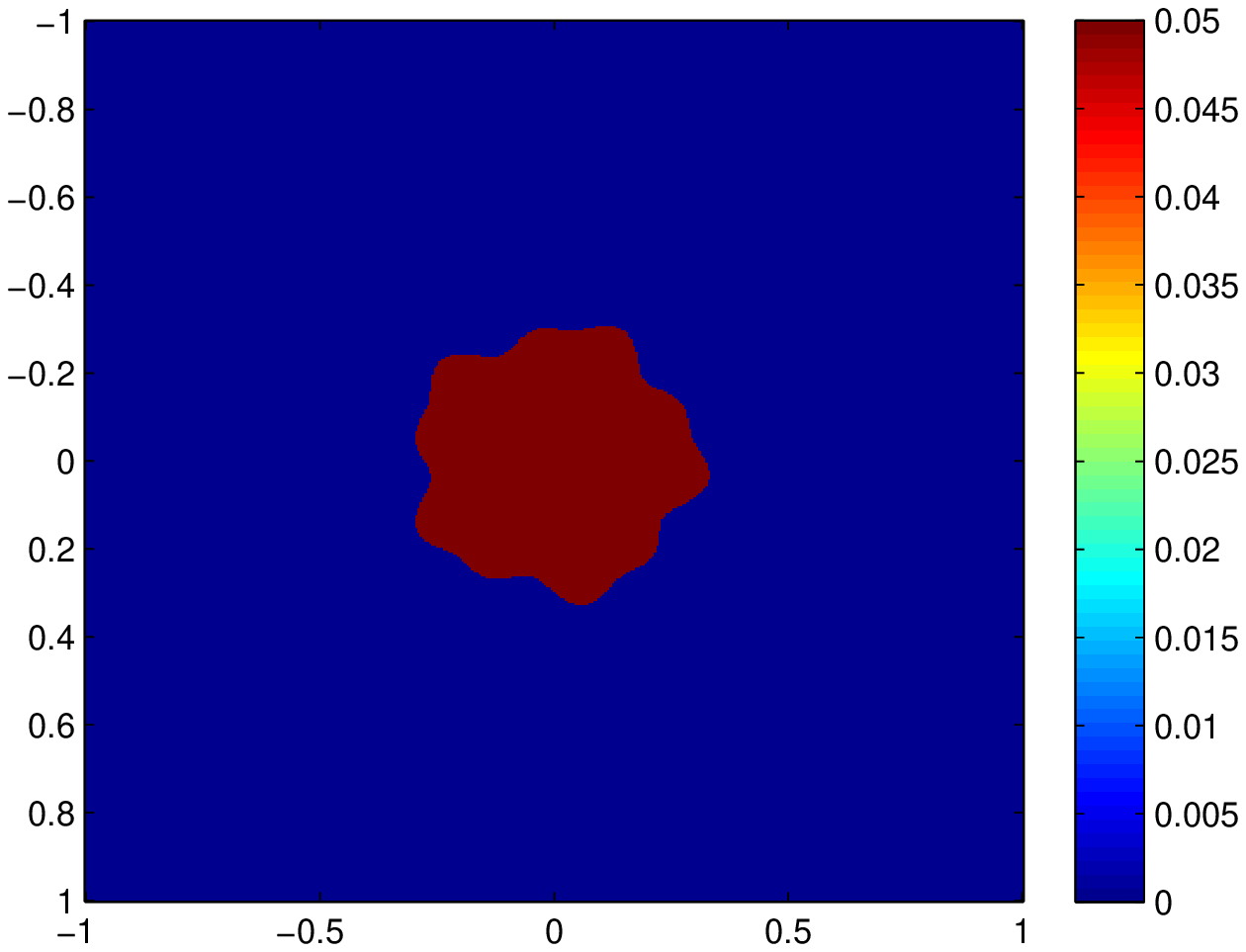}
\includegraphics[width=5cm,height=4cm]{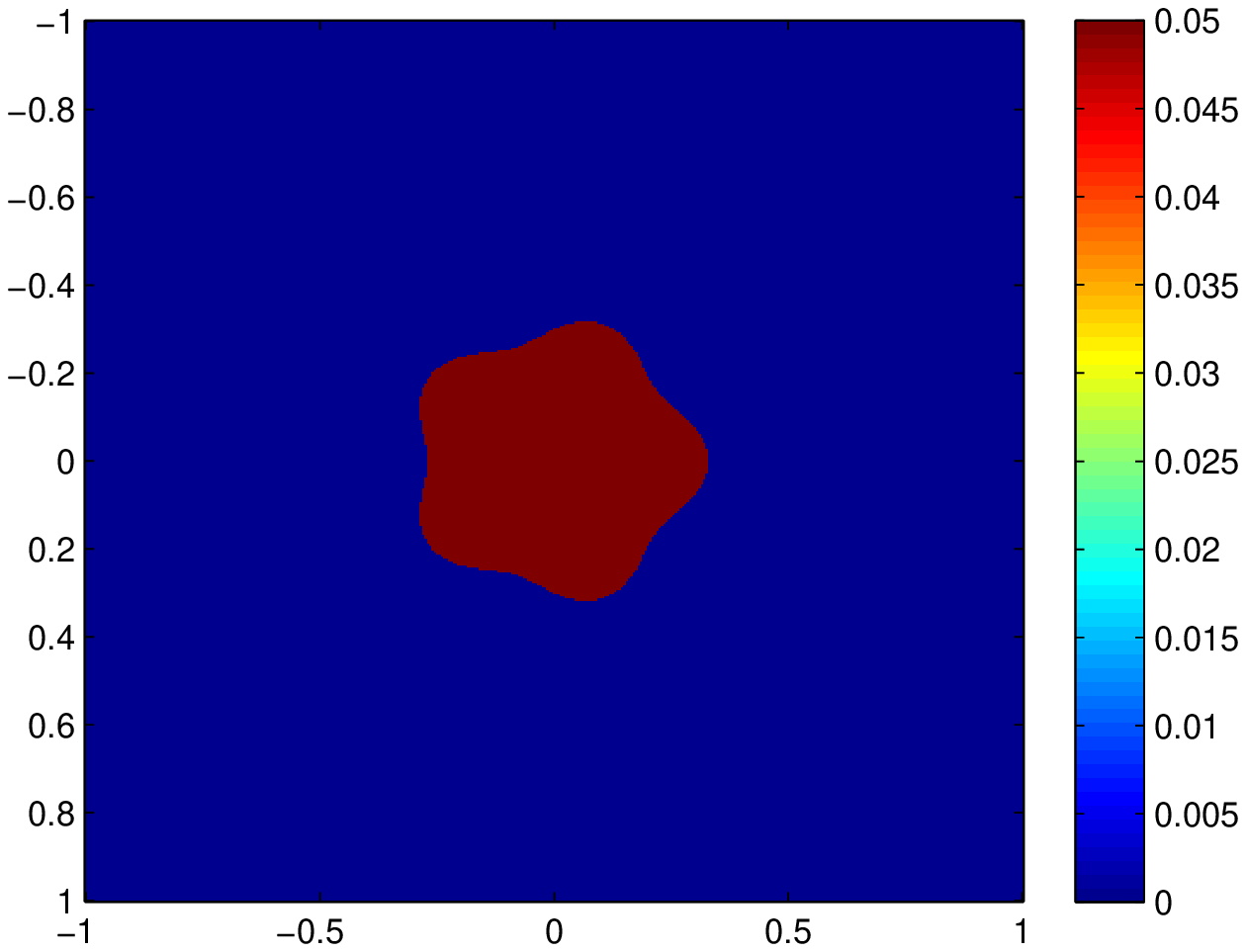} \\
\includegraphics[width=5cm,height=4cm]{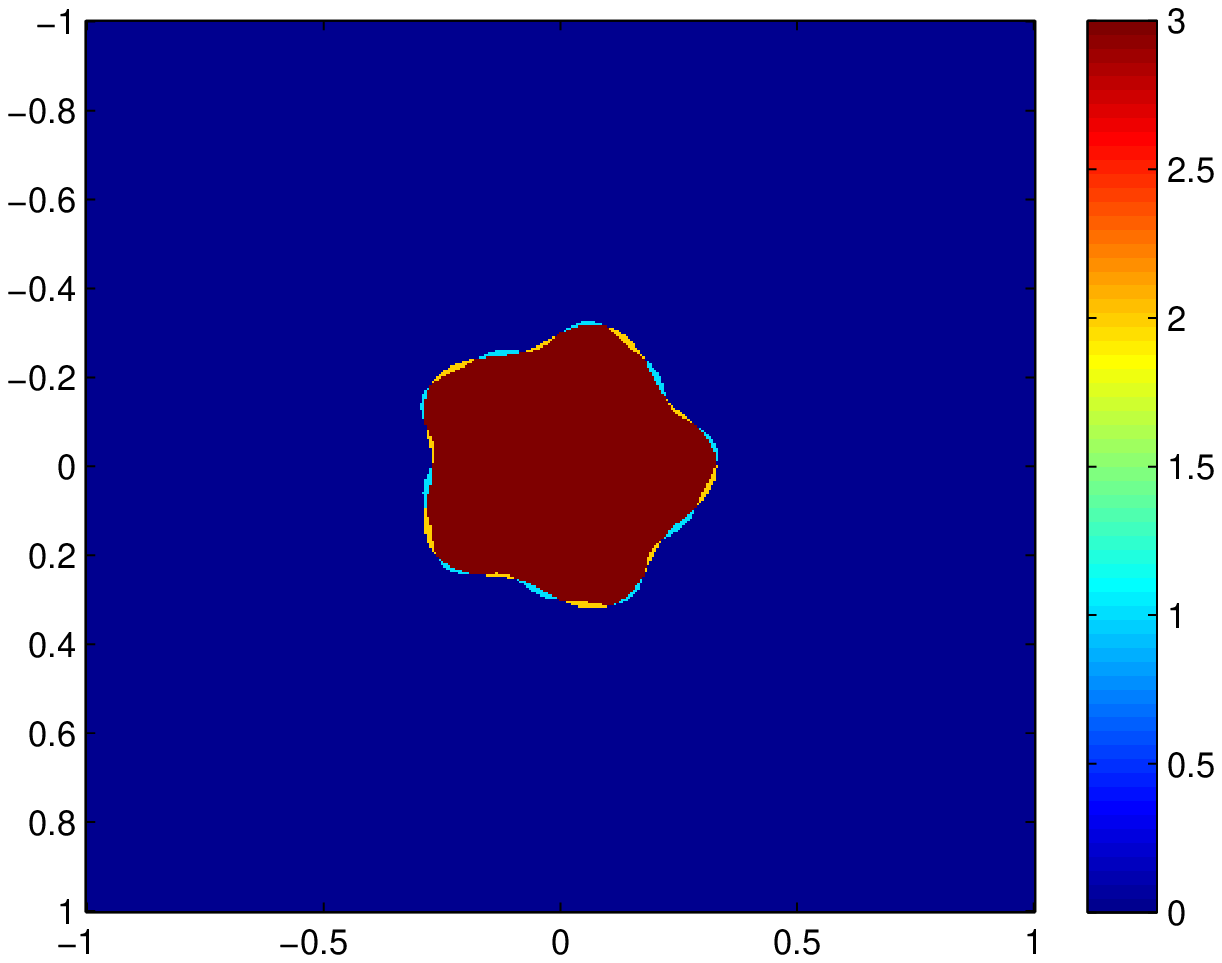}
\includegraphics[width=5cm,height=4cm]{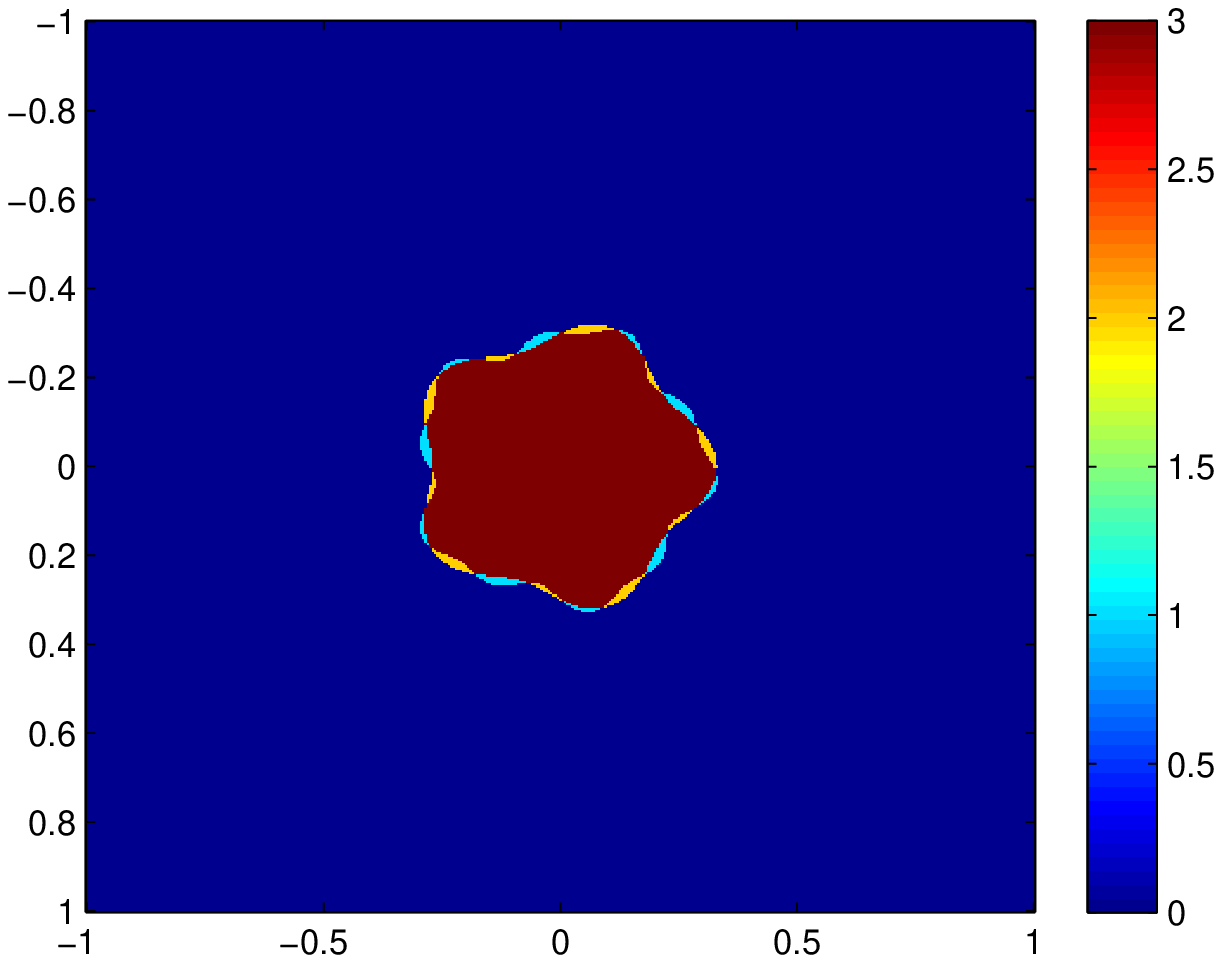}
\includegraphics[width=5cm,height=4cm]{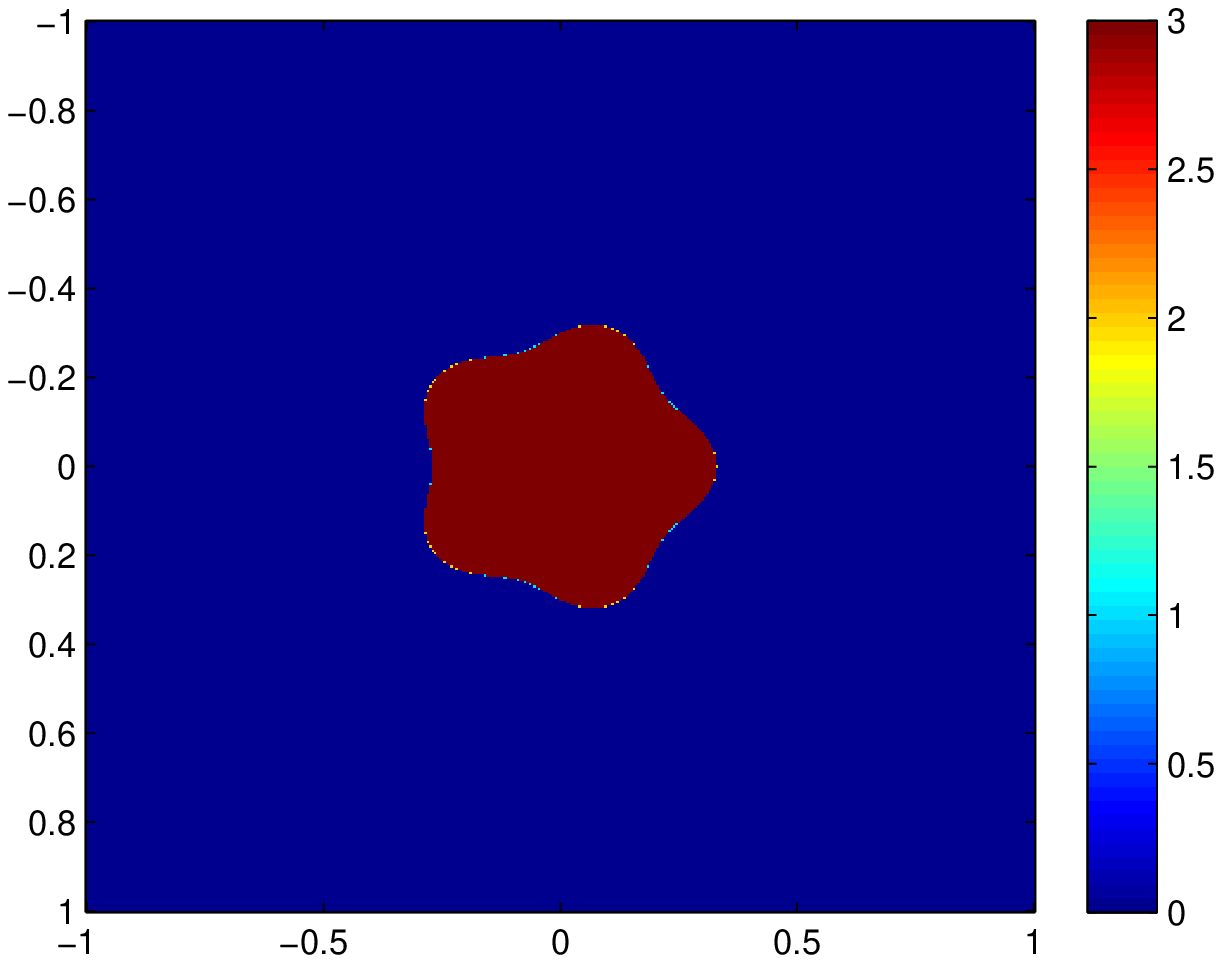} \\
\caption{Reconstructed domain and medium in Example 2 and comparison between the exact and reconstructed domains.  \textbf{Set 1} to \textbf{Set 3} from left to right; 
reconstructed shape, reconstructed inclusion and 
comparison between reconstructed and exact domains from top to bottom.}  \label{test2C}
\end{figurehere}

\textbf{Example 3}.
In this last example, we test a domain of more complicated flori-form shape $D = B^\delta$ described by the following parametric form (with $\delta = 0.1$ and $n = 3$ ):
    \begin{equation}
        r = 0.2 ( 1 + \delta \cos (n \theta) + 2 \delta \cos (2 n \theta) )\, , \quad \theta \in (0,2\pi]\,,
        \label{circle_perturb_form_2}
    \end{equation}
The shape of the domain is given in Figure \ref{test4A} (left) and the contrast of the inhomogeneous medium in Figure \ref{test4A} (right).

\begin{figurehere} \centering
\includegraphics[width=5cm,height=4cm]{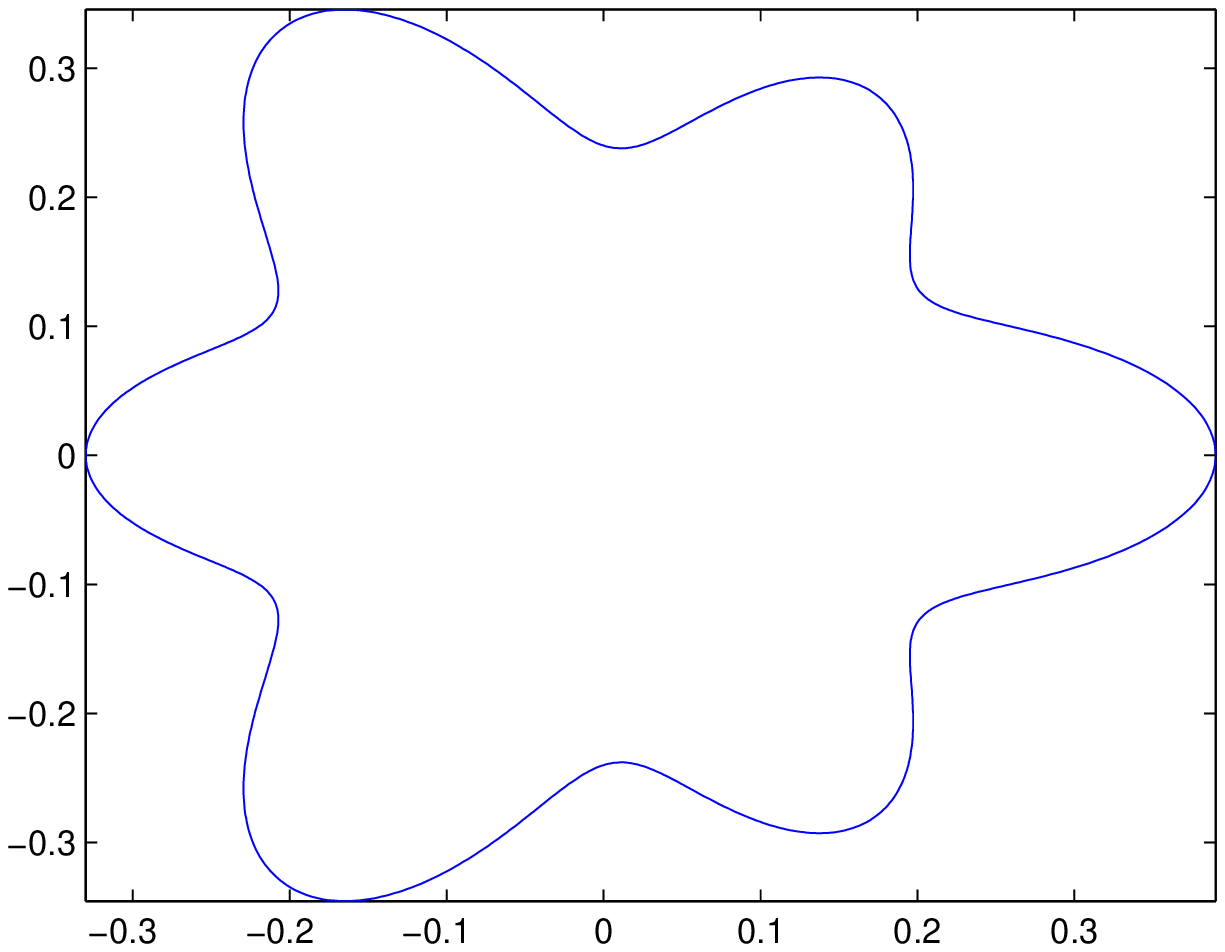}
\includegraphics[width=5cm,height=4cm]{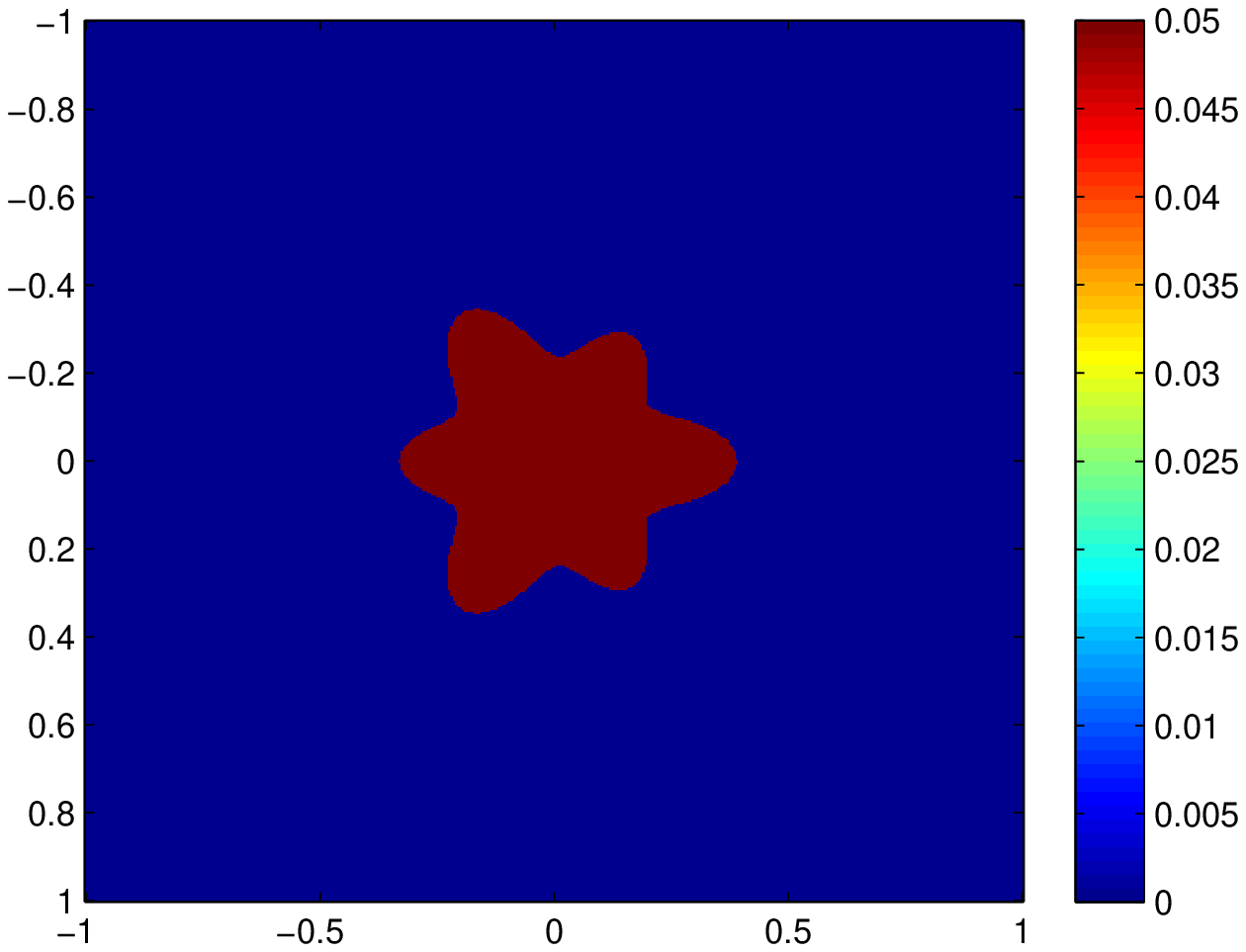}
\caption{Exact inhomogeneous domain (left) and contrast of the inclusion (right) in Example 3.} \label{test4A}
\end{figurehere}

The magnitude of far-field pattern for $6$ wave-numbers are used for shape reconstruction, i.e. $\tilde{C} = 5$, and 
the Fourier coefficients of the reconstructed perturbations using the respective measurement sets are 
shown in Figure \ref{test4B}.

\begin{figurehere} \centering
\includegraphics[width=5cm,height=4cm]{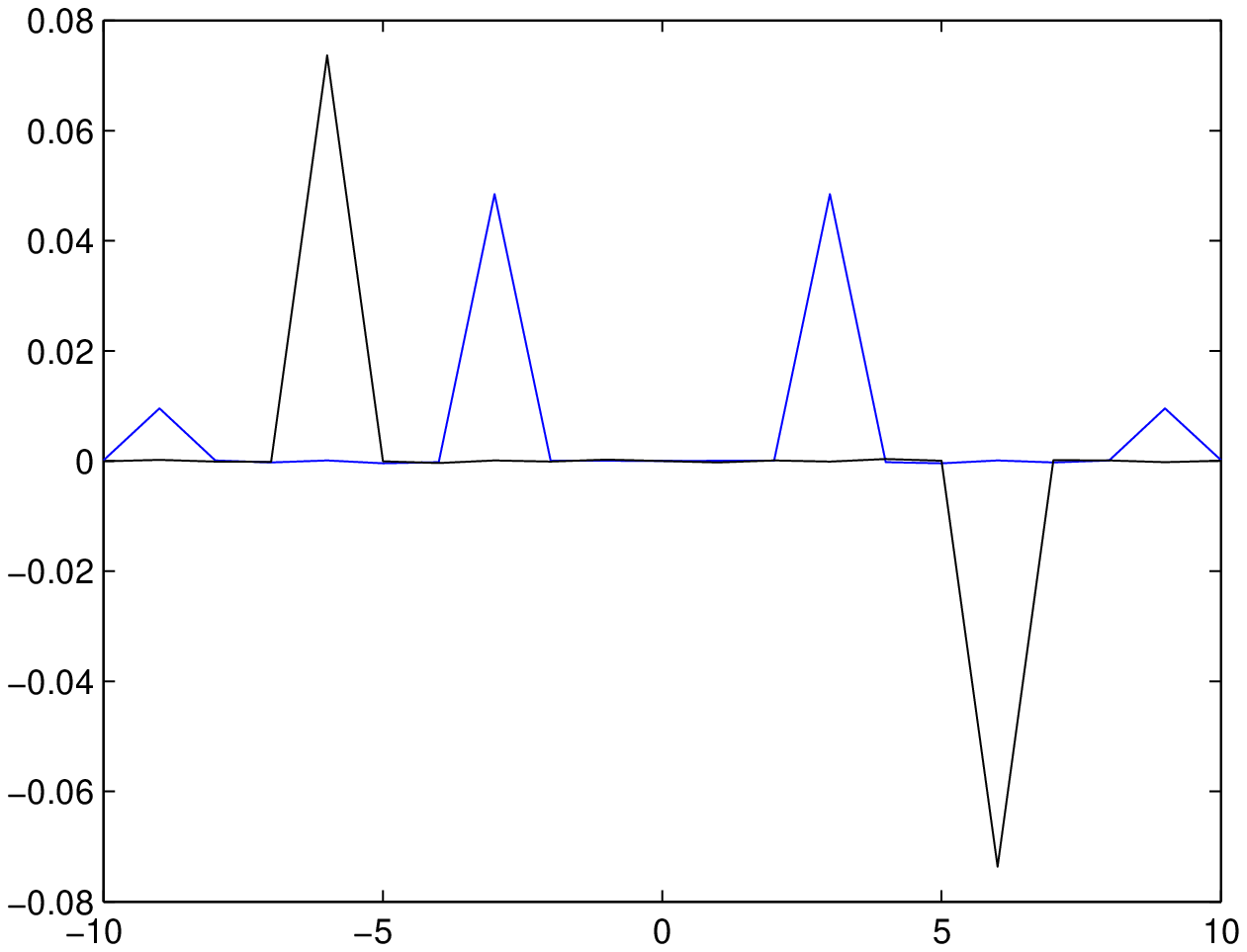}
\includegraphics[width=5cm,height=4cm]{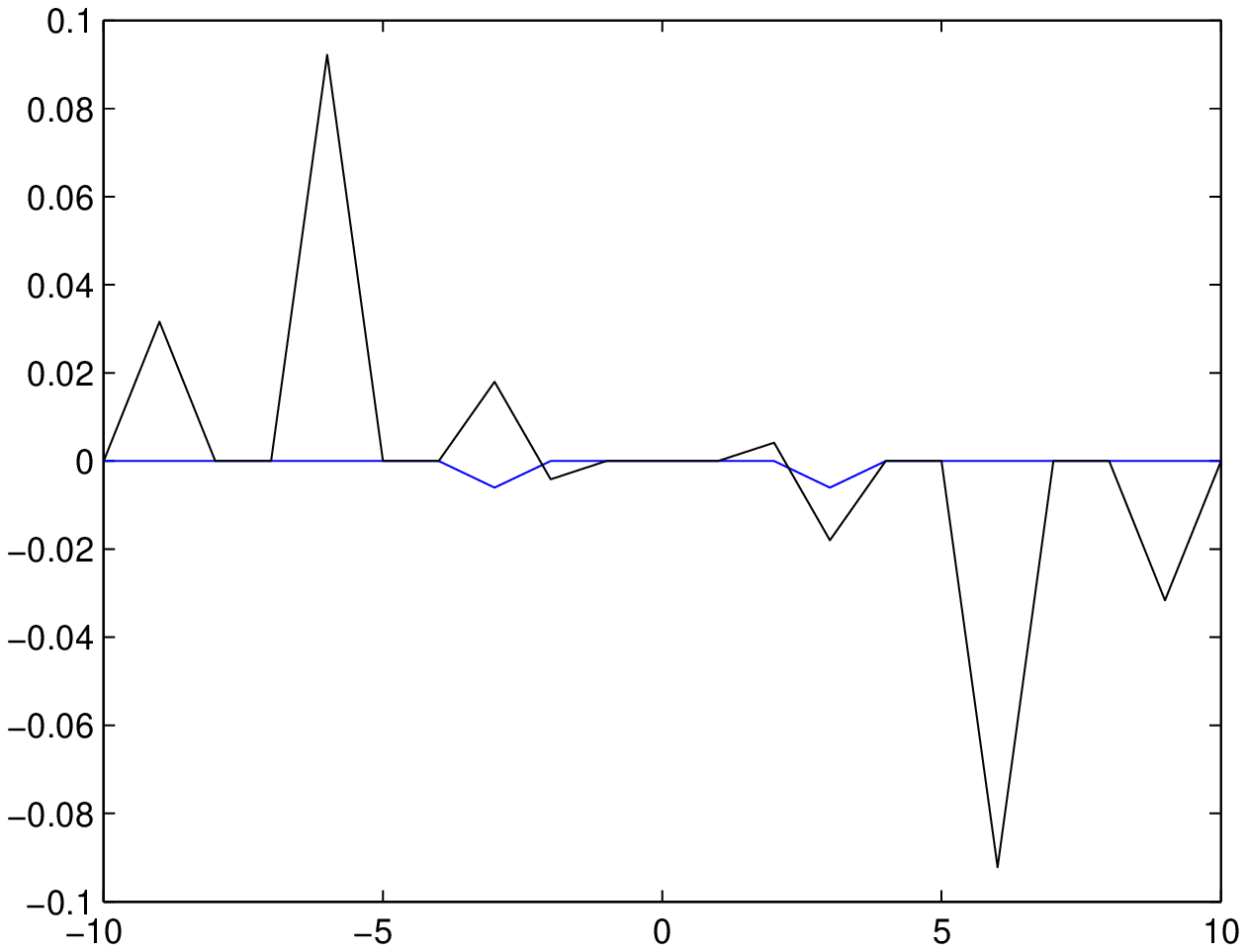}
\includegraphics[width=5cm,height=4cm]{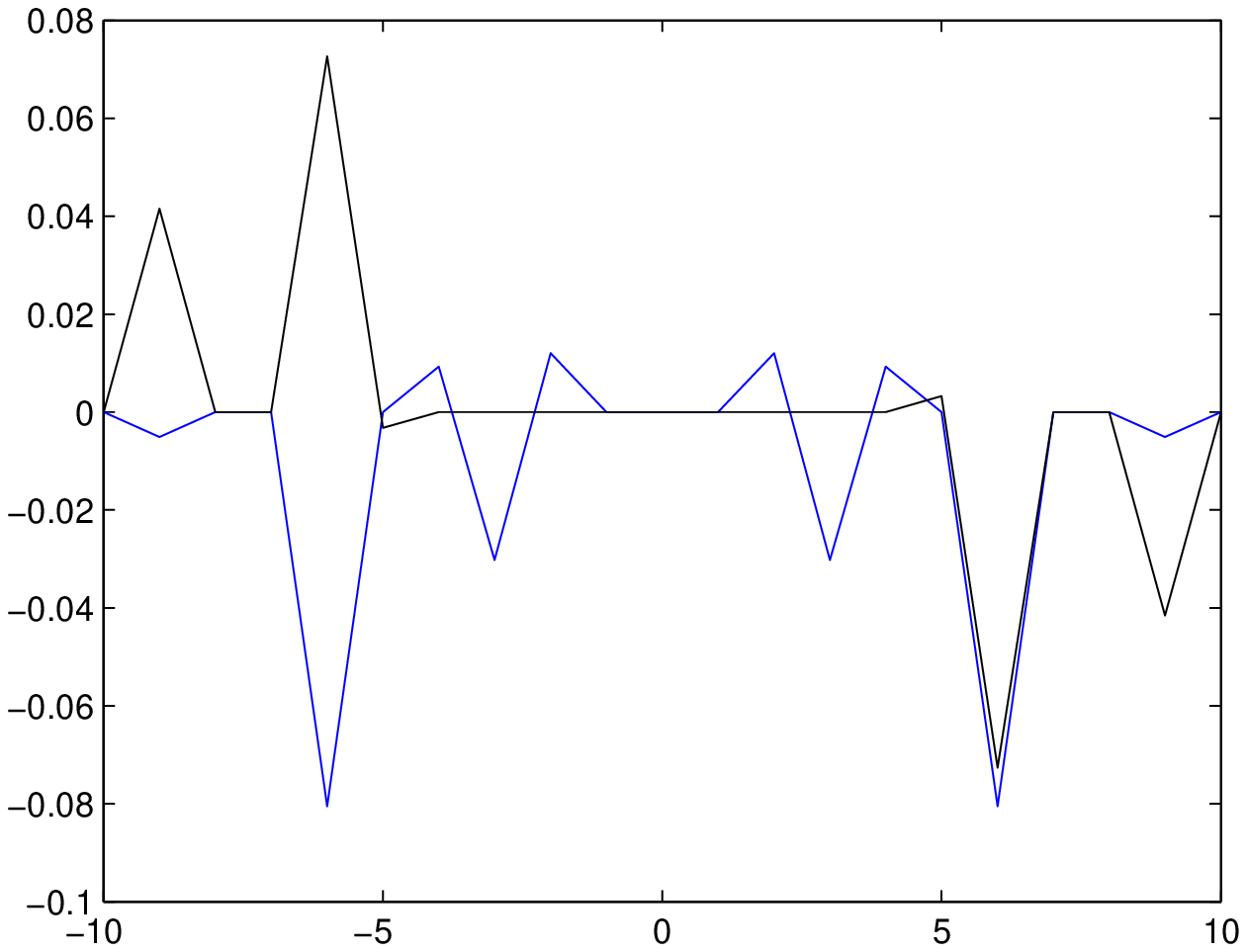}
\caption{Fourier coefficients of reconstructed perturbations in Example 3; \textbf{Set 1} to \textbf{Set 3} from left to right;   
blue: real part; black: imaginary part.} \label{test4B}
\end{figurehere}

In this example, as we can see, the reconstructions from \textbf{Set 1} is the best, with both the peak Fourier modes and their magnitudes quite close to the exact one, although with some phase shifts.
Reconstruction from \textbf{Set 2} is still reasonable.  The magnitude of the $6$-th Fourier modes is closer to the exact one, however that of the $3$-th mode deviates further from the exact one, and they have more phase shifts.
Reconstruction from \textbf{Set 3} is the worst, with great deficiency from the exact perturbation, considering the fact this reconstruction gives us many modes that do not exist in the exact perturbation.

In Figure \ref{test4C} (top), (middle) and (right), 
the shapes of reconstructed domains, 
the contrast of the reconstructed media and 
a comparison between the reconstructed domains $D^{\text{approx}}$ and exact domain $D$
are presented respectively
The relative $L^2$ errors of the reconstructions for \textbf{Set 1} to 
\textbf{Set 3} are respectively $11.61 \%$,
$  13.48 \%$ and
$  14.84 \%$. This indicates that the reconstruction for \textbf{Set 1} is the best, that for 
\textbf{Set 2} is still good, and that for \textbf{Set 3} is the worst. This goes with the theory we discussed in section \ref{sec71}.

\begin{figurehere} \centering
\includegraphics[width=5cm,height=4cm]{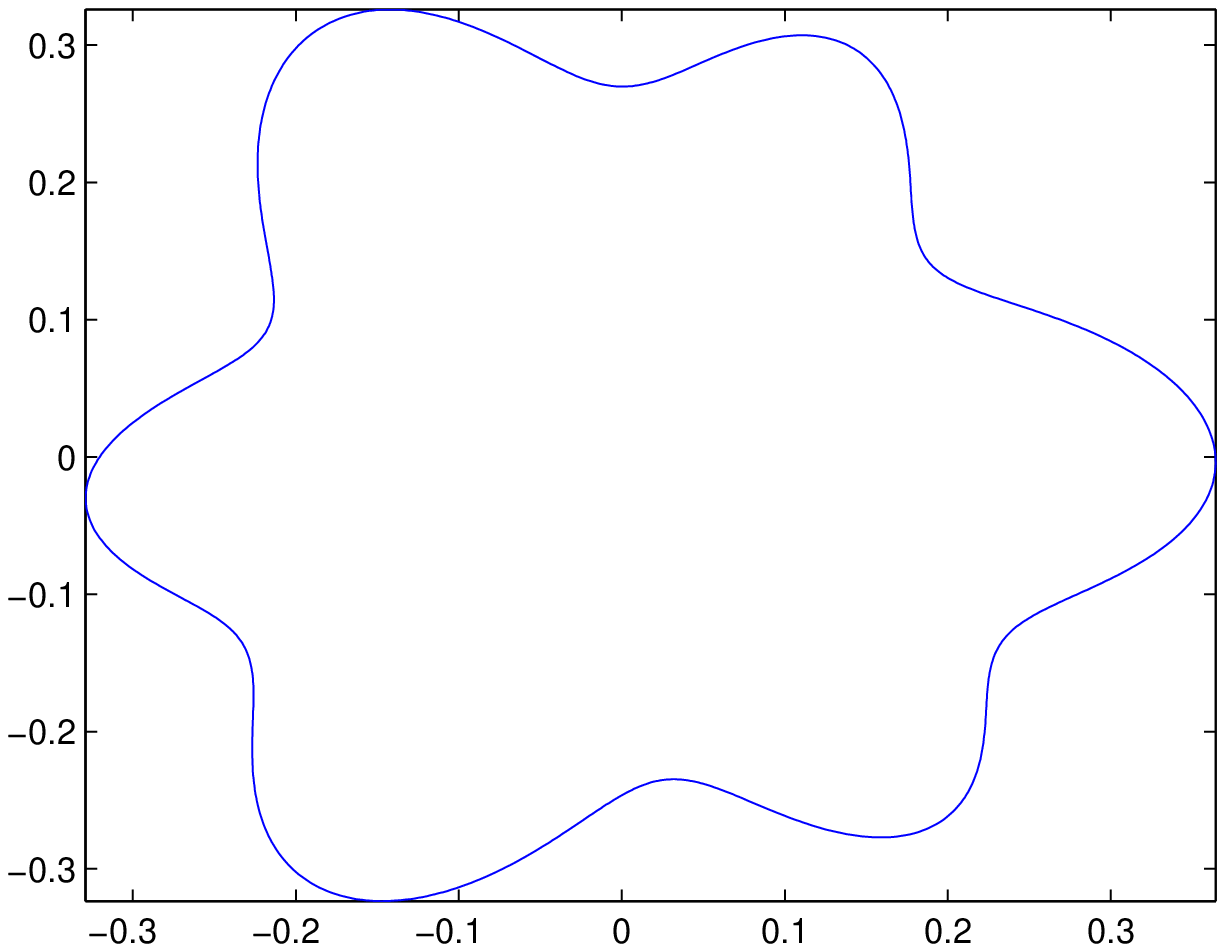}
\includegraphics[width=5cm,height=4cm]{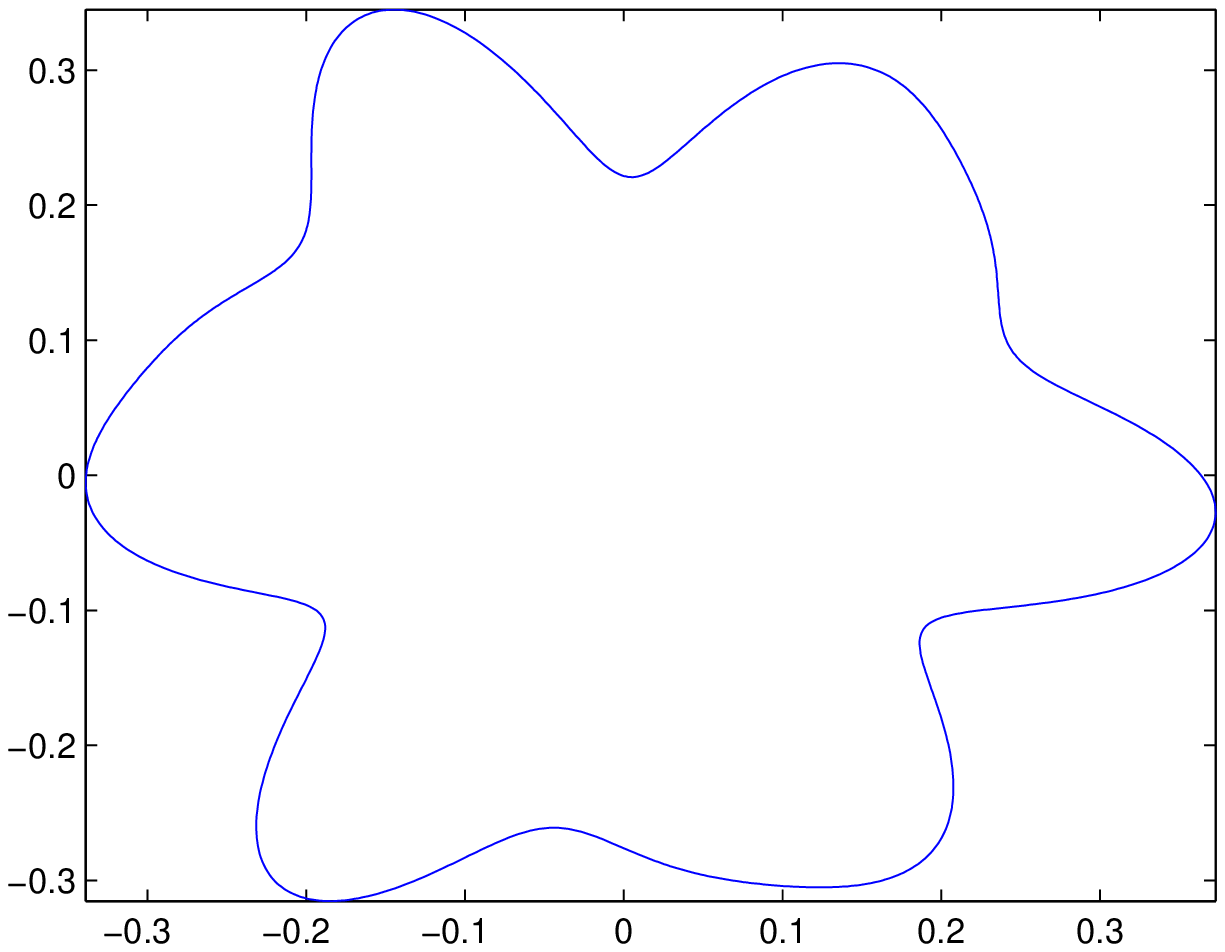}
\includegraphics[width=5cm,height=4cm]{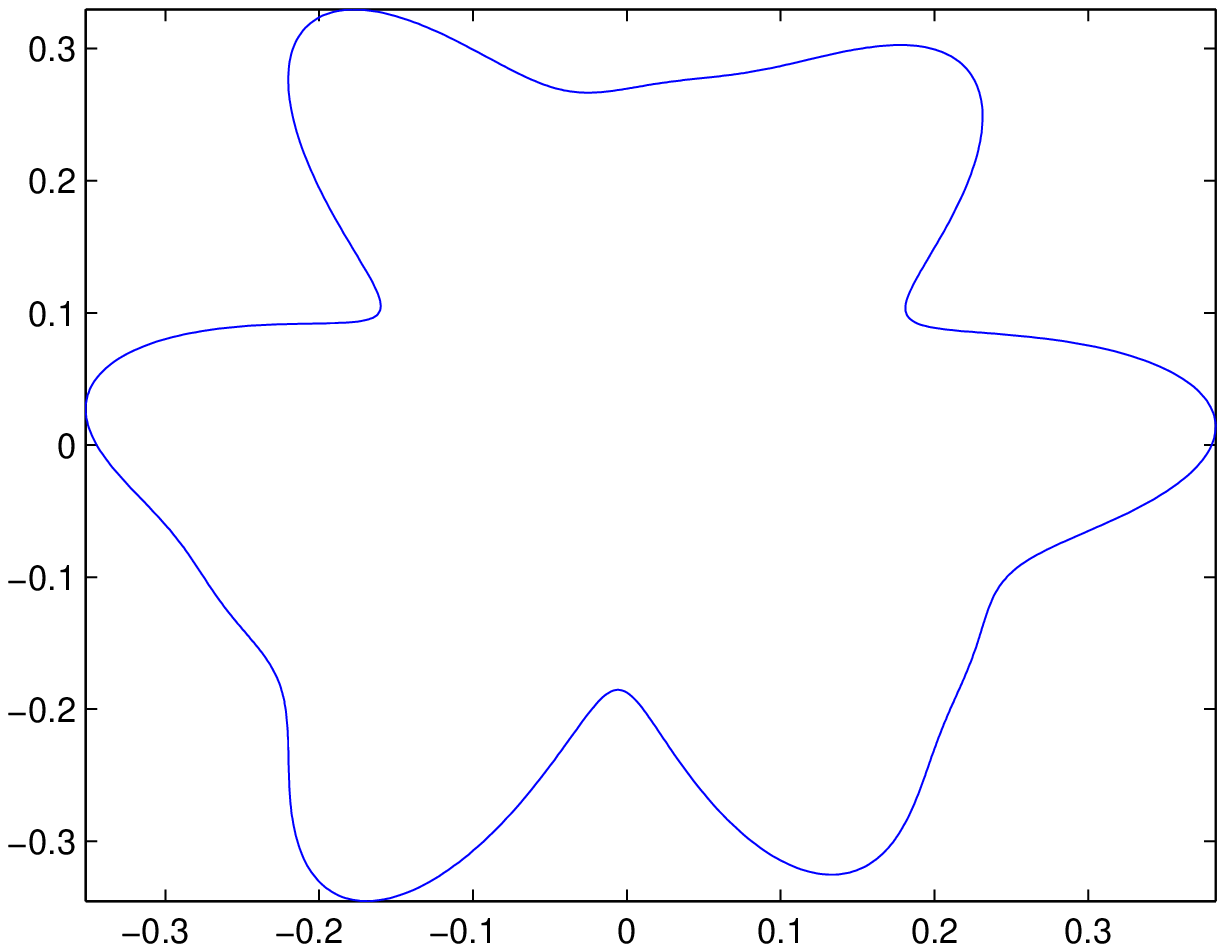} \\
\includegraphics[width=5cm,height=4cm]{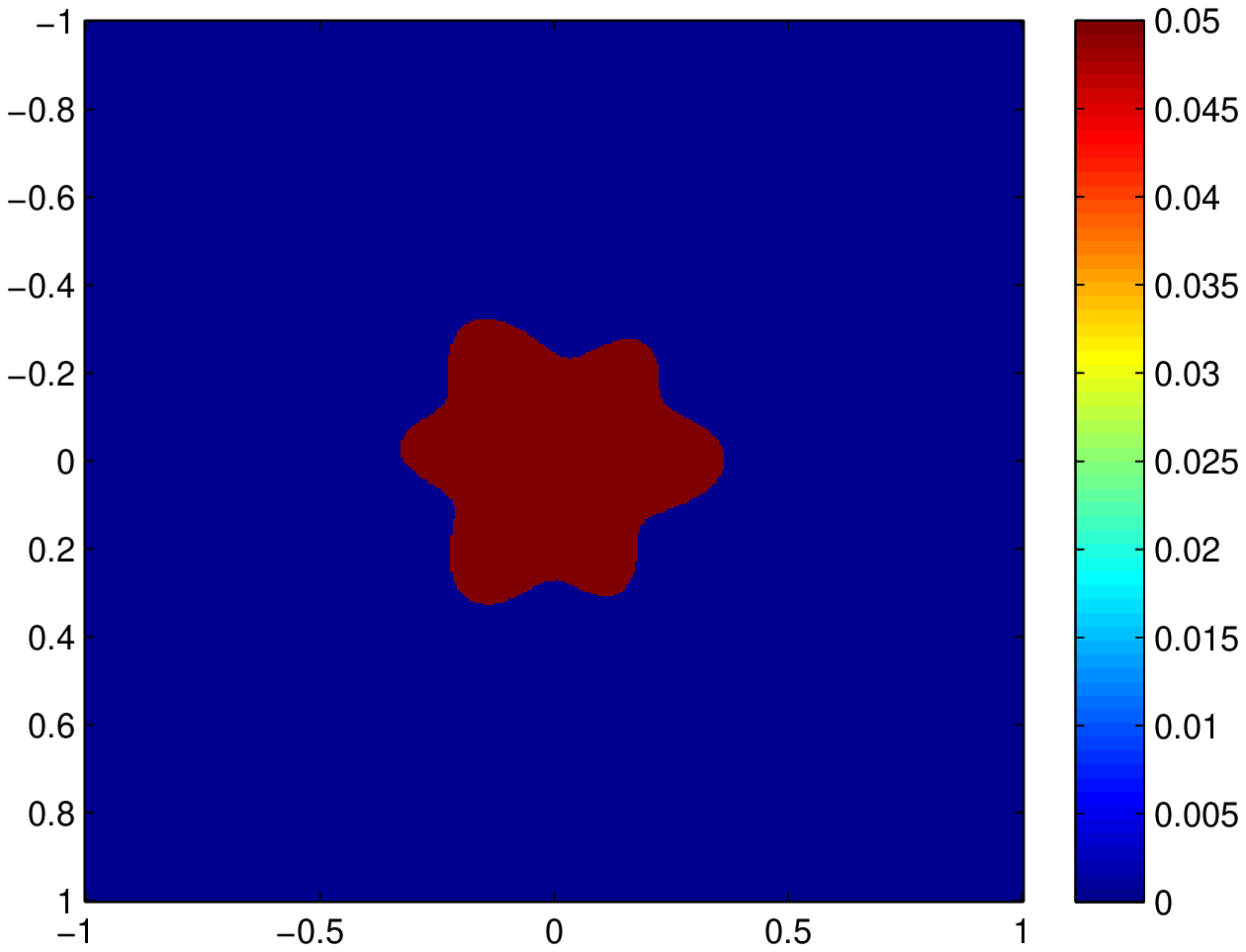}
\includegraphics[width=5cm,height=4cm]{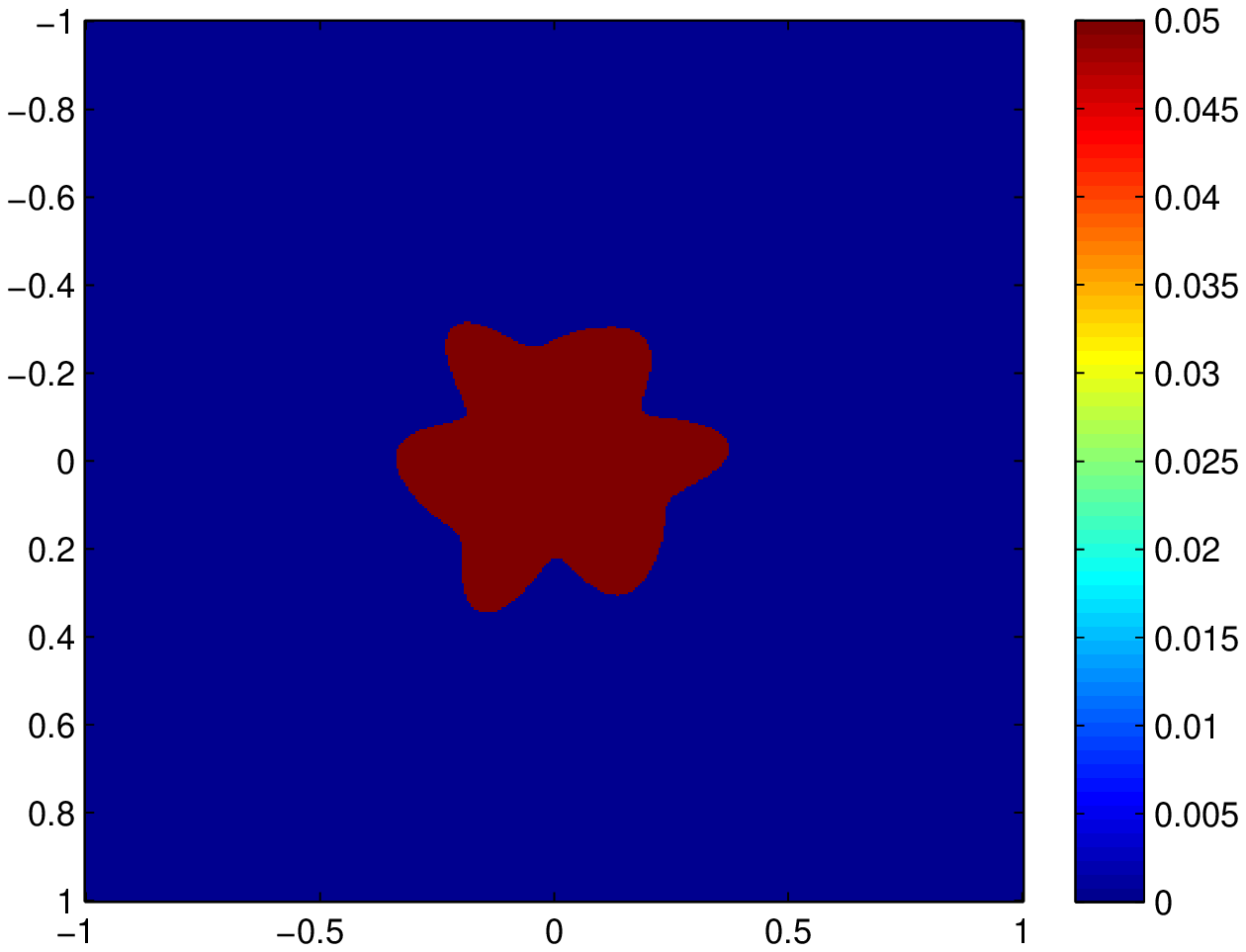}
\includegraphics[width=5cm,height=4cm]{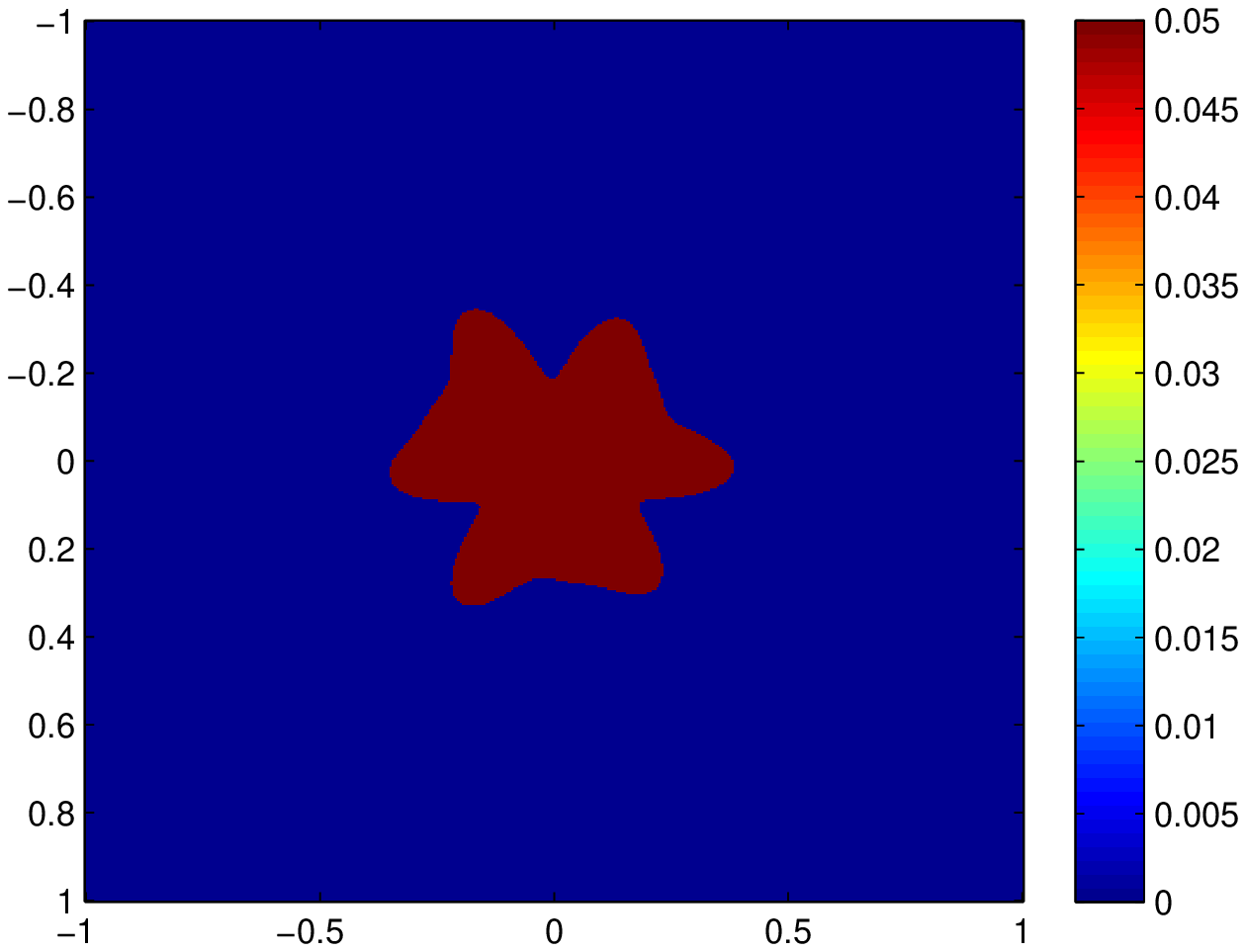} \\
\includegraphics[width=5cm,height=4cm]{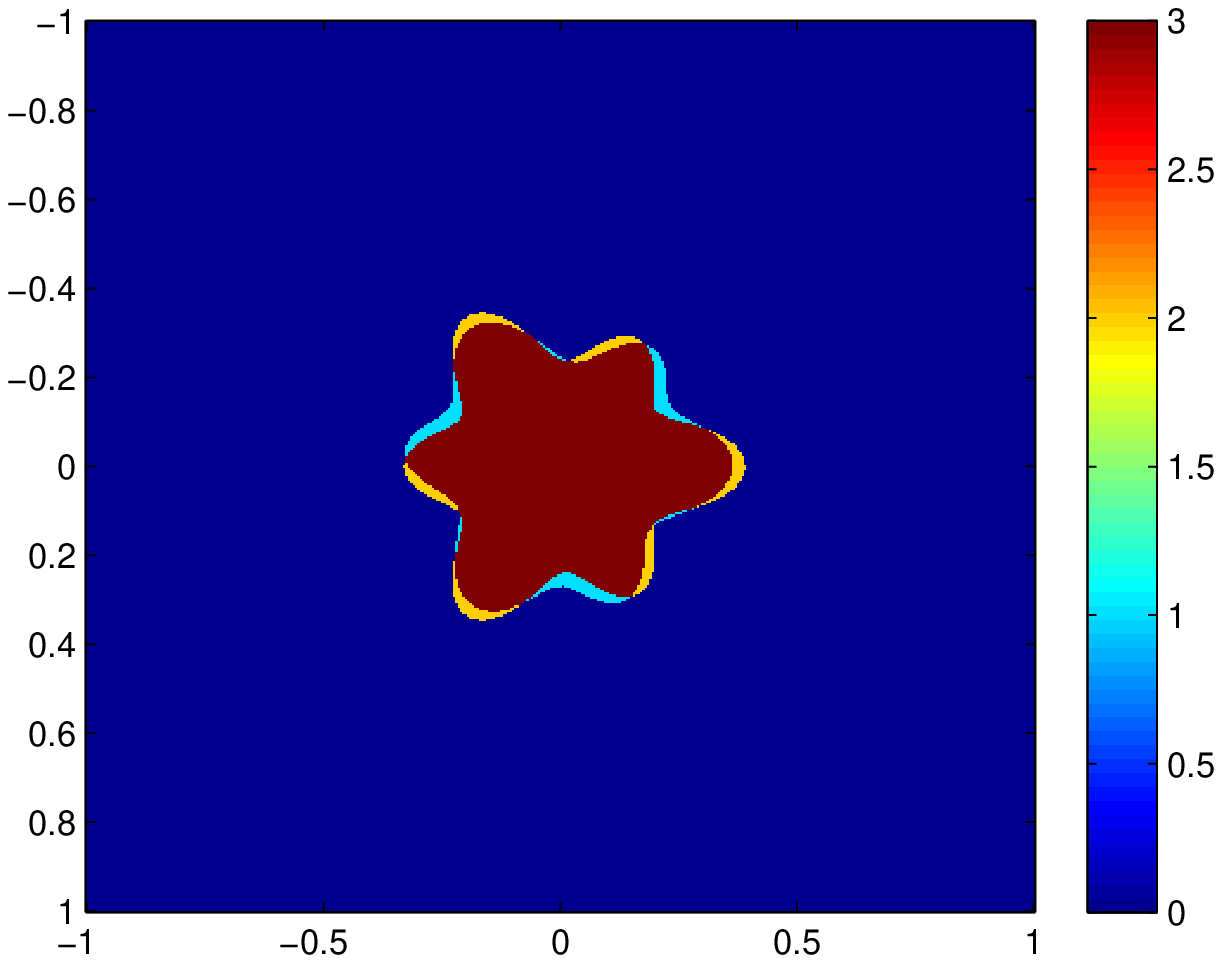}
\includegraphics[width=5cm,height=4cm]{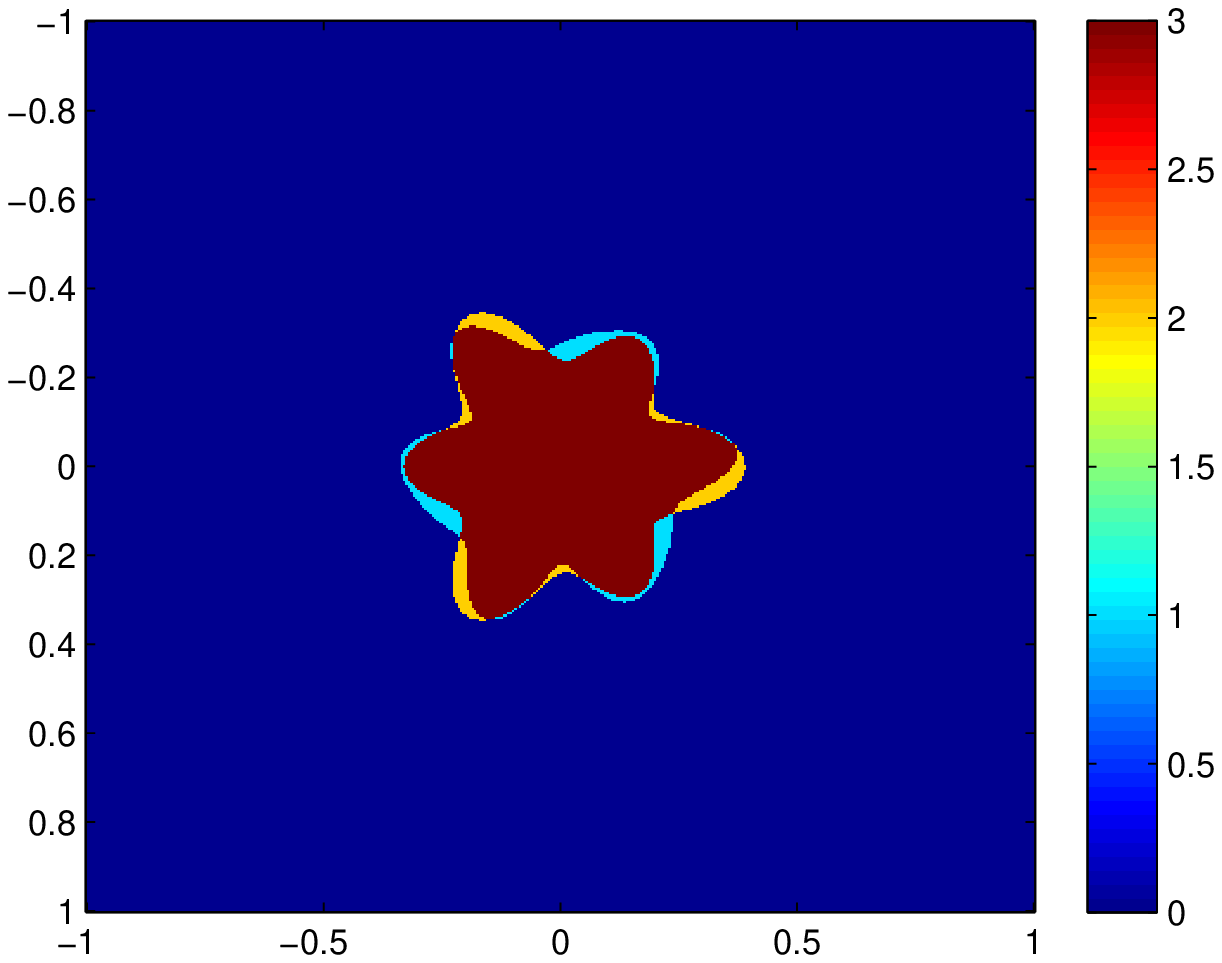}
\includegraphics[width=5cm,height=4cm]{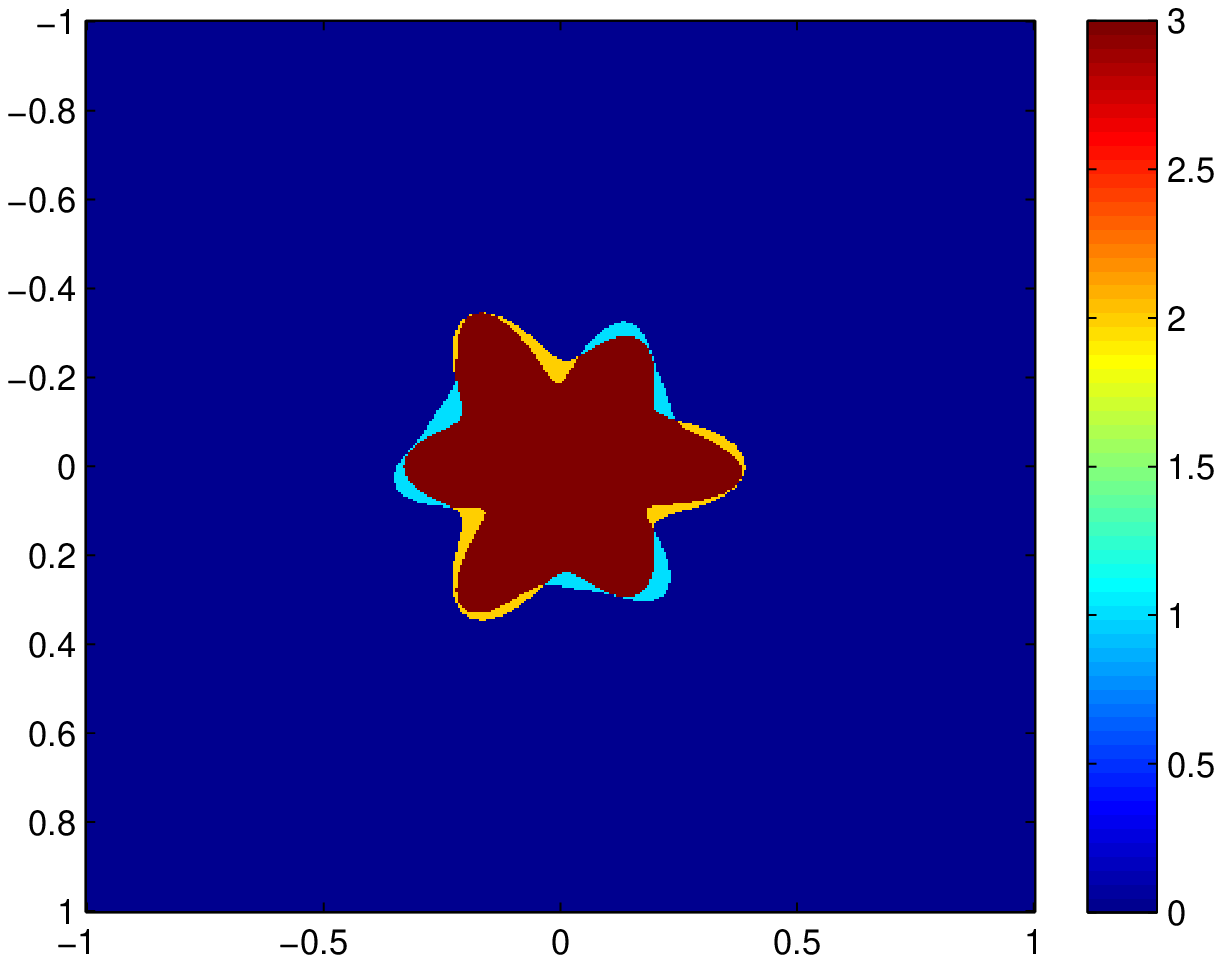} \\
\caption{Reconstructed domain and medium in Example 3 and comparison between the exact and reconstructed domains;   \textbf{Set 1} to \textbf{Set 3} from left to right; reconstructed shape, reconstructed inclusion and 
comparison between reconstructed and exact domains from top to bottom.}  \label{test4C}
\end{figurehere}

The reconstructions for \textbf{Set 1} (over-abundant number of measurements) and \textbf{Set 2} (critical number of measurements) are quite reasonable, considering the severe ill-posedness of the phaseless reconstruction problem 
and a $5\%$ percent measurement noise.

\end{document}